\documentclass[journal]{IEEEtran}
\usepackage[pdftex,dvipsnames]{xcolor}
\usepackage{tikz}
\usepackage{tikz}
\usepackage{array}
\newcolumntype{P}[1]{>{\centering\arraybackslash}p{#1}} 
\ifCLASSOPTIONcompsoc
\usepackage[caption=false,font=normalsize,labelfont=sf,textfont=sf]{subfig}
\else
\usepackage[caption=false,font=footnotesize]{subfig}
\fi
\usepackage{textcomp}
\usepackage{stfloats}
\usepackage{verbatim}
\usepackage{dirtytalk}
\usepackage{multirow}
\usepackage{rotating}
\usepackage{booktabs}
\usepackage{tikz}
\usepackage{epstopdf}
\usepackage{threeparttable} 
\usepackage{transparent}
\usepackage{graphicx}
\usepackage{amsmath,amsfonts}
\usepackage{bm}
\usepackage{hyperref}
\usepackage{url}
\usepackage{cite}
\usepackage{subfig}
\usepackage{algorithm} 
\usepackage{algpseudocode} 
\usepackage{tikz-network}
\usetikzlibrary{arrows,automata,positioning}

\usepackage{tikz}
\usetikzlibrary{trees}
\usetikzlibrary{chains, fit}
\usepackage{listings}
\usepackage{xcolor,colortbl}
\usetikzlibrary{shapes,arrows}
\usetikzlibrary{shadings}
\usetikzlibrary{automata,topaths}
\usetikzlibrary{shapes}
\usetikzlibrary{shapes.geometric, arrows}
\usepackage{tikz-3dplot}
\usetikzlibrary{matrix}
\usepackage{comment}

\usepackage{pgfplots}
\usetikzlibrary{pgfplots.groupplots}

\newcommand{\R}{{\rm I\!R}}
 
\newcommand{\addbb}[1]{{{\color{blue!0!black}#1}}}
\newcommand{\addb}[1]{{{\color{blue!0!black}#1}}}

\usetikzlibrary{quotes}

\def\BibTeX{{\rm B\kern-.05em{\sc i\kern-.025em b}\kern-.08em
    T\kern-.1667em\lower.7ex\hbox{E}\kern-.125emX}}
\usepackage{balance}
  
\newcommand{\boundellipse}[3]
{(#1) ellipse (#2 and #3)
}

\newcommand{\cf}{{\emph{cf.~}}}
\newcommand{\tran}{^{\mbox{\scriptsize T}}}
\usepackage{pifont}
\newtheorem{theorem}{Theorem}
\newtheorem{lemma}[theorem]{Lemma}
\newtheorem{assump}{Assumption}

\numberwithin{assump2}{section} 

\newtheorem{remark}{Remark}

\newtheorem{prop}{Proposition}

\numberwithin{prop2}{section} 

\newtheorem{definition}{Definition}
\newtheorem{corollary}{Corollary}
\newtheorem{example}{Example}

\makeatletter
\newcommand{\nosemic}{\renewcommand{\@endalgocfline}{\relax}}
\newcommand{\dosemic}{\renewcommand{\@endalgocfline}{\algocf@endline}}
\let\oldnl\nl
\newcommand{\nonl}{\renewcommand{\nl}{\let\nl\oldnl}}
\makeatother
\RequirePackage{etex}
\usepackage{amssymb}
\newcommand{\QE}{\hfill $\blacksquare$}

\begin{document}
\title{On the Characteristics of the Conjugate Function Enabling Effective Dual Decomposition Methods}
 \author{Hansi~Abeynanda,
        Chathuranga~Weeraddana, and~Carlo~Fischione,~\IEEEmembership{Fellow,~IEEE}
\thanks{H. Abeynanda and C. Fischione are with the School of Electrical Engineering and Computer Science, KTH Royal Institute of Technology, Sweden (e-mail: hkab@kth.se and carlofi@kth.se, respectively).}
\thanks{C. Weeraddana is with the Centre for Wireless Communication,
University of Oulu, Finland (e-mail: chathuranga.weeraddana@oulu.fi).}
}


\maketitle

\begin{abstract}
We investigate a novel characteristic of the conjugate function associated to a generic convex optimization problem, which can subsequently be leveraged for efficient dual decomposition methods. In particular, under mild assumptions, we show that there is a specific region in the domain of the conjugate function such that for any point in the region, there is always a ray originating from that point along which the gradients of the conjugate remain constant. We refer to this characteristic as a \emph{fixed gradient over rays} (FGOR). We further show that this characteristic is inherited by the corresponding dual function. Then we provide a thorough exposition of the application of the FGOR characteristic to dual subgradient methods. More importantly, we leverage FGOR to devise a simple stepsize rule that can be prepended with state-of-the-art stepsize methods enabling them to be more efficient. Furthermore, we investigate how the FGOR characteristic is used when solving the global consensus problem\addb{, a prevalent formulation in diverse application domains}. We show that FGOR can be exploited not only to expedite the convergence of the dual decomposition methods but also to reduce the communication overhead. \addbb{FGOR is extended to nonconvex formulations, and its advantages in stochastic optimization are demonstrated.} Numerical experiments using quadratic objectives and a regularized least squares regression with \addbb{real datasets} are conducted. The results show that FGOR can significantly improve the performance of existing stepsize methods \addb{and outperform the state-of-the-art splitting methods on average in terms of both convergence behavior and communication efficiency. }


\end{abstract}




\section{Introduction}\label{sec:Introduction}

\IEEEPARstart{D}{istributed} algorithms for optimization problems have become necessary and pervasive in various application domains such as signal processing, machine learning, wireless communications, control systems, robotics, and many others~\cite{Cao-2012-overview,di_2020,Nedic-Distributed-Optimization-for-Control-2018,Yang-survey-of-distributed-optimization-2019,Palomar-book-2010,Carlo-Machine-Learning,condat_2013_splitting_methods,Boyd-Parikh-Chu-Peleato-Eckstein-2010,He_2021,Bai_2024}. These problems are typically of very large scale since they deal with thousands, millions, or even more variables. In this respect, effective deployment of such algorithms requires an appeal to light communication among subsystems involved in solving the optimization problem and less computational effort per subsystem, and those inevitably raise questions about how fast the deployed algorithms converge. \addbb{Although second-order algorithms such as Newton’s methods \cite{Hintermuller_2010}, \cite[\S~10]{Boyd-Vandenberghe-04} exhibit fast convergence properties, they are generally less appealing for distributed large-scale optimization due to their substantial computational overhead (e.g., the need for large matrix inversions). As such, in many real-world applications, the most commonly employed method to solve optimization problems is the (sub) gradient method which is a first-order algorithm.} The choice of (sub) gradient method sparked numerous efforts by many researchers in designing various stepsize rules and applying them to yield faster convergences~\cite{Polyak_Intro_Opt,goldstein_1962,armijo_1966,beck_2009,bello_2016,zhang2004nonmonotone,Asl_2020,khirirat_2023,Thinh_2018,Thinh_2021,zhor_1968_1,zhor_1968_2,goffin_1977,shor_2012,Damek_2018,zhu_2019,aybat_2019,wang_2021,chen_2019,krizhevsky_2012,He_2016,ge_2019,polyak_1969,hazan_2022,loizou_2021,wang_2023,barzilai_1988,dai_2002,burdakov_2019,Yura_2020}.

The existing literature that designs adaptive stepsizes primarily considers subgradient methods within the primal domain. These results can be readily applied to distributed algorithms that are based on duality~\cite{Sindri_2020,R1-Liu-2021,Necoara-inexact-dual-decomposition-2014,Yifan-Asynchrony-and-Inexactness-2023}. More specifically, in the case of dual decomposition algorithms for distributed optimization, the master problem that coordinates the subproblems at subsystems is often solved by using subgradient methods~\cite{Boyd-EE364b-PrimDualDecomp-07} where the earlier results for stepsizes can be adopted. However, these stepsize methods are applied to the dual domain without exploiting possible useful characteristics of such domain. Thus, we now pose the natural question: can we exploit any specific characteristics of the dual function itself to further scrutinize the design of stepsize rules that go beyond those designs within the primal domain for effective convergence of the subgradient method associated to the master problem? 

To answer such a question, let us formalize the setting and start with the following generic primal problem: 
\begin{equation} \label{eq:optimization_prob_main}
\begin{array}{ll}
\mbox{minimize} & f_0(\mathbf{y})\\
\mbox{subject to} & \mathbf{y}\in\mathcal{Y}\\
\mbox{} & \mathbf{A}\mathbf{y} = \mathbf{b},
\end{array}
\end{equation}
where the variable is $\mathbf{y}\in\R^n$, $f_0:\R^{n}\rightarrow \R $ is the convex objective function, $\mathcal{Y}$ is a convex set in $\R^n$, and $\mathbf{b}\in\R^{m}$. \addb{Without loss of generality, we assume $\mathbf{A}\in\R^{m\times n}$ with $m<n$ and \texttt{rank} $\mathbf{A}=m$}. Equivalently, one may reformulate problem~\eqref{eq:optimization_prob_main} as 
\begin{equation} \label{eq:optimization_prob_main_extended}
\begin{array}{ll}
\mbox{minimize} & f(\mathbf{y})=f_0(\mathbf{y}) + \delta_\mathcal{Y}(\mathbf{y})\\
\mbox{subject to} & \mathbf{A}\mathbf{y} = \mathbf{b},
\end{array}
\end{equation}
where $\delta_{\mathcal{Y}}$ denotes the indicator function of the set $\mathcal{Y}$, \cf~Definition~\ref{Definition-Indicator-function}. Note that in a decomposition setting, $f_0$ and constraints are endowed with structural properties such as separability, enabling the application of the dual decomposition~\cite{Boyd-EE364b-PrimDualDecomp-07}.

The convergence of dual decomposition methods applied to problem~\eqref{eq:optimization_prob_main_extended}  hinges predominantly on the characteristics of the problem's associated dual function. Note that the dual function of~\eqref{eq:optimization_prob_main_extended} with $\mathbf{A}\mathbf{y} = \mathbf{b}$ being the constraint to which Lagrange multipliers are associated, is intimately connected with the conjugate function of $f$ \cite[\S~5.1.6]{Boyd-Vandenberghe-04}. As such, to answer our posed problem, we investigate a novel characteristic of the conjugate function of~$f$. In particular, we show that the variational geometry of~$\mathcal{Y}$, \cf~\cite[Fig.~6-8]{Rockafellar-98} and nondifferentiable properties of $f_0$ induce a specific characteristic on the conjugate function. We effectively use this characteristic to devise a precursory yet simple stepsize rule that can be prepended with other stepsizes~\cite{Polyak_Intro_Opt,goldstein_1962,armijo_1966,beck_2009,bello_2016,zhang2004nonmonotone,Asl_2020,khirirat_2023,Thinh_2018,Thinh_2021,zhor_1968_1,zhor_1968_2,goffin_1977,shor_2012,Damek_2018,zhu_2019,aybat_2019,wang_2021,chen_2019,krizhevsky_2012,He_2016,ge_2019,polyak_1969,hazan_2022,loizou_2021,wang_2023,barzilai_1988,dai_2002,burdakov_2019,Yura_2020} while improving their convergences. The proposed precursory stepsize rule enables not only distributed but also efficient implementation of the algorithms with a light communication overhead. \addbb{Our exposition begins with convex problems, as they provide a fundamental basis for developing insights that enable extensions to more general settings (e.g., nonconvex problems and stochastic methods). Furthermore, in the machine learning domain, these convex formulations correspond to simple and computationally lightweight models, such as linear/logistic regression and support vector machines, and are also appealing from an explainability and interpretability standpoint, properties that are particularly important in sensitive application domains such as medicine and bioinformatics~~\cite{karim2023explainable,karim2022interpreting}.
}

\subsection{Related Work}\label{sec:Related-Work}

Employing dual decomposition techniques to solve a primal problem of the form~\eqref{eq:optimization_prob_main_extended} with adequate structural properties solely relies on solving the associated dual problem~\cite[\S~2]{Boyd-EE364b-PrimDualDecomp-07}. Recall that such a dual problem is often solved by using subgradient methods. The traditional stepsize methods~\cite{goldstein_1962,armijo_1966,beck_2009,bello_2016,Asl_2020,khirirat_2023,Thinh_2018,Thinh_2021,zhor_1968_1,zhor_1968_2,goffin_1977,shor_2012,Polyak_Intro_Opt,Damek_2018,zhu_2019}, that have already been designed for subgradient methods are constant stepsize, line search, polynomially decay stepsizes, and geometrically decay stepsizes. These stepsizes directly apply to solving the dual problem, provided the dual function conforms to necessary regularity conditions. We give a short overview of these schemes in the following.

Given the dual function satisfies some gradient Lipschitz conditions and the Lipschitz constant is known, then a \emph{constant stepsize} can be chosen appropriately to ensure convergence of the algorithms~\cite[\S~1.4]{Polyak_Intro_Opt}. The concept of line search dates back to the 1960s \cite{goldstein_1962,armijo_1966} and is used in numerous application domains~\cite{beck_2009,bello_2016,zhang2004nonmonotone,Asl_2020}. It is a commonly used method for determining an appropriate stepsize, especially when the gradient Lipschitz constants are unknown. \addbb{However, unlike constant stepsize approaches, incorporating line search techniques does not lend itself to an efficient distributed implementation of dual decomposition methods. This limitation arises because, in each iteration of the dual decomposition algorithm, additional Lagrangian minimization steps must be coordinated among subsystems to compute the corresponding dual function and gradient values, which are subsequently used to verify the line search conditions. Furthermore, implementing even a single iteration of the line search procedure introduces significant signaling overhead among subsystems~\cite[p.~464]{Boyd-Vandenberghe-04},~\cite[p.~1044]{zhang2004nonmonotone}.}

The \emph{polynomially decay stepsizes} such as $\gamma_0/k$ and $\gamma_0/\sqrt{k}$, where $\gamma_0>0$ is an appropriately chosen constant \cite[\S~5.3]{Polyak_Intro_Opt}, \cite{khirirat_2023,Thinh_2018,Thinh_2021} are commonplace and highly popular in distributed optimization. The reason is that they can be readily applied to respective distributed algorithms because of their simple dependency on the iteration index $k$.  
Nonetheless, their drawback lies in the fact that, as they progressively decrease with each iteration, they tend to yield slow convergence~\cite{Bazaraa_1981,Bianchi_2022}.

The \emph{geometrically decay stepsize} strategy for subgradient methods originated from \cite{zhor_1968_1,zhor_1968_2,goffin_1977}. A comprehensive exposition is found in \cite[\S~2.3]{shor_2012}. The resulting linear rate of convergences associated to such stepsizes has gained attention in recent studies \cite{Damek_2018,zhu_2019}. However, linear convergence is guaranteed under stringent conditions restricting the use of such stepsizes in practice~\cite[Theorem~2.7]{shor_2012}.

The stepsize rules considered in first-order methods~\cite{aybat_2019,wang_2021,chen_2019,krizhevsky_2012,He_2016,ge_2019,polyak_1969,hazan_2022,loizou_2021,wang_2023,barzilai_1988,dai_2002,burdakov_2019,Yura_2020} are different from the direct application of any traditional stepsizes discussed above. More importantly, these stepsize strategies are designed to further improve the convergences, capture more general problem formulations, and handle problems when the gradient Lipschitz constants are unknown.
The rules therein can be classified as step decay stepsize~\cite{aybat_2019,wang_2021,chen_2019,krizhevsky_2012,He_2016,ge_2019} and running attributes-dependent stepsizes~\cite{Yura_2020,polyak_1969,hazan_2022,loizou_2021,wang_2023,barzilai_1988,dai_2002,burdakov_2019}. We overview them in the following.

The \emph{step decay stepsize} is widely adopted in stochastic non-convex optimization \cite{aybat_2019,wang_2021,chen_2019} and is used as a potential model for training deep neural networks \cite{krizhevsky_2012,He_2016}. This stepsize rule is characterized by maintaining constant stepsizes within stages of the algorithm iterations while decreasing it at each subsequent stage. The basic idea is to start stepping more aggressively in the initial stages and gradually become less aggressive in subsequent stages. The rule intrinsically includes two decisions: 1) how the stepsize is reduced and 2) how the stage length is changed, over the stage count. Typically, the rule reduces the stepsize by a fixed factor at each stage, i.e., the stepsize is $\alpha^t$, where $\alpha\in(0,1)$ and $t$ is the stage count. This policy, where the stepsize is reduced by a factor after a predetermined number of epochs, instead of at every iteration, is viewed as a variant of the geometrically decay stepsize. In contrast, the stepsize schedule employed in \cite{chen_2019} is proportional to $1/t$, where $t$ is the stage count. On the other hand, stage lengths are changed either with linear growth~\cite{chen_2019} or exponential growth~\cite{aybat_2019} or sometimes are kept constant~\cite{wang_2021,ge_2019}. In general, almost all step decay stepsize rules can be employed for dual decomposition methods.

The running attributes-dependent stepsizes that rely on algorithms' runtime data, such as decision variables, gradients, and function values, are extensively studied in~\cite{polyak_1969,hazan_2022,loizou_2021,wang_2023,barzilai_1988,dai_2002,burdakov_2019,Yura_2020}. \emph{Polyak's stepsize} \cite{polyak_1969} is an appealing rule demonstrating linear convergence rates for general strongly convex functions, even under non-differentiable settings. However, it requires the optimal value of the underlying problem or initial guess of it with subsequent refinement to compute the stepsize in every iteration. Despite such restrictions, recent studies have revived interest in Polyak's stepsize due to its linear convergence properties, \cf~\cite{hazan_2022,loizou_2021,wang_2023}.

Another running attributes-dependent stepsize is the \emph{Barzilai-Borwein} (BB) \emph{rule} \cite{barzilai_1988,dai_2002,burdakov_2019}. Unlike Polyak's stepsize, the BB rule doesn't require the knowledge of the optimal value of the underlying problem. Initial BB results are limited to two-dimensional quadratic objective functions with linear convergence guarantees~\cite{barzilai_1988}. Linear convergence with BB stepsizes under strictly convex quadratic functions, in general, has been established in \cite{dai_2002}. Extensions of similar convergence results under more general settings, in particular, for smooth and strongly convex functions have been established in \cite{burdakov_2019}, however under more restricted assumptions. In \cite{Yura_2020}, a more appealing stepsize is proposed without the need for any restricted assumptions like those in Polyak's stepsize and BB rule. More specifically, the stepsize proposed in \cite{Yura_2020} requires no specific assumptions beyond the local gradient Lipschitz continuity. In a dual decomposition setting, all the above stepsize rules~~\cite{polyak_1969,hazan_2022,loizou_2021,wang_2023,barzilai_1988,dai_2002,burdakov_2019,Yura_2020} can be deployed, given the necessary runtime data in each iteration is communicated to the entity that is updating the dual variables, e.g., a parameter server. 

When solving the dual problem in a dual decomposition setting, one can always employ first-order algorithms~\cite{Polyak_Intro_Opt,goldstein_1962,armijo_1966,beck_2009,bello_2016,Asl_2020,khirirat_2023,Thinh_2018,Thinh_2021,zhor_1968_1,zhor_1968_2,goffin_1977,shor_2012,Damek_2018,zhu_2019,aybat_2019,chen_2019,wang_2021,krizhevsky_2012,He_2016,ge_2019,polyak_1969,hazan_2022,loizou_2021,wang_2023,barzilai_1988,dai_2002,burdakov_2019,Yura_2020} that are devised for subgradient methods. Nevertheless, it is instructive to observe that the dual function of the underlying primal formulation, can be endowed with certain characteristics to design stepsize rules to further improve the convergences.   
In this respect, dualization of strong convexity~\cite[Prop.~12.60]{Rockafellar-98} which entangles the gradient Lipschitzian properties of primal functions and strong convexity of the associated dual functions and vice versa is one of the most classical results that enables the choice of constant stepsize length in a dual decomposition setting, \cf~\cite[\S~1.4]{Polyak_Intro_Opt}. But to the best of the authors' knowledge, there seems to be no other explicit record previously of the employment of dual characteristics towards designing stepsize rules for dual decomposition methods.
\subsection{Our Contribution}
We investigate a novel characteristic of the conjugate function associated to the convex-constrained optimization, which we can leverage to improve the performance of dual subgradient methods. In particular, within problem setting~\eqref{eq:optimization_prob_main}, the main contributions of this paper are given below. 

\begin{itemize}
    \item \emph{Fixed-gradient-over-rays} (FGOR) characteristic of the conjugate function: we show that there is a specific region in the domain of the conjugate function $f^*$ [\cf~Definition~\ref{Definition:Conjugate-function}] of $f$ in problem~\eqref{eq:optimization_prob_main_extended} such that for all points $\nu$ in the region there is a ray originating from $\nu$ along which the gradients of $f^*$ remain constant, \cf~\S~\ref{Sec:Flatness-Conjugate-Function}, Proposition~\ref{Prop:flat-region}. 
    \item FGOR characteristic of the dual function: following the FGOR characteristic of $f^*$, we establish that the domain of the dual function $g$ of problem~\eqref{eq:optimization_prob_main_extended} contains a ray along which the gradients of $g$ remain constant, \cf~\S~\ref{sec:Flatness-of-g}, Corollary~\ref{Prop:flat-region-dual-function}.
    \item Application of FGOR on the dual subgradient method: using established FGOR characteristics, we devise a simple stepsize rule  that can be prepended with state-of-the-art stepsize methods while improving the convergence of the dual subgradient method, \cf~\S~\ref{sec:Dual-Subgradient-Algorithm-The-Proposed- Approach}, Algorithm~\ref{Alg:Dual-Subgradient-Method}.
\item \addb{Application of FGOR to real-world problems}: we explore how FGOR can be exploited when solving the global consensus problem [\cf~\eqref{eq:consensus-problem-distributed}], \addb{one of the most prevalent formulations in large-scale signal processing and machine learning applications.} In particular, we show that by leveraging FGOR, the performance of the standard dual decomposition algorithm can be improved in terms of both speed of convergence and communication efficiency, \cf~\S~\ref{sec:Application:The-General-Consensus-Problem}, Algorithm~\ref{Alg:Dual-decompositoin-FGOR}, Lemma~\ref{prop:consensus-problem-feasible-y-bar}, and \S~\ref{sec:numerical-Results}.
\item \addbb{Extensions to nonconvex and stochastic settings: we extend the established FGOR properties to nonconvex formulations and stochasic methods, \cf~Appendix~\ref{Sec:Appendix:Extensions-nonconvex-case} and Appendix~\ref{subsec:FGOR-for-Stochastic-Subgradient-Method}.} %

    \item \addb{The numerical results highlight that FGOR can substantially improve the performance of existing stepsize methods~\cite{Polyak_Intro_Opt,goldstein_1962,armijo_1966,beck_2009,bello_2016,Asl_2020,khirirat_2023,Thinh_2018,Thinh_2021,zhor_1968_1,zhor_1968_2,goffin_1977,shor_2012,Damek_2018,zhu_2019,aybat_2019,chen_2019,wang_2021,krizhevsky_2012,He_2016,ge_2019,polyak_1969,hazan_2022,loizou_2021,wang_2023,barzilai_1988,dai_2002,burdakov_2019,Yura_2020} and, on average, outperform state-of-the-art splitting methods~\addbb{\cite{condat_2013_splitting_methods,Boyd-Parikh-Chu-Peleato-Eckstein-2010,He_2021,Bai_2024}} with respect to both convergence behavior and communication efficiency, \cf~\S~\ref{sec:numerical-Results}.}
    \item Lipschitzian properties of the dual function $g$: under less restricted assumptions, we establish that the gradients $\nabla g$ of~$g$ are Lipschitz continuous, a more general Lipschitzian property than the existing results, \cf~Appendix~\ref{Appendix:Imp-Assump-Conj-Func}, Proposition~\ref{Prop:Lipscitz-Continuity}, Corollary~\ref{Lemma:UnconstrainedMin-of-f}.
\end{itemize}

\subsection{Notation}

Normal font lowercase letters $x$, bold font lowercase letters~$\mathbf{x}$, bold font uppercase letters $\mathbf{X}$, and calligraphic font $\mathcal{X}$ represent scalars, vectors, matrices, and sets, respectively. The set of real numbers, set of extended real numbers, set of real $n$-vectors, set of real $m\times n$ matrices, set of positive integers, and set of nonnegative integers are denoted by $\R$, $\overline{\R}$, $\R^n$, $\R^{m\times n}$, $\mathbb{Z}_+$, and $\mathbb{Z}^0_+$, respectively. The boundary, the interior, the convex hull, and \addb{the closure of a set $\mathcal{X}$} are denoted by $\texttt{bnd}~\mathcal{X}$, $\texttt{int}~\mathcal{X}$, $\texttt{con}~\mathcal{X}$, and \addb{$\texttt{cl}~\mathcal{X}$}, respectively. The domain of a function $f:\R^n \to \R$ is a subset of $\R^n$ and is denoted by $\texttt{dom}~f$. The set of all the subgradients of a function $f$ at a point $\mathbf{x}\in\R^n$ is denoted by $\partial f(\mathbf{x})$. \addbb{The convex hull of a function $f:\R^n \to \R$ is denoted by $\texttt{con}~f$.} The range \addb{and rank} of a matrix~$\mathbf{X}$ are denoted by $R(\mathbf{X})$ \addb{and $\texttt{rank}~(\mathbf{X})$, respectively.} For a given matrix $\mathbf{X}$, $\mathbf{X}\tran$ denotes the matrix transpose. The positive definite cone is denoted by $\mathbb{S}^ {n}_{++}$. The $n\times n$ identity matrix and the $n$-vectors with all entries equal to one are denoted by $\mathbf{I}_n$ and $\mathbf{1}_n$, respectively. 

\subsection{Organization of the Paper} \label{Sec:Organization-of-the-Paper}
The rest of the paper is organized as follows. In \S~\ref{Sec:Flatness-Conjugate-Function}, we investigate the FGOR characteristic of the conjugate function. The FGOR characteristic of the dual function is explored in \S~\ref{sec:Flatness-of-g}. In \S~\ref{sec:Dual-Subgradient-Algorithm-The-Proposed- Approach}, the application of FGOR on the dual subgradient algorithm is discussed. Employing FGOR to solve the global consensus problem is discussed in \S~\ref{sec:Application:The-General-Consensus-Problem}. Numerical experiments are presented in \S~\ref{sec:numerical-Results}. Finally, we conclude the paper in \S~\ref{sec:conclusion}, followed by the appendices.

\section{FGOR Characteristic of the Conjugate} \label{Sec:Flatness-Conjugate-Function}

In this section, we demonstrate an appealing property of the conjugate function $f^*$ of $f$ [\cf problem~\eqref{eq:optimization_prob_main_extended}], that is, FGOR, which we will employ to devise our stepsize rules. Roughly speaking, we show that there is a specific region in the domain of $f^*$ such that for all points $\nu$ in the region, there is a ray originating from $\nu$ along which the gradients of $f^*$ remain constant. Besides, note that the dual function, $g$ is simply a restriction of $-f^*$. As a result, further, we show that $g$ directly inherits the FGOR characteristic from $-f^*$ if the restriction conforms to certain conditions. Then we leverage the FGOR characteristic of $g$ to design a precursory yet simple and communication-efficient stepsize rule that can be prepended with other existing ones for improved convergences.    

Let us start by making a couple of assumptions about $f_0$, the set $\mathcal{Y}$, and the subdifferential $\partial f(\mathbf{\bar y})$, where $\mathbf{\bar y}\in\texttt{bnd}~\mathcal{Y}$.
\begin{assump}\label{Assumption:PropStCvx}
The function $f_0$ is lower semicontinuous, proper, and strictly convex. The set $\mathcal{Y}$ is convex and compact. 
\end{assump}
\addb{Note that since $f_0$ is strictly convex, $f^*$ is differentiable on $\R^n$, and so is $g$, \cf~Appendix~\ref{Appendix:Imp-Assump-Conj-Func}, Lemma~\ref{Lemma:ConjugateDifferentiability}, Cor.~\ref{Lemma:UnconstrainedMin-of-f}.} We further note that, the compact constraint set is not typically a restriction in practice. For example, it is a common approach in many machine learning problems that impose regularization as a constraint \cite{singh_2018,gouk_2021}. It can be shown that the problem of convex regularized minimization is equivalent to its corresponding regularization-constrained formulation \cite{bach_2012}, \addbb{\cf~Appendix~\ref{Appendix:Regularization-constrained-Formulation}}. Moreover, decision variables pertaining to feasible resources associated to engineering problems are usually finite \cite{Boyd-Vandenberghe-04,Calafiore_2014}. \addbb{Generalization to nonconvex $f_0$ is deferred to Appendix~\ref{Sec:Appendix:Extensions-nonconvex-case}}.

\begin{assump}\label{Assumption:Exclusion-of-Some-Regions}
For all $\mathbf{\bar y}\in\texttt{bnd}~\mathcal{Y}$, there exists a sequence $\{\mathbf{y}_k\}_{k=1}^{\infty}$, $\mathbf{y}_k\in\mathcal{Y}$, with $\lim_{k\to\infty}\mathbf{y}_k=\mathbf{\bar y}$ such that $\lim\sup_{k\rightarrow\infty} \ \|\boldsymbol{\nu}_k\|_2 \neq \infty  \ \mbox{with} \ \boldsymbol{\nu}_k\in\partial f(\mathbf{y}_k)$.
\end{assump}
By the assumption above, one can exclude cases of infinite gradients at the boundary of $\mathcal{Y}$. The function $l=l_0+\delta_{\mathcal{Y}}$ in Example~\ref{ex:y-power-4} [\cf~Appendix~\ref{Appendix:Imp-Assump-Conj-Func}] is a case where the assumption above breaks since the gradient of the function $l_0$ tends to $\infty$ or $-\infty$ as $\bar y$ tends to $1$ or $-1$, respectively. 

\begin{figure}[!t]
\centering
\subfloat[]{\includegraphics[width=0.47\linewidth]{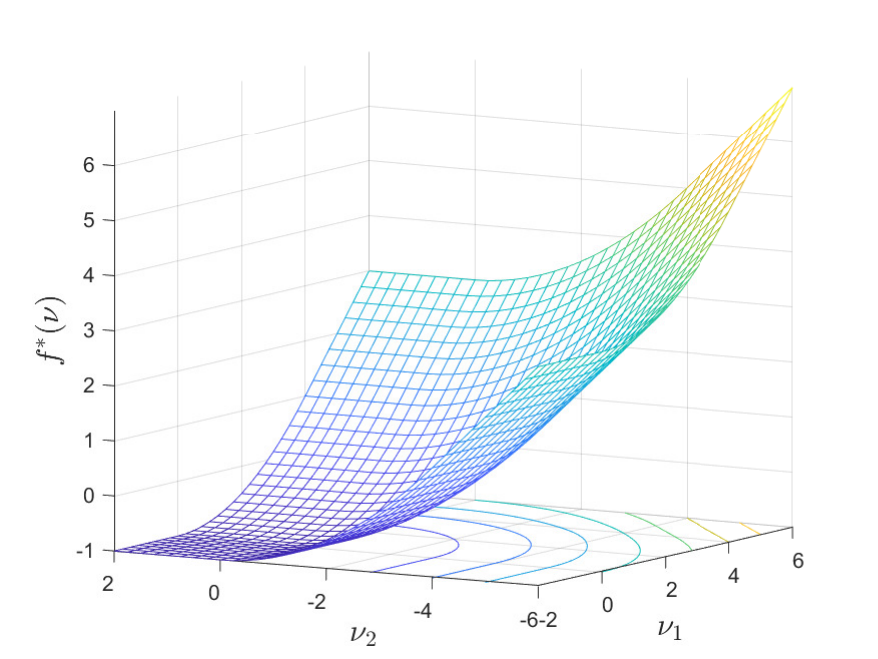}%
\label{subfig:Conjugate-of-f1}}
\hfil
\subfloat[]{\includegraphics[width=0.46\linewidth]
{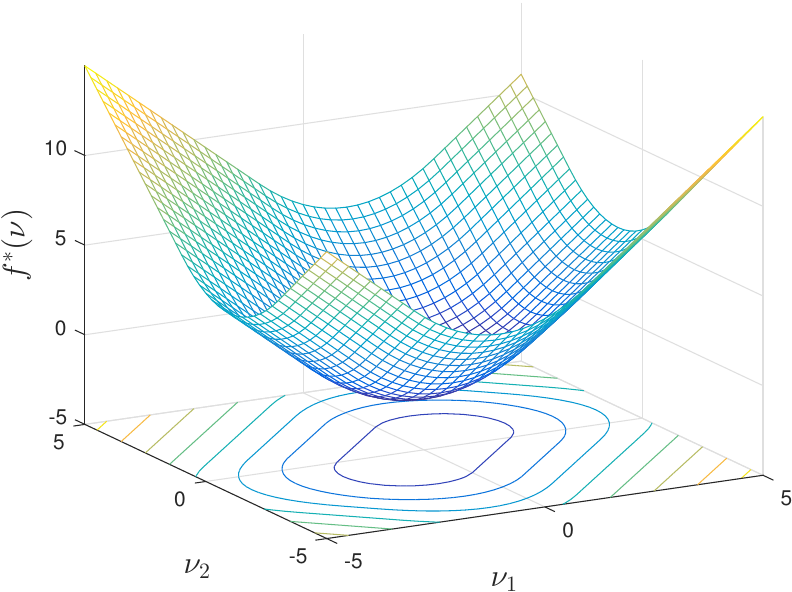}%
\label{subfig:Conjugate-of-f2}}
\\
\subfloat[]{\includegraphics[width=0.47\linewidth]
{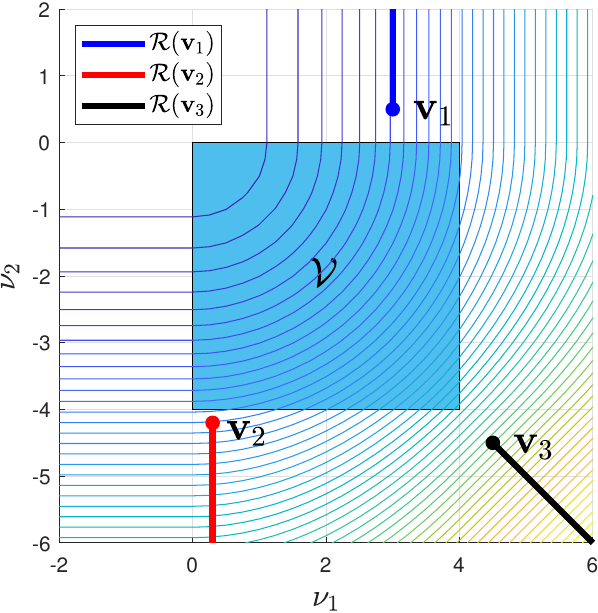}
\label{subfig:LevelSets-of-f1}}
\hfil
\subfloat[]{\includegraphics[width=0.47\linewidth]{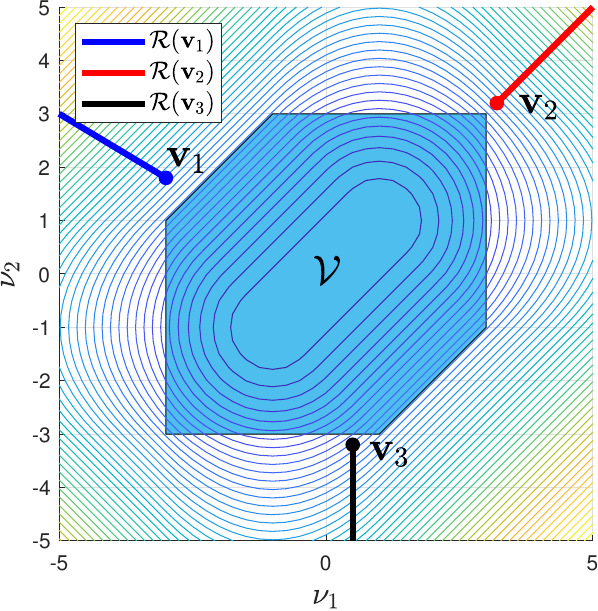}
\label{subfig:LevelSets-of-f2}}
\caption{An illustration of Proposition~\ref{Prop:flat-region}. \protect\subref{subfig:Conjugate-of-f1} $f^*_1$ of $f_1+\delta_{\mathcal{Y}_1}$. \protect\subref{subfig:Conjugate-of-f2} $f^*_2$ of $f_2+\delta_{\mathcal{Y}_2}$. \protect\subref{subfig:LevelSets-of-f1} Level sets of $f^*_1$ and associated set $\mathcal{V}$: $\nabla f_1^*$ over each ray remains intact, where $\mathcal{R}(\mathbf{v}_1)$, $\mathcal{R}(\mathbf{v}_2)$, and $\mathcal{R}(\mathbf{v}_3)$ originates at $\mathbf{v}_1=[3 \ 0.5]\tran$ extending along $\boldsymbol{\eta}_1=[0 \ 1]\tran$, $\mathbf{v}_2=[0.3 \ -4.2]\tran$ extending along $\boldsymbol{\eta}_2=[0 \ -1]\tran$, and $\mathbf{v}_3=[4.5 \ -4.5]\tran$ extending along $\boldsymbol{\eta}_3=[1 \ -1]\tran$, respectively. \protect\subref{subfig:LevelSets-of-f2} Level sets of $f^*_2$ and associated set $\mathcal{V}$: $\nabla f_2^*$ over each ray remains intact, where $\mathcal{R}(\mathbf{v}_1)$, $\mathcal{R}(\mathbf{v}_2)$, and $\mathcal{R}(\mathbf{v}_3)$ originates at $\mathbf{v}_1=[-3 \ 1.8]\tran$ extending along $\boldsymbol{\eta}_1=[-5 \ 3]\tran$, $\mathbf{v}_2=[3.2 \ 3.2]\tran$ extending along $\boldsymbol{\eta}_2=[1 \ 1]\tran$, and $\mathbf{v}_3=[0.5 \ -3.2]\tran$ extending along $\boldsymbol{\eta}_3=[0 \ -1]\tran$,~respectively.}
\label{fig:Illustration_of_mathcal_Y}
\end{figure}

Next, we \addb{present an important lemma and} a few remarks that are useful in deriving the FGOR characteristic of~$f^*$. 
\addb{\begin{lemma}\label{Lemma:intV_mabs_subgrads_of_intY}
Let 
\begin{equation}\label{eq:define-region-V}
    \mathcal{V}=\texttt{cl}~\{\boldsymbol{{\nu}}\in\partial f_0(\mathbf{y}) \ | \ \mathbf{y} \in \texttt{int}~\mathcal{Y}\}.
\end{equation}
Then $\boldsymbol{\nu}_0\in\texttt{int}~\mathcal{V} \iff  \exists~\mathbf{y}\in\texttt{int}~\mathcal{Y}~ \mbox{s.t.,}~\boldsymbol{\nu}_0\in\partial f(\mathbf{y})$.

\end{lemma}
\begin{IEEEproof}
Let $\boldsymbol{\nu}_0\in\texttt{int}~\mathcal{V}$. Then, by the definition of $\mathcal{V}$, $\exists~\mathbf{y}\in\texttt{int}~\mathcal{Y}$ s.t., $\boldsymbol{\nu}_0\in\partial f_0(\mathbf{y})$. Since $f(\mathbf{y})=f_0(\mathbf{y})\,\forall~\mathbf{y}\in\texttt{int}~\mathcal{Y}$, $\boldsymbol{\nu}_0\in\partial f(\mathbf{y})$. Conversely, suppose that $\exists~\mathbf{y}\in\texttt{int}~\mathcal{Y}$ s.t., $\boldsymbol{\nu}_0\in\partial f(\mathbf{y})$. Again, since $f(\mathbf{y})=f_0(\mathbf{y})\,\forall~\mathbf{y}\in\texttt{int}~\mathcal{Y}$, $\boldsymbol{\nu}_0\in\partial f_0(\mathbf{y})$. Then, by the definition of $\mathcal{V}$, together with the strict convexity of $f_0$, $\boldsymbol{\nu}_0\in~\texttt{int}~\mathcal{V}$.    
\end{IEEEproof}}
We illustrate the set $\mathcal{V}$ [\cf~\eqref{eq:define-region-V}] using the following example for clarity. 
\begin{example}[An Illustration of the Set $\mathcal{V}$]\label{ex:An-Illustration-to-V}
Two simple functions $f_1$ and  $f_2$ are considered with constraint sets $\mathcal{Y}_1$ and $\mathcal{Y}_2$, respectively.
\begin{enumerate}
    \item $f_1(\mathbf{y})=2\|\mathbf{y}\|_2^2+1$,~where $\mathcal{Y}_1=\left\{ \mathbf{y}{\in}\R^2 \ | \  -[0 \ 1]\tran\leq \mathbf{y} \leq [1 \ 0]\tran\right\}$.
\item $f_2(\mathbf{y})=0.5\max\left\{\|\mathbf{y}-\mathbf{1}_{2}\|_2^2,\|\mathbf{y}+\mathbf{1}_{2}\|_2^2\right\}$, where $\mathcal{Y}_2=\left\{ \mathbf{y}\in\R^2 \ | \ -[2 \ 2]\tran\leq \mathbf{y} \leq [2 \ 2]\tran\right\}$.
\end{enumerate}  
Then \addb{$\nabla f_1(\mathbf{y})=4\mathbf{y}$, where $\mathbf{y}\in\mathcal{Y}_1$. Thus, it is straightforward to see that} the set $\mathcal{V}$ associated to $f_1$ is given by $\mathcal{V}=\{{\boldsymbol{\nu}\in\R^2\ | \ -[0 \ 4]\tran\leq \boldsymbol{\nu}\leq [4 \ 0]\tran\}}$, \cf~\figurename~\ref{fig:Illustration_of_mathcal_Y}\subref{subfig:LevelSets-of-f1}. \addb{Next, to illustrate the set $\mathcal{V}$ associated to $f_2$, we let $f_{21}=0.5\|\mathbf{y}-\mathbf{1}_2\|_2^2$ and $f_{22}=0.5\|\mathbf{y}+\mathbf{1}_2\|_2^2$. Then $\partial f_2(\mathbf{y})=\texttt{con}\{\nabla f_{2i}(\mathbf{y}) \ | \ f_{2i}(\mathbf{y})=f_{2}(\mathbf{y}), i=1,2\}$, where $\mathbf{y}\in\mathcal{Y}_2$ since $f_2$ is the pointwise maximum of $f_{21}$ and $f_{22}$ \cite{Boyd-EE364b-subgrad-22}. Therefore, it is easily seen that
\begin{equation}\label{eq:subdifferential-f2}
    \partial f_2(\mathbf{y})=\begin{cases}
    \displaystyle \nabla f_{22}(y_1,y_2) \ ;  \ \text{$y_2>-y_1$}\\
         \displaystyle\nabla f_{21}(y_1,y_2) \  ; \ \text{$y_2<-y_1$}\\
         \displaystyle \texttt{con}\{\nabla f_{21}(y_1,y_2),\nabla f_{22}(y_1,y_2)\} \ ;  \ \text{$y_2=-y_1$},
    \end{cases}  
\end{equation}
where $\mathbf{y}=[y_1 \, y_2]\tran\in\mathcal{Y}_2$. Thus \eqref{eq:subdifferential-f2}, together with simple calculations, yields that} the set $\mathcal{V}$ associated to $f_2$ is given by the polyhedron with vertices $-[3 \ 3]\tran$, $[-3 \ 1]\tran$, $[-1 \ 3]\tran$, $[3 \ 3]\tran$, $[3 \ -1]\tran$, and $[1 \ -3]\tran$, \cf~\figurename~\ref{fig:Illustration_of_mathcal_Y}\subref{subfig:LevelSets-of-f2}. 
\end{example}
\begin{remark}\label{Remark: subdifferential-indicator-function-at-boundary}
     Let $\mathbf{\bar y}\in\texttt{bnd}~\mathcal{Y}$. Then  $\partial \delta_{\mathcal{Y}}(\mathbf{\bar y})=N_{\mathcal{Y}}(\mathbf{\bar y})$, where $N_{\mathcal{Y}}(\mathbf{\bar y})$ denotes the normal cone of the set $\mathcal{Y}$ at $\mathbf{\bar y}$.  
\end{remark}

\begin{remark}\label{Remark: subdifferential-f-at-boundary}
Let $\mathbf{\bar y}\in\texttt{bnd}~\mathcal{Y}$. Then $\partial f(\mathbf{\bar y})=\partial f_0(\mathbf{\bar y})+N_{\mathcal{Y}}(\mathbf{\bar y})$.

\end{remark}
Remark~\ref{Remark: subdifferential-f-at-boundary} is an immediate result due to the convexity of $f_0$ and $\delta_{\mathcal{Y}}$. More specifically, since $~{f=f_0+\delta_{\mathcal{Y}}}$ [\cf~\eqref{eq:optimization_prob_main_extended}] the convexity of $f_0$ and $\delta_{\mathcal{Y}}$ suggests that $~{\partial f=\partial f_0+\partial \delta_{\mathcal{Y}}}$. This together with Remark~\ref{Remark: subdifferential-indicator-function-at-boundary} yield Remark~\ref{Remark: subdifferential-f-at-boundary}.   

Finally, the following proposition establishes the FGOR characteristic of $f^*$.

\begin{prop}\label{Prop:flat-region}
 Suppose Assumption~\ref{Assumption:PropStCvx} and Assumption~\ref{Assumption:Exclusion-of-Some-Regions} hold. Then, for all $\boldsymbol{\nu}\in\R^n\setminus\texttt{int~}\mathcal{V}$, there exists $\boldsymbol{\eta}\in\R^n$ such that  for all $\alpha\geq0$, $\boldsymbol{\nu}+\alpha \boldsymbol{\eta}\in\R^n\setminus\texttt{int}~\mathcal{V}$ and $\nabla f^*(\boldsymbol{\nu}+\alpha \boldsymbol{\eta})$ is a constant vector \addb{that corresponds to some $\bar{\mathbf{y}}\in\texttt{bnd}~\mathcal{Y}$.} 
\end{prop}
\begin{IEEEproof}
\addb{Let $\boldsymbol{\nu}\in \R^n\setminus \texttt{int}~\mathcal{V}$. Note that the equivalence established in Lemma~\ref{Lemma:intV_mabs_subgrads_of_intY} is equivalent to the contrapositive
\begin{equation}\label{remark1-contraposition-1}
\boldsymbol{\nu}\in\R^n\setminus\texttt{int}~\mathcal{V}\iff \forall~\mathbf{y}\in\texttt{int}~\mathcal{Y},\,\boldsymbol{\nu}\not\in\partial f(\mathbf{y}). 
\end{equation}
    However, since $\boldsymbol{\nu}\in\texttt{dom}~f^*$, it follows by the inversion rule for subgradient relations \cite[Prop.~11.3]{Rockafellar-98} that $\boldsymbol{\nu}\in\partial f(\mathbf{y})$ for some $\mathbf{y}\in\mathcal{Y}$. Thus, \eqref{remark1-contraposition-1} yields
    \begin{equation}\label{remark1-contraposition-2} \boldsymbol{\nu}\in\R^n\setminus\texttt{int}~\mathcal{V}\iff \exists~\mathbf{\bar y}\in\texttt{bnd}~\mathcal{Y},\,\mbox{s.t.,}~\boldsymbol{\nu}\in\partial f(\mathbf{\bar y}).
    \end{equation}
    Then, from~\eqref{remark1-contraposition-2} we have $\boldsymbol{\nu}\in\partial f(\mathbf{\bar y})$ for some $\mathbf{\bar y}\in\texttt{bnd}~\mathcal{Y}$}, i.e.,
    \begin{equation}\label{eq:FGOR-proof-1}
\boldsymbol{\nu}=\bar{\boldsymbol{\nu}}+\boldsymbol{\eta},
\end{equation}%
where $\bar{\boldsymbol{\nu}}\in\partial f_0(\mathbf{\bar y})$ and $\boldsymbol{\eta}\in N_{\mathcal{Y}}(\mathbf{\bar y})$, \cf~Remark~\ref{Remark: subdifferential-f-at-boundary}. Now consider the ray 
\begin{equation}\label{eq:prop-FGOR-of-conjugate-2}
    \{\boldsymbol{\nu}+\alpha \boldsymbol{\eta} \ | \ \alpha\geq0 \}=\{\bar{\boldsymbol{\nu}}+(1+\alpha) \boldsymbol{\eta} \ | \ \alpha\geq0 \},
\end{equation}
where the equality follows from \eqref{eq:FGOR-proof-1}. Since $\bar{\boldsymbol{\nu}}\in\partial f_0(\bar{\mathbf{y}})$ and $\boldsymbol{\eta}\in N_{\mathcal{Y}}(\bar{\mathbf{y}})$, $\bar{\boldsymbol{\nu}}+(1+\alpha) \boldsymbol{\eta}\in\partial f(\mathbf{\bar{y}})$ for all $\alpha\geq 0$, \cf~Remark~\ref{Remark: subdifferential-f-at-boundary}. Therefore, \cite[Prop.~11.3]{Rockafellar-98}, together with~ \eqref{eq:prop-FGOR-of-conjugate-2}, we have 
\begin{equation}\label{eq:FGOR-of}
\nabla f^*(\boldsymbol{\nu}+\alpha\boldsymbol{\eta})=\mathbf{\bar y} 
\end{equation}
for all $\alpha\geq0$. Moreover, since $\mathbf{\bar y}\in\texttt{bnd}~\mathcal{Y}$, equivalence \addb{\eqref{remark1-contraposition-2}} confirms that for all $\alpha\geq0$, $\boldsymbol{\nu}+\alpha \boldsymbol{\eta}\in\R^n\setminus\texttt{int}~\mathcal{V}$.
\end{IEEEproof}

Proposition~\ref{Prop:flat-region} indicates that for all $\boldsymbol{\nu}\in\R^n\setminus\texttt{int}~\mathcal{V}$ there exists a ray $\mathcal{R}_{\boldsymbol{\nu}} =\{\boldsymbol{\nu}+\alpha \boldsymbol{\eta} \ | \ \alpha\geq0 \}\subseteq\R^n\setminus\texttt{int}~\mathcal{V}$ that originates from $\boldsymbol{\nu}$ along which the gradients of the conjugate $f^*$ remain constant. To illustrate this phenomenon associated to Proposition~\ref{Prop:flat-region}, we provide the following example.

\begin{example}[An Illustration of Proposition~\ref{Prop:flat-region}]\label{ex:Exclusion-of-Affine-Regions}
Consider $f_1$ and $f_2$ in Example~\ref{ex:An-Illustration-to-V}. The conjugate functions $f_1^*$ and $f_2^*$ of $f_1+\delta_{\mathcal{Y}_1}$ and $f_2+\delta_{\mathcal{Y}_2}$ are depicted in \figurename~\ref{fig:Illustration_of_mathcal_Y}\subref{subfig:Conjugate-of-f1} and \figurename~\ref{fig:Illustration_of_mathcal_Y}\subref{subfig:Conjugate-of-f2}, respectively. The level sets of the conjugate functions are shown in \figurename~\ref{fig:Illustration_of_mathcal_Y}\subref{subfig:LevelSets-of-f1} and \figurename~\ref{fig:Illustration_of_mathcal_Y}\subref{subfig:LevelSets-of-f2} with their corresponding sets $\mathcal{V}$, \cf~Lemma~\ref{Lemma:intV_mabs_subgrads_of_intY}. The rays originated from some random points along which the gradients of $f_1^*$ and $f_2^*$ remain constant are shown in the same plots.
\end{example} 

\section{FGOR Characteristic of the Dual Function}\label{sec:Flatness-of-g}  
The FGOR characteristic of the dual function~$g$ for problem~\eqref{eq:optimization_prob_main_extended} is directly inherited from the conjugate function $f^*$, \cf~\S~\ref{Sec:Flatness-Conjugate-Function}. To this end, we start by letting $\boldsymbol{\lambda}\in\R^m$ be the Lagrange multiplier associated to the equality constraint of problem~\eqref{eq:optimization_prob_main_extended}. Then the dual function is given by %
\begin{align}
g(\boldsymbol{\lambda})&= \underset{\mathbf{y}\in \R^n}{\inf} \ \left(f(\mathbf{y}) - \boldsymbol{\lambda}\tran (\mathbf{A}\mathbf{y}-\mathbf{b})\right)\nonumber\\ 
&=-\underset{\mathbf{y}\in \R^n}{\sup} \ \left(\boldsymbol{\lambda}\tran \mathbf{A}\mathbf{y}-f(\mathbf{y})\right)+\boldsymbol{\lambda}\tran\mathbf{b}  \nonumber\\ 
&=-f^*(\mathbf{A}\tran \boldsymbol{\lambda})+\boldsymbol{\mathbf{b}\tran\lambda}. \label{eq:dual-as-a-restriction-3}
\end{align} 
From \eqref{eq:dual-as-a-restriction-3}, it is clear that the dual function $g$ is based on a restriction of $-f^*$ to a linear space. Moreover, from basic calculus rules, we have
\begin{equation}\label{eq:dual-grad-and-conjugate-grad}
    \nabla g(\boldsymbol{\lambda})=-\mathbf{A}\nabla f^*(\mathbf{A}\tran\boldsymbol{\lambda}) + \mathbf{b}.
\end{equation}
The gradient identity \eqref{eq:dual-grad-and-conjugate-grad} suggests that if $\nabla f^*(\mathbf{A}\tran\boldsymbol{\lambda})$ is constant for some $\mathcal{R}\subseteq R(\mathbf{A}\tran)$, then so is the gradient of~$g$, $\nabla g(\boldsymbol{\lambda})$ for all $\{\boldsymbol{\lambda} \ | \ \mathbf{A}\tran\boldsymbol{\lambda}\in\mathcal{R}\}$. This is essentially the condition based on what the FGOR characteristic of $g$ is established. 

We now formalize the preceding arguments. To this end, we start by defining the collection $\mathcal{A}$ of rays in $\R^n\setminus\texttt{int}~\mathcal{V}$, where the gradients $\nabla f^*$ over each ray remain intact,~\cf~Proposition~\ref{Prop:flat-region}. In particular, for any $\boldsymbol{\nu}\in\R^n\setminus\texttt{int}~\mathcal{V}$, a ray originated at $\boldsymbol{\nu}$, denoted $\mathcal{R}(\boldsymbol{\nu})$ over which $\nabla f^*$ remains intact is given~by
\begin{equation}\label{eq:constant-gradient-ray-in-conjugate}
     \mathcal{R}(\boldsymbol{\nu})= \{\boldsymbol{\nu}+\alpha\boldsymbol{\eta}(\boldsymbol{\nu}) \ | \ \alpha\geq 0\},
\end{equation}
where $\boldsymbol{\eta}(\boldsymbol{\nu})\in\R^n$ depends on $\boldsymbol{\nu}$, \cf~Proposition~\ref{Prop:flat-region}. Consequently, the collection of all such rays $\mathcal{A}$ is given by
\begin{equation}\label{set-of-constant-gradient-rays}
\mathcal{A} = \{\mathcal{R}(\boldsymbol{\nu}) \ | \ \boldsymbol{\nu}\in\R^n\setminus\texttt{int}~\mathcal{V} \}.
\end{equation}
It is now necessary to impose the following assumption: 
\begin{assump}\label{Assumption:dual-flat-assump}
 $\exists~\mathcal{R}\in\mathcal{A}$ such that $\mathcal{R}\subseteq R(\mathbf{A}\tran)$.
    \end{assump}
\addb{Assumption~\ref{Assumption:dual-flat-assump} plays a key role when establishing FGOR characteristics of $g$ in general. In specific cases, its validity must be verified within the given context. For instance, the consensus problem, which is widely utilized across various signal processing and machine learning application domains~\cite[p.~49]{Boyd-Parikh-Chu-Peleato-Eckstein-2010} is a case where Assumption~\ref{Assumption:dual-flat-assump} holds, \cf~Remark~\ref{Remark:Assumption-3}, \S~\ref{subsec:Consensus-problem-ell-infinity-norm-constraints}, Lemma~\ref{prop:consensus-problem-feasible-y-bar}. We defer the details to $\S~\ref{sec:Application:The-General-Consensus-Problem}$ while outlining a toy example to demonstrate the idea.} 
\begin{example}[An Illustration of Assumption~\ref{Assumption:dual-flat-assump}] \label{eq:Assumption:dual-flat-assump}
Consider minimizing $f_2$ in Example~\ref{ex:An-Illustration-to-V} with an equality constraint $\mathbf{Ay}=\mathbf{0}$. Then the corresponding optimization problem is 
\begin{equation}
\label{eq:Ex:Assumption:dual-flat-assump}
\begin{array}{ll}
\mbox{minimize} & f_2(\mathbf{y})\\
\mbox{subject to} & \mathbf{y}\in\mathcal{Y}_2\\
\mbox{} & \mathbf{A}\mathbf{y} = \mathbf{0}.
\end{array}
\end{equation}
The sets $R(\mathbf{A}_1\tran)$ and $R(\mathbf{A}_2\tran)$, where $\mathbf{A}_1=[-5 \ 3]$ and $\mathbf{A}_2=[1 \ 1]$ are depicted in \figurename\ref{Fig:example-prop2-dual-functions}\subref{subfig:example-dual-function-restrictions}. Moreover, the rays $\mathcal{R}(\mathbf{v}_1)$ and $\mathcal{R}(\mathbf{v}_2)$ [\cf~\eqref{eq:constant-gradient-ray-in-conjugate}], where $\mathbf{v}_1=[-3 \ 1.8]\tran$ and $\mathbf{v}_2=[3.2 \ 3.2]\tran$, are also shown in the same figure. Clearly, $\mathcal{R}(\mathbf{v}_1)\subseteq R(\mathbf{A}_1\tran)$ and $\mathcal{R}(\mathbf{v}_2)\subseteq R(\mathbf{A}_2\tran)$, respectively. Thus Assumption~\ref{Assumption:dual-flat-assump} holds for problem~\eqref{eq:Ex:Assumption:dual-flat-assump}. 
\end{example}
\addbb{We note that Assumption~3 is inherently problem-dependent and may not hold in all cases. The lack of a general verification criterion for this assumption constitutes a caveat of the proposed approach, which must be addressed separately using problem-specific ingenuity.} We are now ready to show our main result about the FGOR characteristic of the dual function $g$, \cf~\eqref{eq:dual-as-a-restriction-3}.
\begin{corollary}\label{Prop:flat-region-dual-function}
Suppose Assumption~\ref{Assumption:PropStCvx}, Assumption~\ref{Assumption:Exclusion-of-Some-Regions}, and Assumption~\ref{Assumption:dual-flat-assump} hold. Then $\exists$ $\boldsymbol{\lambda}\in \R^m$, $\exists$ $\boldsymbol{\mu}\in\R^m$ such that for all $\alpha\geq0$, $\nabla g(\boldsymbol{\lambda}+\alpha \boldsymbol{\mu})$ is a constant vector.
\end{corollary}

\begin{IEEEproof}
\addb{The proof is a direct consequence of Proposition~\ref{Prop:flat-region}, Assumption~\ref{Assumption:dual-flat-assump}, and $\mathbf{A}\tran\in \R^{n\times m}$ with $m<n$ and $\texttt{rank}~\mathbf{A}=m$.} More specifically, from Assumption~\ref{Assumption:dual-flat-assump}, there exists a ray $\mathcal{R}\in\mathcal{A}$ such that $\mathcal{R}\subseteq R(\mathbf{A}\tran)$. Let $\mathcal{R}~{=\{\boldsymbol{\nu}+\alpha\boldsymbol{\eta} \ | \ \alpha\geq 0\}}$ for some $\boldsymbol{\nu},\boldsymbol{\eta}\in\R^n$. Now consider the inverse image $\mathcal{I}_{\mathcal{R}}$ of $\mathcal{R}$ under the linear mapping $\mathbf{A}\tran$. In particular, we have
\begin{equation}\label{eq:Inverse-Image-of-Conjugate-Ray-1}
     \mathcal{I}_{\mathcal{R}} =   \{\boldsymbol{\kappa} \ | \ \mathbf{A}\tran \boldsymbol{\kappa}\in\mathcal{R}\}= \{\boldsymbol{\lambda}+\alpha\boldsymbol{\mu} \ | \ \alpha\geq 0\},
\end{equation}
where $\boldsymbol{\lambda}=(\mathbf{A}\mathbf{A}\tran)^{-1}\mathbf{A}\mathbf{\boldsymbol\nu}$ and $\boldsymbol{\mu}=(\mathbf{A}\mathbf{A}\tran)^{-1}\mathbf{A}\boldsymbol{\eta}$, and the last equality follows from basic linear algebraic steps. Thus, from \eqref{eq:dual-grad-and-conjugate-grad}, together with Assumption~\ref{Assumption:dual-flat-assump}, we conclude that for all $\alpha\geq 0$, $ \nabla g(\boldsymbol{\lambda}+\alpha\boldsymbol{\mu})$ is a constant vector.   
\end{IEEEproof}
In the following example, we describe two instances of dual functions associated to problem~\eqref{eq:Ex:Assumption:dual-flat-assump} to demonstrate the assertions of Corollary~\ref{Prop:flat-region-dual-function}. 
\begin{example}[An Illustration of
Corollary~\ref{Prop:flat-region-dual-function}]
\figurename~\ref{Fig:example-prop2-dual-functions}\subref{subfig:example-dual-functions} depicts the dual functions $g_1$ and $g_2$ of problem~\eqref{eq:Ex:Assumption:dual-flat-assump} when $\mathbf{A}=\mathbf{A}_1=[-5 \ 3]$ and $\mathbf{A}=\mathbf{A}_2=[1 \ 1]$, respectively. The functions $g_1$ and $g_2$ are based on restrictions of $-f^*$ to $R(\mathbf{A}_1\tran)$ and $R(\mathbf{A}_2\tran)$, respectively, \cf~\figurename~\ref{Fig:example-prop2-dual-functions}\subref{subfig:example-dual-function-restrictions}. The figure shows both $g_1$ and $g_2$ have regions with constant gradients.
\end{example}
Finally, let us remark that the conjugate function and the dual function are both gradient-Lipschitz continuous. Assumption~\ref{Assumption:PropStCvx} is sufficient to yield such Lipschitzian properties, although we won't take this up in our subsequent stepsize designs. We note that the aforementioned Lipschitzian properties are more general than the existing results. Details are deferred
to Appendix~\ref{Appendix:Imp-Assump-Conj-Func} for completeness.

In the next section, we describe how the FGOR characteristic of $g$ [\cf~Corollary~\ref{Prop:flat-region-dual-function}] is used to devise a simple stepsize rule that can precede existing ones for improved convergences of the associated subgradient method.

\section{On the Application of FGOR on Dual Subgradient Algorithm}\label{sec:Dual-Subgradient-Algorithm-The-Proposed- Approach}

FGOR is the prime characteristic we use in the sequel for devising our simple stepsize rule. To formalize the exposition, let us first consider the dual problem associated to problem~\eqref{eq:optimization_prob_main_extended}: 
 \begin{equation}\label{eq:dual-problem}
    \underset{\boldsymbol{\lambda}\in\R^m}{\text{maximize}} \quad \left(g(\boldsymbol{\lambda})=\underset{\mathbf{y}\in \mathcal{Y}}{\inf} \ \left(f_0(\mathbf{y}) - \boldsymbol{\lambda}\tran (\mathbf{A}\mathbf{y}-\mathbf{b})\right)\right).
\end{equation}
Then the subgradient method to solve \eqref{eq:dual-problem} is given by 
\begin{equation}\label{eq:Lambda-Update}
    \boldsymbol{\lambda}^{(k+1)} = \boldsymbol{\lambda}^{(k)} - \gamma_k\mathbf{s}^{{(k)}},
\end{equation}
where $\gamma_k>0$ is the stepsize, $\mathbf{s}^{(k)}$ is a subgradient of~$-g$ at $\boldsymbol{\lambda}^{(k)}\in\R^{m}$, and $k$ represents the iteration index. 

\begin{figure}[!t]
\centering
\subfloat[]{\includegraphics[width=0.49\linewidth]{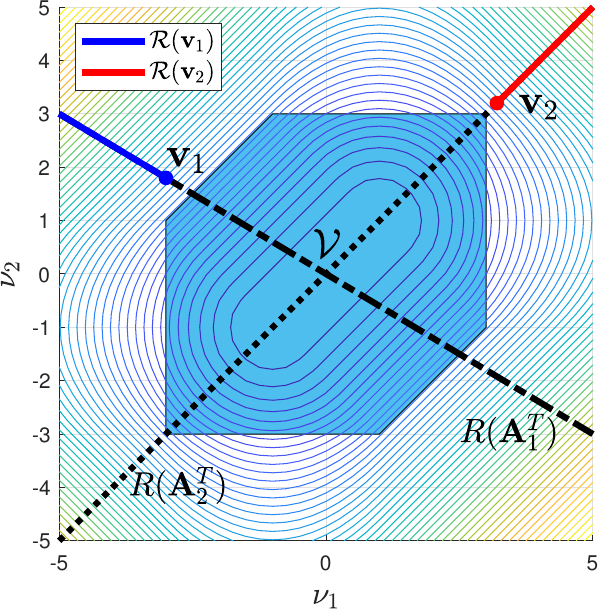}%
\label{subfig:example-dual-function-restrictions}}
 \hfil
 \subfloat[]{\includegraphics[width=0.5\linewidth]{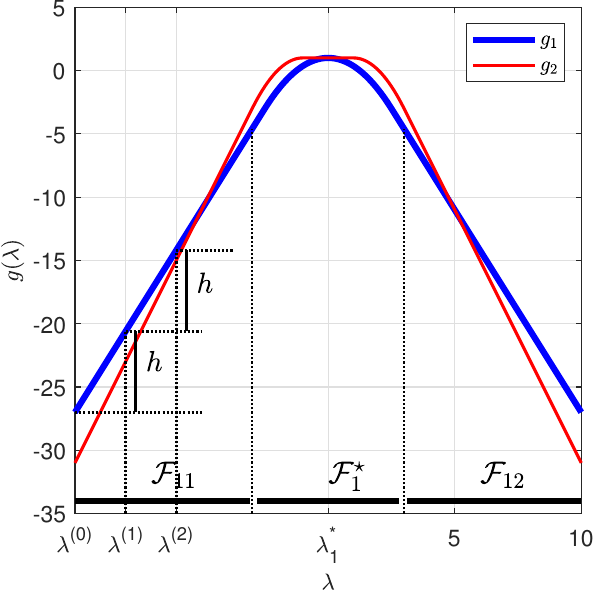}%
 \label{subfig:example-dual-functions}}
\caption{An illustration of Corollary~\ref{Prop:flat-region-dual-function}. \protect\subref{subfig:example-dual-function-restrictions} Level sets of $f^*_2$ and associated set $\mathcal{V}$: ${\mathcal{R}(\mathbf{v}_1){\subseteq}R(\mathbf{A}_1\tran)}$ and ${\mathcal{R}(\mathbf{v}_2){\subseteq} R(\mathbf{A}_2\tran)}$. \protect\subref{subfig:example-dual-functions} Dual functions of problem~\eqref{eq:Ex:Assumption:dual-flat-assump} when $\mathbf{A}=\mathbf{A}_1$ and $\mathbf{A}=\mathbf{A}_2$: FGOR region of $g_1$ is $\mathcal{F}_1=\mathcal{F}_{11}\cup\mathcal{F}_{12}$, and $\mathcal{F}_1^\star$ is the region in which the dual optimal solution $\lambda_1^\star=0$ of $g_1$ resides. The points $\lambda^{(0)}=-10$, $\lambda^{(1)}=-8$, and $\lambda^{(2)}=-6$ of the subgradient method \eqref{eq:Lambda-Update} are due to the constant stepsizes taken in the FGOR region $\mathcal{F}_{11}$.}  
\label{Fig:example-prop2-dual-functions}
\end{figure}

 We now give some insight into how the stepsize $\gamma_k$ values should be chosen by considering problem~\eqref{eq:Ex:Assumption:dual-flat-assump} and its associated dual function $g_1$ depicted in \figurename~\ref{subfig:example-dual-functions}. The idea is very simple. Roughly speaking, the stepsize rule is to use a \emph{constant} stepsize as long as $\lambda^{(k)}$ is in a region of FGOR [\cf~$\mathcal{F}_1=\mathcal{F}_{11}\cup\mathcal{F}_{12}$, \figurename~\ref{subfig:example-dual-functions}], making the subgradient method yield a constant increase in the dual objective function in every iteration. For instance, moving from $\lambda^{(0)}$ to $\lambda^{(1)}$, from $\lambda^{(1)}$ to $\lambda^{(2)}$, and so forth is due to constant stepsizes since $\lambda^{(0)},  \lambda^{(1)}, \lambda^{(2)}\in \mathcal{F}_{11}$, which in turn makes the subgradient method able to have constant increments,~\cf~\figurename~\ref{subfig:example-dual-functions}. The important point is to make the subgradient method~\eqref{eq:Lambda-Update} exit the region $\mathcal{F}_1$ effectively and reach quickly the region ${\mathcal{F}}_1^\star$ in which the dual solution $\lambda_1^\star$ resides. As soon as the subgradient method exits the region $\mathcal{F}_1$, the stepsize can be restored to any desirable existing one to guarantee convergence. Note that the constant stepsize rule we propose can precede existing ones by furnishing potential advantages as we will further explain in subsequent sections.

\subsection{Formalizing the FGOR on Dual Subgradient}

Recall from the preceding discussion that the constant stepsize rule yields a constant increase in the dual objective function in every
iteration of the subgradient method~\eqref{eq:Lambda-Update} as long as $\lambda^{(k)}$ is in a region of FGOR. Although it is not evident in the case of the lower dimensional dual function~$g_1$ of the preceding discussion [\cf~\figurename~\ref{subfig:example-dual-functions}], in the context of a general formalization, an extra condition that guarantees such a constant increase of the dual~function is to be considered. More specifically, among the collection of all the rays in $\texttt{dom}~g$ with the FGOR characteristic [\cf~Corollary~\ref{Prop:flat-region-dual-function}], we seek for a ray $\{\boldsymbol{\lambda}+\alpha\boldsymbol{\mu} \ | \ \alpha\geq 0\}$ for which
\begin{equation}\label{eq:condition-dual-function}
    \nabla g(\boldsymbol{\lambda}+\alpha \boldsymbol{\mu})=-\beta{\boldsymbol\mu} 
\end{equation}
for some $\beta>0$. \addb{In other words, besides being a constant vector, condition~\eqref{eq:condition-dual-function} requires that the gradient $\nabla g$, when evaluated along the ray, aligns with the ray itself. We refer to the set of all rays that satisfy condition~\eqref{eq:condition-dual-function} as the restricted-FGOR (RFGOR) region.} \figurename~\ref{Fig:example-FGOR-dual-grads-parellel-to-line} illustrates the idea of the condition by depicting two rays $\{\boldsymbol{\lambda}_1+\alpha\boldsymbol{\mu}_1 \ | \ \alpha\geq 0\}$ and $\{\boldsymbol{\lambda}_2+\alpha\boldsymbol{\mu}_2 \ | \ \alpha\geq 0\}$. Note that the gradient $\nabla g$ is a constant vector over each ray. However, the condition~\eqref{eq:condition-dual-function} is satisfied only over the first ray.

\begin{figure}[!t]
\centering
\subfloat[]
{\begin{tikzpicture}[scale=0.38]
\tikzset{vertex/.style = {shape=circle,draw,minimum size=0.2cm}}
\tikzset{edge/.style = {->,> = latex}}
\draw (0,0) .. controls (1,-0.1) .. (2,-0.5);
\draw (0,0) .. controls (0.2,-1.5) .. (0,-3); 
\draw (2,-0.5) .. controls (2.2,-2) .. (2,-3.5); 
\draw (0,-3) .. controls (1,-3.1) .. (2,-3.5); 

\draw (0.75,-0.13) .. controls (0.8,0.25) .. (0.75,0.5);
\draw (0.75,0.5) .. controls (1.5,0.4) .. (2.75,0);
\draw (2.75,0) .. controls (2.95,-1.5) .. (2.75,-3.3);
\draw (2.75,-3.3) .. controls (2.4,-3.05) .. (2.1,-3);

\draw (1,-1.2) -- (-3,-2.7);
\draw[thick,dotted] (1,-1.2) -- (2.1,-0.8);
\draw (2.1,-0.8) -- (2.4,-0.65);
\draw[thick,dotted] (2.5,-0.65) -- (2.85,-0.5);
\draw (2.85,-0.5) -- (3.2,-0.35);
\draw[thick,dotted] (3.2,-0.35) -- (3.7,-0.1);
\draw (3.7,-0.1) -- (6.5,1.1);
 \draw[fill=black] (5.5,0.65) circle (0.1cm);
  \node[draw=none,align=center, font=\scriptsize,text width = 1.25cm] at (5,1) {$\boldsymbol{\lambda}_1$};
\draw[-stealth] (5.5,1) -- (6.7,1.5);
\node[draw=none,align=center, font=\scriptsize,text width = 1.25cm] at (7.2,1.8) {$\boldsymbol{\mu}_1$};

\draw (-3,-3.2) -- (1.5,-1.75); 
\draw[thick,dotted] (1.5,-1.75) -- (2.1,-1.55);
\draw (2.1,-1.55) -- (2.5,-1.4);
\draw[thick,dotted] (2.5,-1.4) -- (2.9,-1.3);
\draw (2.9,-1.3) -- (3.3,-1.15);
\draw[thick,dotted] (3.3,-1.15) -- (3.5,-1.1);
\draw (3.5,-1.1) -- (7,0);
\draw[fill=black] (5.7,-0.4) circle (0.1cm);
\node[draw=none,align=center, font=\scriptsize,text width = 1.25cm] at (5.3,-1) {$\boldsymbol{\lambda}_2$};
\draw[-stealth] (5.7,-0.74) -- (7.1,-0.29);
\node[draw=none,align=center, font=\scriptsize,text width = 1.25cm] at (7.6,-0.1) {$\boldsymbol{\mu}_2$};

\draw (1.25,0.43) .. controls (1.3,0.75) .. (1.25,1.1);
\draw (1.25,1.1) .. controls (1.5,1.1) and (3,0.8).. (3.5,0.4);
\draw (3.5,0.4) .. controls (3.7,0) and (3.7,-2)  .. (3.5,-3);
\draw (3.5,-3) .. controls (3,-2.7) .. (2.8,-2.7);
\end{tikzpicture}\label{subfig:grad-g-parallel-to-line}}%
 \hfil
 \subfloat[]{\begin{tikzpicture}[scale=0.3]
\tikzset{vertex/.style = {shape=circle,draw,minimum size=0.2cm}}
\tikzset{edge/.style = {->,> = latex}}
\draw (-3.2,-2.5) -- (6.9,2.7);
\draw (-0.5,-0.6) -- (-0.9,-0.9); 
\draw (-0.9,-0.9) -- (-0.7,-1.2);
\draw (-0.5,-0.6) -- (-0.9,-0.9); 
\draw (-0.1,-0.4) -- (0.35,-0.1);
\draw (-0.1,-0.4) -- (0.2,-0.75);

\draw (0.8,0.05) -- (1.2,0.3);
\draw (0.8,0.05) -- (1,-0.3);
\draw[-stealth,line width=0.1mm] (1.5,-0.1) -- (0.6,-0.45); 
\draw[-stealth,line width=0.1mm] (0.6,-0.5) -- (-0.2,-0.9);

\draw[-stealth,line width=0.1mm] (-0.3,-1) -- (-1.15,-1.4);

\node[draw=none,align=center, font=\scriptsize,text width = 0.2cm] at (-2,-0.8) {$\nabla g$};

\draw[-stealth,line width=0.1mm] (2.15,-1.8) -- (1.2,-2.01); 
\draw[-stealth,line width=0.1mm] (1,-2.2) -- (0.3,-2.3);
\draw[-stealth,line width=0.1mm] (0.1,-2.45) -- (-0.8,-2.55);
\node[draw=none,align=center, font=\scriptsize,text width = 0.2cm] at (-1.8,-2.55) {$\nabla g$};

 \draw[fill=black] (4.9,1.7) circle (0.1cm);
 \node[draw=none,align=center, font=\scriptsize,text width = 1.25cm] at (4.8,0.8) {$\boldsymbol{\lambda}_1$};
\draw[-stealth] (5,1.4) -- (6.7,2.3);
\node[draw=none,align=center, font=\scriptsize,text width = 1.25cm] at (7.3,2.1) {$\boldsymbol{\mu}_1$};

\draw(-3,-3.5) --  (7,-0.3);
\draw[fill=black] (5.7,-0.73) circle (0.1cm);
\node[draw=none,align=center, font=\scriptsize,text width = 1.25cm] at (5.3,-1.4) {$\boldsymbol{\lambda}_2$};
\draw[-stealth] (5.7,-1) -- (7.1,-0.55);
\node[draw=none,align=center, font=\scriptsize,text width = 1.25cm] at (7.8,-0.3) {$\boldsymbol{\mu}_2$};

\draw (-1,0) .. controls (-0.2,-1) and (0.6,-1.1) .. (0,-4.5);
\draw (-0.3,0.7) .. controls (1,-1) and (1.3,-1.1) .. (1,-4);
\draw (0.5,1.2) .. controls (2,-1) and (2.3,-1.1) .. (2,-3.5);

\end{tikzpicture}\label{subfig:grad-g-not-parallel-to-line-plane-view}}\caption{An illustration of condition~\eqref{eq:condition-dual-function}: the condition is satisfied only over the ray $\{\boldsymbol{\lambda}_1+\alpha\boldsymbol{\mu}_1 \ | \ \alpha\geq 0\}$. \protect\subref{subfig:grad-g-parallel-to-line} Level curves of a dual function $g~{:\R^3\to\R}$ and two rays with the FGOR characteristic. \protect\subref{subfig:grad-g-not-parallel-to-line-plane-view} Corresponding plan view.}  
\label{Fig:example-FGOR-dual-grads-parellel-to-line}
\end{figure}
Now it is straightforward to see that under condition~\eqref{eq:condition-dual-function}, with constant stepsizes, the dual subgradient method~\eqref{eq:Lambda-Update} initialized at a point $\boldsymbol\lambda^{(0)}$ in the ray $\{\boldsymbol{\lambda}+\alpha\boldsymbol{\mu} \ | \ \alpha\geq 0\}$ continues to remain in the same ray yielding a constant increase of the dual function $g$ at each iteration, until $\boldsymbol\lambda^{(k)}\notin \{\boldsymbol{\lambda}+\alpha\boldsymbol{\mu} \ | \ \alpha\geq 0\}$ for some iteration index $k$. When $\boldsymbol\lambda^{(k)}\notin \{\boldsymbol{\lambda}+\alpha\boldsymbol{\mu} \ | \ \alpha\geq 0\}$, it is natural to switch to any other existing stepsize rule. More specifically, we lay out the following stepsize rule:
\begin{align}\label{eq:stepsize-rule-original-1}
  &\mathrm{If} \ \boldsymbol\lambda^{(k)}\in \{\boldsymbol{\lambda}+\alpha\boldsymbol{\mu} \ | \ \alpha\geq 0\}, \ \mathrm{chose} \ \gamma_k = \gamma, \\ \label{eq:stepsize-rule-original-2}
  &\mathrm{Otherwise,} \  \mathrm{switch \ to \ an \ existing \ stepsize \ rule}.
\end{align}

\subsection{\addb{Determining a Ray Ensuring (18)}}

\begin{figure}[!t]
\centering
\subfloat[]{\includegraphics[width=0.50\linewidth]{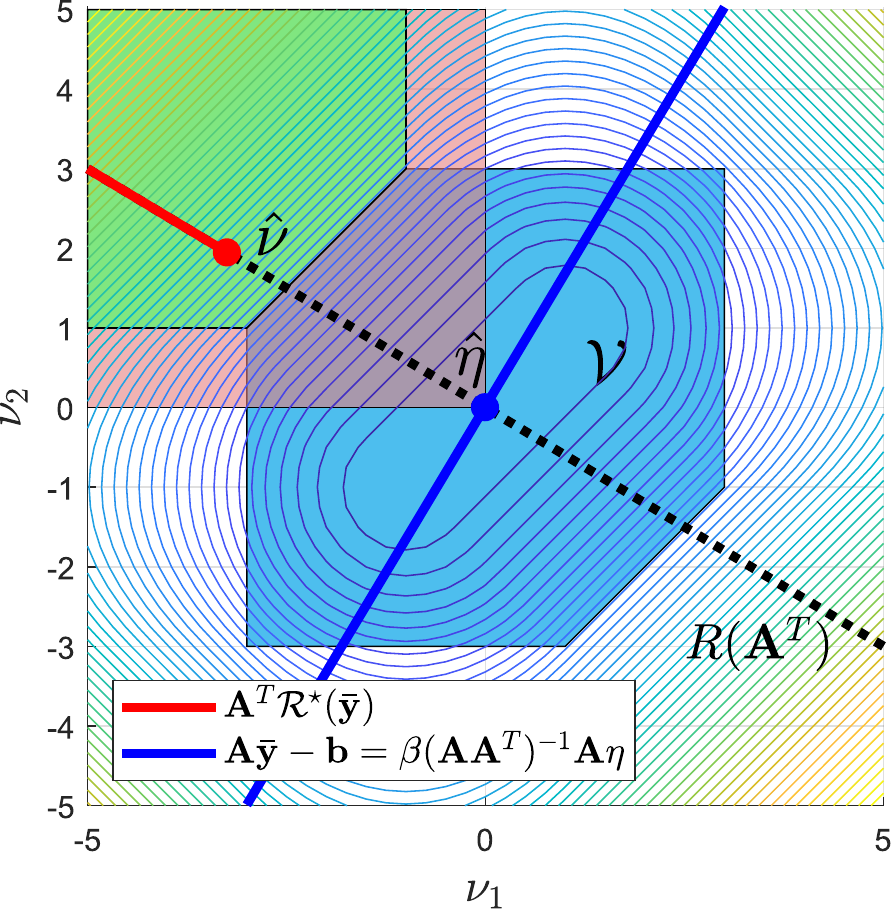}%
\label{subfig:feasible-problem1-1}}
 \hfil
 \subfloat[]{\includegraphics[width=0.50\linewidth]{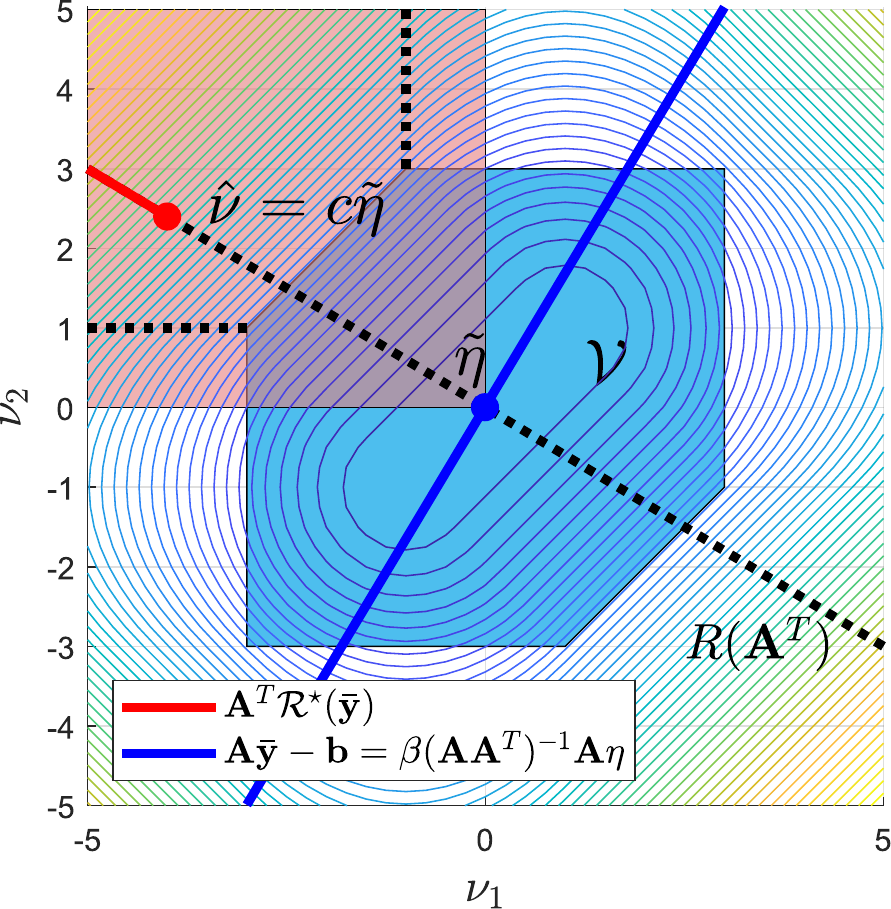}%
 \label{subfig:feasible-problem1-2}}
\caption{An illustration of a ray ${\mathcal{R}}^{\star}(\mathbf{\bar y})$ [\cf~\eqref{eq:FGOR-Ray}] associated to problem~\eqref{eq:Ex:Assumption:dual-flat-assump}, where $\mathbf{A}=[-0.05 \ 0.03]$, $\mathbf{b}=\mathbf{0}$, $\mathbf{\bar y}=[-2 \ 2]^\text{T}$, and $\beta=1$. The sets $N_{\mathcal{Y}_2}(\mathbf{\bar y})$ and $\partial f_2(\mathbf{\bar y})+N_{\mathcal{Y}_2}(\mathbf{\bar y})$ are depicted by the regions shaded in pink and green, respectively. \protect\subref{subfig:feasible-problem1-1} ${\mathcal{R}}^{\star}(\mathbf{\bar y})$ given by solving \eqref{eq:find-ray-direction-eta---}: $\hat{\boldsymbol{\nu}}=[-3.25 \ 1.95]^\text{T}$, $\hat{\boldsymbol{\eta}}=[-0.008 \ 0.0048]^\text{T}$, $\hat{\boldsymbol{\lambda}}=65$, $\hat{\boldsymbol{\mu}}=0.16$.  \protect\subref{subfig:feasible-problem1-2} ${\mathcal{R}}^{\star}(\mathbf{\bar y})$ given by solving \eqref{eq:find-ray-direction-eta-2}: $\boldsymbol{\tilde\eta}=[-0.008 \ 0.0048]^\text{T}$,~$\boldsymbol{\tilde\mu}=0.16$. Then $\hat{\boldsymbol{\nu}}=c\boldsymbol{\tilde\eta}$, $\hat{\boldsymbol{\eta}}=\boldsymbol{\tilde\eta}$, $\hat{\boldsymbol{\lambda}}=c\boldsymbol{\tilde\mu}$, $\hat{\boldsymbol{\mu}}=\boldsymbol{\tilde\mu}$, where $c=500$, \cf~\eqref{eq:Solution-for-Relaxed-Problem}.}   
\label{Fig:feasible-problem}
\end{figure}

\addb{The computation of a ray in the form $\{\boldsymbol{\lambda}+\alpha\boldsymbol{\mu} \ | \ \alpha\geq 0\}$, over which condition ~\eqref{eq:condition-dual-function} holds, is carried out through a properly formulated \emph{feasibility problem}, as detailed in the sequel. %
Let us first note the following equivalent conditions for \eqref{eq:condition-dual-function}:
\begin{align} \nonumber
   & \forall \alpha\geq 0,~\nabla  g(\boldsymbol{\lambda}  +\alpha \boldsymbol{\mu} )=-\beta {\boldsymbol\mu}  \\ \label{eq:equivalence-gradient-parallel-assumption-1}
    & \iff \forall \alpha\geq 0,~- \mathbf{A} \nabla f^*\big(\mathbf{A}\tran(\boldsymbol{\lambda}+\alpha \boldsymbol{\mu})\big)=-\beta{\boldsymbol\mu}-\mathbf{b} \\ \label{eq:equivalence-gradient-parallel-assumption-2}
    & \iff \forall \alpha\geq 0,~- \mathbf{A} \nabla f^*(\boldsymbol{\nu}+\alpha \boldsymbol{\eta})=-\beta(\mathbf{A}\mathbf{A}\tran)^{-1}\mathbf{A}{\boldsymbol\eta}-\mathbf{b},
\end{align}
where \eqref{eq:equivalence-gradient-parallel-assumption-1} follows from \eqref{eq:dual-grad-and-conjugate-grad} and \eqref{eq:equivalence-gradient-parallel-assumption-2} follows from that $\boldsymbol{\lambda}$ and $\boldsymbol{\mu}$ are inverse images of some $\boldsymbol{\nu}$ and $\boldsymbol{\eta}$ under the linear mapping $\mathbf{A}\tran$ [\cf~\eqref{eq:Inverse-Image-of-Conjugate-Ray-1}]. Consequently, the feasibility of the optimization problem
\begin{equation} \label{eq:feasibility-prob-1}
\begin{array}{ll}
\mbox{minimize} & 0\\
\mbox{subject to} & \forall \alpha\geq 0,~- \mathbf{A} \nabla f^*(\boldsymbol{\nu}+\alpha \boldsymbol{\eta})\\
& \hspace{6em} =-\beta(\mathbf{A}\mathbf{A}\tran)^{-1}\mathbf{A}{\boldsymbol\eta}-\mathbf{b}\\
& \mathbf{A}\tran\boldsymbol{\lambda}=\boldsymbol{\nu}\\
& \mathbf{A}\tran \boldsymbol{\mu}= \boldsymbol{\eta}
\end{array}
\end{equation}
with variables ${\boldsymbol\lambda}$, ${\boldsymbol\mu}$, ${\boldsymbol\nu}$, and ${\boldsymbol\eta}$, is \emph{sufficient} for \eqref{eq:condition-dual-function}.

At this point, we rely on Proposition~\ref{Prop:flat-region} to yield a more constrained formulation of \eqref{eq:feasibility-prob-1}. More specifically, Proposition~\ref{Prop:flat-region} claims that there are rays in $\R^n\setminus\texttt{int}~\mathcal{V}$ over which the gradient $\nabla f^*$ is constant and is characterized by some $\bar{\mathbf{y}}\in\texttt{bnd}~\mathcal{Y}$ [\cf~\eqref{eq:FGOR-of}]. Thus, in line with this claim, let us impose an additional constraint $\forall~\alpha\geq 0, \ \bar{\mathbf{y}}= \nabla f^*(\boldsymbol{\nu}+\alpha \boldsymbol{\eta})$, %
where $\bar{\mathbf{y}}\in\texttt{bnd}~\mathcal{Y}$. As a result, we have the following problem, which is more constrained than \eqref{eq:feasibility-prob-1}:
\begin{equation} \label{eq:feasibility-prob-2}
\begin{array}{ll}
\mbox{minimize} & 0\\
\mbox{subject to} & - \mathbf{A} \bar{\boldsymbol{y}} =-\beta(\mathbf{A}\mathbf{A}\tran)^{-1}\mathbf{A}{\boldsymbol\eta}-\mathbf{b}\\
& \forall \alpha\geq 0,~\nabla f^*(\boldsymbol{\nu}+\alpha \boldsymbol{\eta}) = \bar{\boldsymbol{y}}\\
& \mathbf{A}\tran\boldsymbol{\lambda}=\boldsymbol{\nu}\\
& \mathbf{A}\tran \boldsymbol{\mu}= \boldsymbol{\eta},
\end{array}
\end{equation}
where the variables are ${\boldsymbol\lambda}$, ${\boldsymbol\mu}$, ${\boldsymbol\nu}$, and ${\boldsymbol\eta}$. By using \cite[Prop.~11.3]{Rockafellar-98}, \eqref{eq:feasibility-prob-2} can equivalently be reformulated as 
\begin{equation} \label{eq:feasibility-prob-3}
\begin{array}{ll}
\mbox{minimize} & 0\\
\mbox{subject to} & - \mathbf{A} \bar{\boldsymbol{y}} =-\beta(\mathbf{A}\mathbf{A}\tran)^{-1}\mathbf{A}{\boldsymbol\eta}-\mathbf{b}\\
& \{\boldsymbol{\nu}+\alpha \boldsymbol{\eta} \ | \ \alpha\geq 0 \}\subseteq\partial f(\bar{\mathbf{y}})\\
& \mathbf{A}\tran\boldsymbol{\lambda}=\boldsymbol{\nu}\\
& \mathbf{A}\tran \boldsymbol{\mu}= \boldsymbol{\eta},
\end{array}
\end{equation}
where the second constraint of \eqref{eq:feasibility-prob-2} is updated accordingly.
}
\addb{It is straightforward to see that the second constraint in \eqref{eq:feasibility-prob-3} can break into two constraints ${\boldsymbol\nu} \in \partial f_0(\bar{\mathbf{y}}) + N_{\mathcal{Y}}(\bar{\mathbf{y}})$ and ${\boldsymbol\eta} \in N_{\mathcal{Y}}(\bar{\mathbf{y}})$ [\cf proof of Proposition~\ref{Prop:flat-region}, \eqref{eq:FGOR-proof-1}, \eqref{eq:prop-FGOR-of-conjugate-2}]. Thus, problem~\eqref{eq:feasibility-prob-3} is equivalent to
\begin{equation} \label{eq:find-ray-direction-eta---}
\begin{array}{ll}
\mbox{minimize} & 0\\
\mbox{subject to} & \mathbf{A}\bar{\mathbf{y}}- \mathbf{b} = \beta(\mathbf{A}\mathbf{A}\tran)^{-1}\mathbf{A}{\boldsymbol\eta}\\
& {\boldsymbol\nu} \in \partial f_0(\bar{\mathbf{y}}) + N_{\mathcal{Y}}(\bar{\mathbf{y}})\\
& {\boldsymbol\eta} \in N_{\mathcal{Y}}(\bar{\mathbf{y}})\\
& \mathbf{A}\tran \boldsymbol{\lambda} =  {\boldsymbol\nu} \\
& \mathbf{A}\tran \boldsymbol{\mu} =  {\boldsymbol\eta},
\end{array}
\end{equation}
where the variables are ${\boldsymbol\nu}$, ${\boldsymbol\eta}$, ${\boldsymbol\lambda}$, and ${\boldsymbol\mu}$. 
It is worth noting that
\begin{equation}\nonumber
    \mbox{feasibility of }~\eqref{eq:find-ray-direction-eta---}\implies \mbox{feasibility of}~\eqref{eq:feasibility-prob-1}\implies \eqref{eq:condition-dual-function}.
\end{equation}
If problem~\eqref{eq:find-ray-direction-eta---} is solved with solution $(\hat{\boldsymbol\nu},\hat{\boldsymbol\eta},\hat{\boldsymbol\lambda},\hat{\boldsymbol\mu})$, the ray we seek is simply given~by 
\begin{equation}\label{eq:FGOR-Ray}
    {\mathcal{R}}^{\star}(\bar{\mathbf{y}})=\{\hat{\boldsymbol{\lambda}}+\alpha \hat{\boldsymbol{\mu}} \ | \ \alpha\geq 0 \}.
\end{equation}
An illustration is given in \figurename~\ref{Fig:feasible-problem}\subref{subfig:feasible-problem1-1}. If problem~\eqref{eq:find-ray-direction-eta---} is infeasible, one may choose a different $\bar{\mathbf{y}}\in\texttt{bnd}~\mathcal{Y}$.}

In practice, one should have available the knowledge of $ \partial f_0(\bar{\mathbf{y}})$ to solve problem~\eqref{eq:find-ray-direction-eta---}. Note that in a distributed optimization setting, the global objective function~$f_0$ is often based on local objective functions of the involved subsystems. Consequently, $\partial f_0(\bar{\mathbf{y}})$, unlike other problem data, cannot be determined based on only the subsystems' local information. Therefore, this issue has to be dealt with through coordination among the subsystems.
If such coordination is to be avoided, an alternative is to relax problem~\eqref{eq:find-ray-direction-eta---} by dropping the second constraint that entails the main coupling among subsystems. Consequently, the fourth constraint $\mathbf{A}\tran \boldsymbol{\lambda} =  {\boldsymbol\nu}$ becomes obsolete and the resulting relaxed formulation is given by
\begin{equation} \label{eq:find-ray-direction-eta-2}
\begin{array}{ll}
\mbox{minimize} & 0\\
\mbox{subject to} & \mathbf{A}\bar{\mathbf{y}}- \mathbf{b} = \beta(\mathbf{A}\mathbf{A}\tran)^{-1}\mathbf{A}{\boldsymbol\eta}\\
& {\boldsymbol\eta} \in N_{\mathcal{Y}}(\bar{\mathbf{y}})\\
& \mathbf{A}\tran \boldsymbol{\mu} =  {\boldsymbol\eta},
\end{array}
\end{equation}
where the variables are ${\boldsymbol\eta}$ and ${\boldsymbol\mu}$. A careful examination of problem~\eqref{eq:find-ray-direction-eta-2} shows that it can be solved in two steps: 1) solve the equality constraints to uniquely yield the solutions $\tilde{\boldsymbol\mu}=(1/\beta)(\mathbf{A}\mathbf{\bar y}-\mathbf{b})$ and $\tilde{\boldsymbol\eta}=\mathbf{A}\tran\tilde{\boldsymbol\mu}$, 2) check whether the condition $\tilde{\boldsymbol\eta} \in N_{\mathcal{Y}}(\bar{\mathbf{y}})$ holds. If the condition is affirmative, then the unique solution of problem~~\eqref{eq:find-ray-direction-eta-2} is $(\tilde{\boldsymbol\mu},\tilde{\boldsymbol\eta})$. Otherwise, a different $\bar{\mathbf{y}}\in\texttt{bnd}~\mathcal{Y}$ may be chosen and the process is~repeated. \addb{We note that providing a point $\mathbf{\bar y}\in\texttt{bnd}~\mathcal{Y}$ such that problems~\eqref{eq:find-ray-direction-eta---} or \eqref{eq:find-ray-direction-eta-2} are feasible depends on the specific problem instance. In particular, one has to carefully examine the structure of the problem in the primal domain to identify such a point. For example, in Section~\ref{subsec:Consensus-problem-ell-infinity-norm-constraints}, we provide a point $\mathbf{\bar y}\in\texttt{bnd}~\mathcal{Y}$, such that problem~\eqref{eq:find-ray-direction-eta-2} associated to the global consensus problem, is feasible, \cf~Lemma~\ref{prop:consensus-problem-feasible-y-bar}. }

\addb{Note that if problem~\eqref{eq:find-ray-direction-eta---} is feasible, then so is problem~\eqref{eq:find-ray-direction-eta-2}, which follows trivially. Conversely, if \eqref{eq:find-ray-direction-eta-2} is feasible, then so is \eqref{eq:find-ray-direction-eta---} if the affine dimension of the normal cone $N_{\mathcal{Y}}(\mathbf{\bar y})$ is $n$. Specifically, if $(\boldsymbol{\tilde \mu},\boldsymbol{\tilde \eta})$ is feasible for \eqref{eq:find-ray-direction-eta-2}, then 
\begin{equation}\label{eq:Solution-for-Relaxed-Problem}
\hspace{-2mm}\hat{\boldsymbol{\mu}}=\boldsymbol{\tilde\mu},\ \hat{\boldsymbol{\eta}}=\boldsymbol{\tilde\eta},\ \hat{\boldsymbol{\nu}}=c\boldsymbol{\tilde\eta},\  \hat{\boldsymbol{\lambda}}=c(\mathbf{AA}\tran)^{-1}\mathbf{A}\boldsymbol{\tilde\eta}=c\boldsymbol{\tilde{\mu}}
\end{equation}
are feasible for~\eqref{eq:find-ray-direction-eta---} for sufficiently large $c\in\R$, where the last equality follows from \eqref{eq:Inverse-Image-of-Conjugate-Ray-1}. To see this, first note that the constraints $1$, $3$, and $5$ of~\eqref{eq:find-ray-direction-eta---} are trivially satisfied by $\hat{\boldsymbol{\mu}},\hat{\boldsymbol{\eta}}$. Second, we note that the constraint $4$, $\mathbf{A}\tran \hat{\boldsymbol{\lambda}} =  \hat{\boldsymbol\nu}$, is also straightforwardly satisfied. Finally, the constraint $2$ is also satisfied from that $\hat{\boldsymbol{\nu}}=c\boldsymbol{\tilde\eta}$ is in the set $\mathbf{d}_0+N_{\mathcal{Y}}(\mathbf{\bar y})$, when $c$ is sufficiently large, where $\mathbf{d}_0\in\R^n$ is arbitrary. In other words, if the condition $\boldsymbol{\tilde\eta}\in N_{\mathcal{Y}}(\mathbf{\bar y})$ holds [\cf~constraint $2$ of \eqref{eq:find-ray-direction-eta-2}], then for any translation of the cone $N_{\mathcal{Y}}(\mathbf{\bar y})$, there exists a scaling of $\boldsymbol{\tilde\eta}$ such that the scaled 
$\boldsymbol{\tilde\eta}$ remains within the translated cone [\cf~constraint $2$ of \eqref{eq:find-ray-direction-eta---}]. Consequently, by using the solution $(\boldsymbol{\tilde\mu},\boldsymbol{\tilde\eta})$ obtained by solving \eqref{eq:find-ray-direction-eta-2}, we can compute the required RFGOR ray $\mathcal{R}^{\star}(\mathbf{\bar y})$ [\cf~\eqref{eq:FGOR-Ray}] by replacing $\hat{\boldsymbol{\lambda}}$ and $\hat{\boldsymbol{\mu}}$ from \eqref{eq:Solution-for-Relaxed-Problem}, \cf~\figurename~\ref{Fig:feasible-problem}\subref{subfig:feasible-problem1-2}.} 

\addb{Note that problem~\eqref{eq:find-ray-direction-eta-2} can be solved independently by each subsystem, unlike \eqref{eq:find-ray-direction-eta---}. This is because it is quite common in practice that the parameters $\mathbf A$ and $\mathcal{Y}$ of problem~\eqref{eq:optimization_prob_main} are known by subsystems as far as distributed optimization settings are concerned (e.g., consensus problem). Moreover, it is worth noting that both problems~\eqref{eq:find-ray-direction-eta-2} and \eqref{eq:find-ray-direction-eta---} do not require knowing in advance whether Assumption~\ref{Assumption:dual-flat-assump} is satisfied. To see this, let us first remark the following:
\begin{remark}\label{Remark:Assumption-3}
The feasibility of problem \eqref{eq:find-ray-direction-eta-2} [or \eqref{eq:find-ray-direction-eta---}] confirms the satisfaction of Assumption~\ref{Assumption:dual-flat-assump}. 
\end{remark}
\begin{IEEEproof}
The feasibility of \eqref{eq:find-ray-direction-eta-2} [or \eqref{eq:find-ray-direction-eta---}] ensures the existence of an RFGOR ray $\mathcal{R}^\star(\mathbf{\bar y})\subseteq \texttt{dom}~g$, \cf~\eqref{eq:FGOR-Ray},\eqref{eq:Solution-for-Relaxed-Problem}. Consequently, $\nabla f^*$, when evaluated over the ray $\mathbf{A}\tran \mathcal{R}^\star(\mathbf{\bar y})$ [which is trivially in ${R}(\mathbf{A}\tran)$], becomes a constant vector [\cf~\eqref{eq:dual-grad-and-conjugate-grad}], yielding the FGOR characteristic. 
\end{IEEEproof}
As such, even without explicitly knowing whether Assumption~\ref{Assumption:dual-flat-assump} holds, one can simply solve problem~\eqref{eq:find-ray-direction-eta-2} and compute an RFGOR ray $\mathcal{R}^\star(\mathbf{\bar y})$, if feasible. Then we provide the initialization point $\boldsymbol{\lambda}^{(0)}$ from $\mathcal{R}^\star(\mathbf{\bar y})$ for the dual subgradient method. Let us next discuss the corresponding algorithm.}

\subsection{Algorithm}
\addb{We solve \eqref{eq:find-ray-direction-eta-2} and choose the initialization point $\boldsymbol{\lambda}^{(0)}$ for the dual subgradient method~\eqref{eq:Lambda-Update} as follows:
\begin{equation}\label{eq:initialization-point-subgrad-method} \boldsymbol{\lambda}^{(0)}=\hat{\boldsymbol{\lambda}}=c\boldsymbol{\tilde\mu},
\end{equation}
where $c\in\R$, \cf~\eqref{eq:Solution-for-Relaxed-Problem}. Then the corresponding algorithm, together with the proposed stepsize rule [\cf~\eqref{eq:stepsize-rule-original-1}, \eqref{eq:stepsize-rule-original-2}], is outlined in Algorithm~\ref{Alg:Dual-Subgradient-Method}. Note that the proposed $\boldsymbol{\lambda}^{(0)}$ [\cf~\eqref{eq:initialization-point-subgrad-method}] serves as a warm start for the dual subgradient method~\eqref{eq:Lambda-Update}, as it enables \eqref{eq:Lambda-Update} to use a constant stepsize and produce constant increments in the dual function $g$ as long as $\boldsymbol{\lambda}^{(k)}\in\mathcal{R}^{\star}(\mathbf{\bar y})$ [\cf~steps~2-3] allowing $\boldsymbol{\lambda}^{(k)}$ to quickly approach the dual optimal solution $\boldsymbol{\lambda}^{\star}$. 

}
\begin{algorithm}[t]
	\caption{Dual Subgradient Method with FGOR} 
	\begin{algorithmic}[1]
            \Require the solution $\tilde{\boldsymbol{\mu}}$ of problem~\eqref{eq:find-ray-direction-eta-2}, $c >0$, and $\beta>0$.  
	    \State $\boldsymbol{\lambda}^{(0)}=c\tilde{\boldsymbol{\mu}}$, $k=0$.
        \If{$\mathbf{A}\tran\boldsymbol{\lambda}^{(k)}{\in}\partial f_0(\bar{\mathbf{y}}) + N_{\mathcal{Y}}(\bar{\mathbf{y}})$} \Comment{in RFGOR}

                \State Perform \eqref{eq:Lambda-Update} with $\gamma_k=\gamma$ and $\mathbf{s}^{(k)}=\beta\boldsymbol{\tilde\mu}$.
                \State Set $k:=k+1$ and go to step~2.
            \Else \Comment{RFGOR is exited}
                \State Go to step~8. \Comment{to an existing stepsize rule}
            \EndIf 
            \Repeat 
            \State Perform \eqref{eq:Lambda-Update} with an existing stepsize rule.
            \State $k:=k+1$.
        \Until{a stopping criterion true}
	\end{algorithmic} 
	\label{Alg:Dual-Subgradient-Method}
\end{algorithm} 

Note that the condition to be checked at step~2 of Algorithm~\ref{Alg:Dual-Subgradient-Method} is unfavorable in a distributed optimization setting. Therefore, we may wish to have at our disposal an alternative, yet equivalent representation that is favorable in a distributed setting. This is the subject of the following remark.
\begin{remark}\label{rem:Equivalent-FGOR-for-Distribution}
Let $\mathbf{y}^{(k)}$ be the minimizer of the Lagrangian associated to problem~\eqref{eq:optimization_prob_main_extended} with Lagrange multiplier $\boldsymbol{\lambda}^{(k)}$, i.e., $\mathbf{y}^{(k)}=   \arg\min_{\mathbf{y}\in \mathcal{Y}} \ f_0(\mathbf{y}) - {\boldsymbol{\lambda}}^{(k)\mbox{\scriptsize T}} (\mathbf{A}\mathbf{y}-\mathbf{b})$. Then 
\begin{equation}\label{eq:rem:Equivalent-FGOR-for-Distribution-1}
\mathbf{A}\tran\boldsymbol{\lambda}^{(k)}\in\partial f_0(\bar{\mathbf{y}}) + N_{\mathcal{Y}}(\bar{\mathbf{y}}) \iff \mathbf{y}^{(k)}=\bar{\mathbf{y}} 
\end{equation}
\end{remark}
\begin{IEEEproof} 
\addb{See Appendix~\ref{App:Equivalent-FGOR-for-Distribution}.}
\end{IEEEproof}

\addb{We establish next the convergence of Algorithm~\ref{Alg:Dual-Subgradient-Method}.
\begin{remark}\label{rem:conv-of-the-algorithm}
Suppose Assumption~\ref{Assumption:PropStCvx}, Assumption~\ref{Assumption:Exclusion-of-Some-Regions}, and Assumption~\ref{Assumption:dual-flat-assump} hold. Then, Algorithm~\ref{Alg:Dual-Subgradient-Method} converges to an optimal solution $\boldsymbol{\lambda}^{\star}$ of the dual problem~\eqref{eq:dual-problem}.

\end{remark}
\begin{IEEEproof}
If $\nabla g(\boldsymbol{\lambda}^{(0)})=\mathbf{0}$, then $\boldsymbol{\lambda}^{(0)}$ is an optimal solution of problem~\eqref{eq:dual-problem}, and hence the convergence is trivial. Otherwise, either
1) $\boldsymbol{\lambda}^{(0)}\notin \textrm{RFGOR}$ [\cf~step~5] or 2) $\boldsymbol{\lambda}^{(0)}\in\textrm{RFGOR}$ [\cf~step~2]. In the first case, the control goes to step~6, after which the algorithm performs the standard subgradient method~\eqref{eq:Lambda-Update} with an existing stepsize rule [\cf~steps~8-11], guaranteeing the convergence. In the latter case, Algorithm~\ref{Alg:Dual-Subgradient-Method} performs the subgradient method~\eqref{eq:Lambda-Update} with a constant stepsize~$\gamma$ and with the gradient $\beta\boldsymbol{\tilde \mu}$. %
In this respect, we note that there is a support to $g$ at all $\boldsymbol{\lambda}^{(k)}\in\textrm{RFGOR}$, denoted~$h$, which is given by the affine function $h(\boldsymbol{\lambda})=-\beta\boldsymbol{\tilde \mu}\tran\boldsymbol{\lambda}+c$ for some $c\in\R$, which doesn't depend on $\boldsymbol{\lambda}^{(k)}$. Thus, at each $\boldsymbol{\lambda}^{(k)}\in\textrm{RFGOR}$, $g(\boldsymbol{\lambda}^{(k)})=h(\boldsymbol{\lambda}^{(k)})$. Moreover, since $g$ is concave and has an optimal solution $\boldsymbol{\lambda}^{\star}$~\footnote{\addb{According to the assumptions, primal problem~\eqref{eq:optimization_prob_main} has a unique solution~${\mathbf{y}}^\star$ and strong duality holds.}}, $h$ is an upper-bound on $g$ such that $h(\boldsymbol{\lambda})>g(\boldsymbol{\lambda})$ for some $\boldsymbol{\lambda}$. 
Thus, there must be an index $\bar k>0$ such that $h(\boldsymbol{\lambda}^{(k)})=g(\boldsymbol{\lambda}^{(k)})$ for all $k\leq \bar k$ and $h(\boldsymbol{\lambda}^{(\bar k +1)})>g(\boldsymbol{\lambda}^{(\bar k+1)})$. This implies that $\boldsymbol{\lambda}^{(\bar k+1)}\not\in\textrm{RFGOR}$, and thus the algorithm switches to the standard subgradient method, guaranteeing the convergence according to steps~5,6,8, and 9.
\end{IEEEproof}
}

In the sequel, we place greater emphasis on the application of FGOR to the global consensus problem~\cite[\S~7]{Boyd-Parikh-Chu-Peleato-Eckstein-2010} since it is a widely used formulation in numerous application domains~\cite{Chen_consensus_2012,Carlo-Machine-Learning,Yang-survey-of-distributed-optimization-2019,Ball_consensus_2023,halsted-Multiple-robots-2021}.

\section{Global Consensus Problem}\label{sec:Application:The-General-Consensus-Problem}

Let us now consider the problem
\begin{equation} \label{eq:consensus-problem}
\begin{array}{ll}
\mbox{minimize} & (1/m)\sum_{i{=}1}^{m}f_i(\mathbf{z}) \\
\mbox{subject to} & \mathbf{z}\in\mathcal{Z}, 
\end{array}
\end{equation}
where the variable is $\mathbf{z}\in\R^n$. We assume that the functions $f_i:\R^{n}\rightarrow \R$, $i=1,\ldots,m$ and the constraint set $\mathcal{Z}\subseteq\R^n$ are conforming to Assumption~\ref{Assumption:PropStCvx} and Assumption~\ref{Assumption:Exclusion-of-Some-Regions}. %
A commonly used distributed solution method for solving \eqref{eq:consensus-problem} is based on dual decomposition~\cite{Boyd-EE364b-PrimDualDecomp-07}.

\subsection{The Standard Dual Decomposition Method}

Problem~\eqref{eq:consensus-problem} is equivalent to the global consensus problem
\begin{equation} \label{eq:consensus-problem-distributed}
\begin{array}{ll}
\mbox{minimize} & f_0(\mathbf{y})=(1/m)\sum_{i{=}1}^{m}f_i(\mathbf{y}_{i}) \\
\mbox{subject to} & \mathbf{y}_i\in\mathcal{Z},\ i\in \mathcal{S}\\ 
& \mathbf{y}_i=\mathbf{y}_{i+1}, i\in\mathcal{S}\setminus \{m\},
\end{array}
\end{equation}
where $\mathbf{y}_i\in\R^n$, $i\in\mathcal{S}$, are the local versions of $\mathbf{z}$, $\mathbf{y}=[\mathbf{y}_{1}\tran \ \ldots \ \mathbf{y}_{m}\tran]\tran$, and $\mathcal{S}=\{1,\ldots,m\}$ is the set of agents. The equality constraint $\mathbf{y}_i=\mathbf{y}_{i+1}$, $i\in\mathcal{S}\setminus \{m\}$ is introduced to impose the consistency among the local variables $\mathbf{y}_i$s. Then the dual function $g:\R^{n(m-1)}\rightarrow \overline{\R}$ of problem~\eqref{eq:consensus-problem-distributed} is given~by 
\begin{align} \label{eq:consensus-problem-dual-function}
    g(\boldsymbol{\lambda)}& = {\underset{\mathbf{y}_{i}\in\mathcal{Z},\ i\in\mathcal{S}}{\inf}}\left[\sum_{i=1}^{m}f_i(\mathbf{y}_{i}) {+} \sum_{i=1}^{m{-}1} \boldsymbol{\lambda}_i\tran(\mathbf{y}_{i}{-}\mathbf{y}_{{i+1}})\right],\allowdisplaybreaks 
    \end{align}
where $\boldsymbol{\lambda}_{i}\in\R^n$ denotes the Lagrange multiplier associated to the constraint $\mathbf{y}_{i}=\mathbf{y}_{{i+1}}$, $i\in\mathcal{S}\setminus \{m\}$ and $\boldsymbol{\lambda}=[\boldsymbol{\lambda}_1\tran \ \ldots \ \boldsymbol{\lambda}_{m-1}\tran]\tran$. %
Moreover, the dual problem associated to \eqref{eq:consensus-problem-distributed} is given by
\begin{equation}\label{eq:dual-problem-consensus-problem}
\underset{\boldsymbol{\lambda}\in\R^{n(m-1)}}{\text{maximize}} \quad g(\boldsymbol{\lambda}).
\end{equation}
The standard dual decomposition algorithm for solving \eqref{eq:dual-problem-consensus-problem} is given by Algorithm~\ref{Alg:Dual-Decomposition-Algorithm}.
\begin{algorithm}[t]
	\caption{Dual Decomposition for Global Consensus} 
	\begin{algorithmic}[1]
	    \Require $\boldsymbol{\lambda}\in\R^{n(m-1)}$.
	    \State $k=0$, $\boldsymbol{\lambda}^{(0)}=\boldsymbol{\lambda}$. $\boldsymbol{\lambda}_0^{(j)}=\boldsymbol{\lambda}_m^{(j)}=\mathbf{0}\in\R^n, \ j\in\mathbb{Z}^0_+$.
	    \State Central server (CS) broadcasts $\boldsymbol{\lambda}^{(0)}$ to agents. 
		\Repeat 
            \State $\forall i$, agent $i$ computes:
            \begin{equation}\nonumber
            \mathbf{y}_i^{(k)}=\underset{\mathbf{y}_i\in\mathcal{Z}}{\arg\min}~f_i(\mathbf{y}_i)+{\big(\boldsymbol{\lambda}^{(k)}_i{-}\boldsymbol{\lambda}^{(k)}_{i-1}\big)}\tran \mathbf{y}_i.
             \end{equation}
 \State $\forall i$, agent $i$ transmits $\mathbf{y}_i^{(k)}$ to CS. %
 
            \State CS computes ${\mathbf{s}}^{(k)}{=}[(\mathbf{y}^{(k)}_2{-}{\mathbf{y}}^{(k)}_1)\tran{\ldots}(\mathbf{y}^{(k)}_{m}{-}\mathbf{y}^{(k)}_{m{-}1})\tran]\tran$. 
            \State CS computes $\boldsymbol{\lambda}^{(k+1)} = \boldsymbol{\lambda}^{(k)}- \gamma_k\mathbf{s}^{(k)}$. 
             \State $\forall i$, CS broadcasts $\boldsymbol{\lambda}_{i-1}^{(k+1)}$ and $\boldsymbol{\lambda}_i^{(k+1)}$ to agent $i$. %
            \State $k:=k+1$.
        \Until{a stopping criterion true}
	\end{algorithmic} 
	\label{Alg:Dual-Decomposition-Algorithm}
\end{algorithm}

In Algorithm~\ref{Alg:Dual-Decomposition-Algorithm}, CS performs the dual variable update [\cf~step~7] using any standard stepsize rule \cite{Polyak_Intro_Opt,goldstein_1962,armijo_1966,beck_2009,bello_2016,Asl_2020,khirirat_2023,Thinh_2018,Thinh_2021,zhor_1968_1,zhor_1968_2,goffin_1977,shor_2012,Damek_2018,zhu_2019,aybat_2019,wang_2021,chen_2019,krizhevsky_2012,He_2016,ge_2019,polyak_1969,hazan_2022,loizou_2021,wang_2023,barzilai_1988,dai_2002,burdakov_2019,Yura_2020}. It is worth noting that when the FGOR characteristic of $g$ is utilized, the communication steps~5 and 8 in Algorithm~\ref{Alg:Dual-Decomposition-Algorithm} are not needed as long as $\boldsymbol{\lambda}^{(k)}$ resides in an RFGOR region. In particular, utilizing the FGOR characteristic of $g$, we can improve the performance of Algorithm~\ref{Alg:Dual-Decomposition-Algorithm} not only in terms of the speed of the convergence but also in terms of communication efficiency. We discuss this idea more formally in the following subsection.

\subsection{Using FGOR Characteristic of $g$ on Algorithm~\ref{Alg:Dual-Decomposition-Algorithm}} \label{subsec:using-FGOR-characteristics-on-DDA}
For our exposition of using the FGOR characteristic of $g$ on Algorithm~\ref{Alg:Dual-Decomposition-Algorithm}, we start by noting that the constraints $\mathbf{y}_i\in\mathcal{Z},\ i\in \mathcal{S}$ and $\mathbf{y}_i=\mathbf{y}_{i+1}, i\in\mathcal{S}\setminus \{m\}$ are equivalent with $\mathbf{y}\in\mathcal{Z}^m$ and $\mathbf{Ay}=\mathbf{0}$, respectively, where $\mathbf{y}=[\mathbf{y}_{1}\tran \ \ldots \ \mathbf{y}_{m}\tran]\tran$ and $\mathcal{Z}^m$ is the $m$-fold Cartesian product of $\mathcal{Z}$. The matrix $\mathbf{A}\in\R^{n(m-1)\times nm}$ has the form
\begin{equation}\label{eq:distributed-consenses-problem-matrix-A}
\mathbf{A}=\begin{bmatrix} 
\mathbf{I}_n & -\mathbf{I}_{n} & \mathbf{0} & \cdots & \mathbf{0}\\  
\mathbf{0} & \mathbf{I}_n & -\mathbf{I}_{n} & \ddots & \vdots \\ %
\vdots & \ddots & \ddots & \ddots &  \mathbf{0}\\ %
\mathbf{0} & \cdots & \mathbf{0} & \mathbf{I}_n & -\mathbf{I}_n
\end{bmatrix}  .
\end{equation} 
We let $\mathcal{Z}^m=\mathcal{Y}$ for clarity. Suppose that problem~\eqref{eq:find-ray-direction-eta-2} associated to problem~\eqref{eq:consensus-problem-distributed} is feasible for some $\mathbf{\bar y}=[\mathbf{\bar y}_{1}\tran \ \ldots \ \mathbf{\bar y}_{m}\tran]\tran\in\texttt{bnd}~\mathcal{Y}$, where $\mathbf{\bar y}_{i}\in\mathcal{Z}$, $\forall i\in\mathcal{S}$~\addb{\footnote{\addb{The feasibility of problem~\eqref{eq:find-ray-direction-eta-2} associated to problem~\eqref{eq:consensus-problem-distributed} is shown in \S~\ref{subsec:Consensus-problem-ell-infinity-norm-constraints}, \cf~Lemma~\ref{prop:consensus-problem-feasible-y-bar}.}}}. Then Algorithm~\ref{Alg:Dual-decompositoin-FGOR} outlines how the FGOR characteristic of the dual function $g$ is leveraged to solve problem \eqref{eq:consensus-problem-distributed} using dual decomposition.
\begin{algorithm}[t]
	\caption{Dual Decomposition with FGOR for Global Consensus} 
	\begin{algorithmic}[1]
	     \Require the solution $\tilde{\boldsymbol{\mu}}$ of problem~\eqref{eq:find-ray-direction-eta-2}, $c >0$, and $\beta>0$.  
	    \State $\boldsymbol{\lambda}^{(0)}=c\tilde{\boldsymbol{\mu}}$, $k=0$, $\boldsymbol{\lambda}_0^{(j)}=\boldsymbol{\lambda}_m^{(j)}=\mathbf{0}\in\R^n, \ j\in\mathbb{Z}^0_+$.
      \State $\forall i$, agent $i$ computes:  $$\mathbf{y}_i^{(k)}=\underset{\mathbf{y}_i\in\mathcal{Z}}{\arg\min}~f_i(\mathbf{y}_i)+ {\big(\boldsymbol{\lambda}^{(k)}_i-\boldsymbol{\lambda}^{(k)}_{i-1}\big)}\tran \mathbf{y}_i.  $$ 
\State $\forall i$, agent $i$ transmits CS one-bit information $b_i$, where
\begin{equation}\label{eq:Agent-communication-about-boundary}
b_i=\begin{cases}
    \displaystyle 0  \ ; & \ \text{$\mathbf{y}_i^{(k)}=\mathbf{\bar y}_i$}\\
         \displaystyle 1 \  ; & \ \text{Otherwise}.
    \end{cases}   
\end{equation} 
\State CS computes $\bar b=\max_{i\in\mathcal{S}} b_i$.
\State $\forall i$, CS broadcasts $\bar b$ to agent $i$. 

\If {$\bar b=0$,} \Comment{in RFGOR region} 
\State $\forall i$, agent $i$ and CS perform \eqref{eq:Lambda-Update} with 
\begin{equation} \nonumber
    \gamma_k=\gamma \ \mathrm{and} \ \mathbf{s}^{(k)}=\beta\boldsymbol{\tilde\mu}.
\end{equation}
\State Set $k:=k+1$ and go to step~2.
\Else  \Comment{RFGOR region is exited}
\State Go to step~12.
\EndIf 
\State $\forall i$, agent $i$ transmits $\mathbf{y}_i^{(k)}$ to CS. 
\State CS performs \eqref{eq:Lambda-Update} using $\gamma_k$ in Algorithm~\ref{Alg:Dual-Decomposition-Algorithm}, where
\begin{equation} \nonumber
    {\mathbf{s}}^{(k)}{=}[(\mathbf{y}^{(k)}_2{-}{\mathbf{y}}^{(k)}_1)\tran \ldots  (\mathbf{y}^{(k)}_{m}{-}\mathbf{y}^{(k)}_{m-1})\tran]\tran.
\end{equation}
\State Switch to Algorithm~\ref{Alg:Dual-Decomposition-Algorithm} with $\boldsymbol{\lambda}=\boldsymbol{\lambda}^{(k+1)}$.
	\end{algorithmic} 
	\label{Alg:Dual-decompositoin-FGOR}
\end{algorithm}
Note that Algorithm~\ref{Alg:Dual-decompositoin-FGOR} is a blend of Algorithm~\ref{Alg:Dual-Subgradient-Method} [\cf~\S~\ref{sec:Dual-Subgradient-Algorithm-The-Proposed- Approach}] and Algorithm~\ref{Alg:Dual-Decomposition-Algorithm} in the sense that we exploit the steps described in Algorithm~\ref{Alg:Dual-Subgradient-Method} to leverage FGOR characteristics on Algorithm~\ref{Alg:Dual-Decomposition-Algorithm}. In particular, steps~2-6 of Algorithm~\ref{Alg:Dual-decompositoin-FGOR} are corresponding to step~2 of Algorithm~\ref{Alg:Dual-Subgradient-Method}, where step~6 is the criterion we set in Algorithm~\ref{Alg:Dual-decompositoin-FGOR} to check whether $\boldsymbol\lambda^{(k)}$ is in an RFGOR region. Roughly speaking, if $\bar b=0$, it is straightforward to see that $\mathbf{y}^{(k)}=\mathbf{\bar y}$ [\cf~\eqref{eq:Agent-communication-about-boundary}], and thus $\boldsymbol\lambda^{(k)}$ is in an RFGOR region, \cf~Remark~\ref{rem:Equivalent-FGOR-for-Distribution}. Moreover, steps~3-4 and steps~5-11 of Algorithm~\ref{Alg:Dual-Subgradient-Method} are corresponding to steps~7-8 and steps~9-14 of Algorithm~\ref{Alg:Dual-decompositoin-FGOR}, respectively.

It is worth noting that each agent $i$ and the CS can individually solve \eqref{eq:find-ray-direction-eta-2}. If feasible, the solution $(\boldsymbol{\tilde\mu},\boldsymbol{\tilde\eta})$ obtained by solving \eqref{eq:find-ray-direction-eta-2} is the same for all agents and CS. Thus, both the CS and each agent $i$ know $\boldsymbol{\lambda}^{(0)}$ [\cf~\eqref{eq:initialization-point-subgrad-method}]. Moreover, they also know that the dual subgradient is $\mathbf{s}^{(k)}{=}\beta\boldsymbol{\tilde\mu}$ as long as $\boldsymbol{\lambda}^{(k)}$ resides in the respective RFGOR region, \cf~steps~6-7. This enables each agent $i$ to perform the dual variable update \eqref{eq:Lambda-Update} individually and in parallel until $\boldsymbol{\lambda}^{(k)}$ exits the RFGOR region, \cf~steps~6-9. However, the CS also performs the same computation to keep track of $\boldsymbol{\lambda}^{(k)}$ to use it after the algorithm steps exit the RFGOR region, \cf~steps~6-13.
We discuss the efficiency of Algorithm~\ref{Alg:Dual-decompositoin-FGOR} compared to Algorithm~\ref{Alg:Dual-Decomposition-Algorithm} in the sequel in detail.
\subsubsection{Fast Convergence with FGOR} Rather than choosing $\boldsymbol{\lambda}^{(0)}$ arbitrarily as in Algorithm~\ref{Alg:Dual-Decomposition-Algorithm}, Algorithm~\ref{Alg:Dual-decompositoin-FGOR} is initialized with $\boldsymbol{\lambda}^{(0)}$ [\cf~\eqref{eq:initialization-point-subgrad-method}] by exploiting the FGOR characteristic of the dual function $g$, \cf~problem~\eqref{eq:find-ray-direction-eta-2}. Unlike Algorithm~\ref{Alg:Dual-Decomposition-Algorithm}, $\boldsymbol{\lambda}^{(0)}$ in Algorithm~\ref{Alg:Dual-decompositoin-FGOR} permits one to use fixed stepsizes and makes the subgradient algorithm [\cf~\eqref{eq:Lambda-Update}] evolve along the specified ray yielding constant increments of $g$ as long as $\boldsymbol{\lambda}^{(k)}$ resides in an RFGOR region. Intuitively speaking, such a choice is fast in the sense that it directly drives $\boldsymbol{\lambda}^{(k)}$ toward $\boldsymbol{\lambda}^{\star}$ instead of potential detours where diminishing stepsizes such as $\gamma_0/k$, $\gamma_0/\sqrt{k}$ can be ineffective. %

\subsubsection{Efficient Communication with FGOR}\label{subsec:Efficient-Communication-using-Algorithm3} Note that Algorithm~\ref{Alg:Dual-Decomposition-Algorithm} requires each agent $i$ to transmit $n$ real numbers, i.e., $\mathbf{y}_i^{(k)}\in\R^n$ in step~5. In addition, it also requires CS to broadcast $n(m-1)$ real numbers to agents, \cf~step~8. Therefore, if $b$ bits are used to represent each real number, then Algorithm~\ref{Alg:Dual-Decomposition-Algorithm} requires $mnb+n(m-1)b=nb(2m-1)$ bits per epoch~\footnote{In the case of $64$-bit double precision floating point number format $b=64$.}. However, at the corresponding steps of Algorithm~\ref{Alg:Dual-decompositoin-FGOR} [\cf~step~3 and step~5], only \emph{one-bit} communication is needed and it does not depend on the length $n$ of $\mathbf{y}_i^{(k)}$. More specifically, in each epoch, Algorithm~\ref{Alg:Dual-decompositoin-FGOR} requires only $m+1$ bits, i.e., $m$ bits at step~3 and one bit at step~5, as long as $\boldsymbol{\lambda}^{(k)}$ resides in an RFGOR region. It is worth noting that engineering problems in real-world applications can consist of decision variables whose size $n$ can be in the order of several thousands, or even more \cite{Tibshirani_2009}. Thus, by leveraging FGOR characteristics, Algorithm~\ref{Alg:Dual-decompositoin-FGOR} can yield a considerable reduction in the communication overhead. 
\subsection{Problem~\eqref{eq:consensus-problem-distributed} with $\ell_{\infty}$-norm Constraints}\label{subsec:Consensus-problem-ell-infinity-norm-constraints}
In this section, we exploit specific structural properties of the constraint set $\mathcal{Z}$ in problem~\eqref{eq:consensus-problem-distributed} to provide a point $\mathbf{\bar y}\in\texttt{bnd}~\mathcal{Y}$, such that problem~\eqref{eq:find-ray-direction-eta-2} associated to~\eqref{eq:consensus-problem-distributed}, is feasible, where $\mathcal{Y}=\mathcal{Z}^m$. In particular, when $\mathcal{Z}$ conforms to $\ell_{\infty}$-norm constraints~\footnote{Typically, $\ell_{\infty}$-norm constraints appear in machine learning applications. For example, the $l_{\infty}$ regularization can be considered as a form of constraint on the $l_{\infty}$ norm of parameters.}, we show that problem~\eqref{eq:find-ray-direction-eta-2} associated to \eqref{eq:consensus-problem-distributed} is feasible for a specific $\mathbf{\bar y}\in\mathcal{Y}$. We first present the following remark, which will be useful in the sequel for our exposition.
\begin{remark}\label{Remark:normal-cone-of-Y-at-y-bar}
  Let $\mathcal{X}=\{\mathbf{x}\in\R^n \ | \ \|x\|_{\infty}\leq c\}$, where $c\in\R$. Moreover, let $\mathcal{V}$ be the set of all vertices of $\mathcal{X}$ and $\mathbf{\bar x}=[x_1 \ \ldots \ x_n]\tran\in\mathcal{V}$, where $\forall i=1,\ldots,n$, $x_i\in\R$. Then
  \begin{equation}\label{eq:normal-cone-at-x-bar}
  N_{\mathcal{X}}(\mathbf{\bar x})=\big\{\mathbf{p}=[p_1 \ \ldots \ p_n]\tran \ | \ \forall i \ p_i\in\R, \ p_i\lesseqgtr 0 \ \ \text{if} \ x_i\lesseqgtr 0 \big\}.   
  \end{equation}
\end{remark}
Then the following lemma substantiates the feasibility of problem~\eqref{eq:find-ray-direction-eta-2} associated to problem~\eqref{eq:consensus-problem-distributed}.
\begin{lemma}\label{prop:consensus-problem-feasible-y-bar}
Let Assumption~\ref{Assumption:PropStCvx} and Assumption~\ref{Assumption:Exclusion-of-Some-Regions} hold. Suppose $\mathcal{Z}=\{\mathbf{z}\in\R^n \ | \ \|\mathbf{z}\|_{\infty}\leq a\}$, where $a\in\R$. Let $\bar{\mathbf{y}}=[\bar{\mathbf{y}}_1\tran \ \ldots \ \bar{\mathbf{y}}_m\tran]\tran$, where $\forall i\in\mathcal{S}$,
\begin{equation}\label{eq:consensus-problem-feasibile-y-bar}
 \bar{\mathbf{y}}_i=
    \begin{cases}
    \displaystyle a\mathbf{1}_n \ \  \  \ \, ; & \ \text{$i$ odd}\\
         \displaystyle -a\mathbf{1}_n \ \ \, ; & \ \text{$i$ even}.
    \end{cases}       
\end{equation}
Then problem~\eqref{eq:find-ray-direction-eta-2} associated to problem~\eqref{eq:consensus-problem-distributed}, is feasible.    \end{lemma}
\begin{IEEEproof}
Clearly, $\mathbf{\bar y}\in\texttt{bnd}~\mathcal{Y}$, where $\mathcal{Y}=\mathcal{Z}^m$. By solving the equality constraints in \eqref{eq:find-ray-direction-eta-2}, we uniquely obtain the solutions $\boldsymbol{\tilde \mu}=(1/\beta)\mathbf{A}\mathbf{\bar y}$ and $\boldsymbol{\tilde \eta}=\mathbf{A}\tran\boldsymbol{\tilde \mu}$. Thus, we have
    \begin{equation}\label{eq:feasible-mu-for-consensus-problem}
     \boldsymbol{\tilde\mu}=(1/\beta)[2a\mathbf{1}_n\tran \ -2a\mathbf{1}_n\tran \ 2a\mathbf{1}_n\tran \ \ldots]\tran\in\R^{n(m-1)},   
    \end{equation}
    \begin{equation*}\label{eq:consensus-problem-feasibile-eta}
 \text{and }\boldsymbol{\tilde \eta}{=}
    \begin{cases}
    \displaystyle (1/\beta)[2a\mathbf{1}_n\tran \ -2a\mathbf{1}_n\tran]\tran \  ; & \hspace{-4cm} \text{$m=2$}\\
         \displaystyle (1/\beta)[2a\mathbf{1}_n\tran \ -4a\mathbf{1}_n\tran \ 4a\mathbf{1}_n\tran \ -4a\mathbf{1}_n\tran \ \ldots \ 2a\mathbf{1}_n\tran]\tran ; & \\ \ \hspace{4cm} \text{$m$ odd, $m\geq 3$}\\
         \displaystyle (1/\beta)[2a\mathbf{1}_n\tran \ -4a\mathbf{1}_n\tran \ 4a\mathbf{1}_n\tran \ -4a\mathbf{1}_n\tran \ \ldots \ -2a\mathbf{1}_n\tran]\tran  ; & \\ \ \hspace{4cm} \text{$m$ even, $m\geq 4$}.
    \end{cases}       
\end{equation*}
Clearly, $\boldsymbol{\tilde\eta}\in N_{\mathcal{Y}}(\mathbf{\bar y})$, \cf Remark~\ref{Remark:normal-cone-of-Y-at-y-bar}. %
Thus, problem~\eqref{eq:find-ray-direction-eta-2} associated to problem~\eqref{eq:consensus-problem-distributed}, is feasible.  
\end{IEEEproof}

\addb{Note that the feasibility of problem~\eqref{eq:find-ray-direction-eta-2} associated to the consensus problem \eqref{eq:consensus-problem-distributed}, confirms the satisfaction of Assumption~\ref{Assumption:dual-flat-assump} for \eqref{eq:consensus-problem-distributed}, \cf~Remark~\ref{Remark:Assumption-3}}. In the sequel, we numerically evaluate the performance of Algorithm~\ref{Alg:Dual-decompositoin-FGOR}. 

\section{Numerical Results}\label{sec:numerical-Results}
In this section, we numerically evaluate the benefits of FGOR for solving the consensus problem~\eqref{eq:consensus-problem-distributed}, \cf~\S~\ref{sec:Application:The-General-Consensus-Problem}. In particular, we compare the performance of Algorithm~\ref{Alg:Dual-decompositoin-FGOR} with Algorithm~\ref{Alg:Dual-Decomposition-Algorithm} \addbb{and splitting methods~\cite{condat_2013_splitting_methods,Boyd-Parikh-Chu-Peleato-Eckstein-2010,He_2021,Bai_2024}.} \addb{The experiments are conducted} for $\mathcal{Y}=\mathcal{Z}^m$ and $\mathcal{Z}$ is defined by an $\ell_{\infty}$-norm constraint, \cf~Lemma~\ref{prop:consensus-problem-feasible-y-bar}. For Algorithm~\ref{Alg:Dual-decompositoin-FGOR}, we choose $\boldsymbol{\lambda}^{(0)}=\addb{c}\boldsymbol{\tilde \mu}$, where $c\in\R$, for some $\beta\in\R$, \cf~\eqref{eq:feasible-mu-for-consensus-problem}. Furthermore, the constant stepsize $\gamma$ in Algorithm~\ref{Alg:Dual-decompositoin-FGOR} [\cf~step~7] is chosen to be sufficiently small to avoid undesired stepping over the region $\mathcal{V}$, \cf~\eqref{eq:define-region-V}. Note that the initialization point $\addb{c}\boldsymbol{\tilde\mu}$ of Algorithm~\ref{Alg:Dual-decompositoin-FGOR} is deterministic and is computed in advance, \cf~problem~\eqref{eq:find-ray-direction-eta-2}, \eqref{eq:feasible-mu-for-consensus-problem}. However, since the initialization\addb{s} for Algorithm~\ref{Alg:Dual-Decomposition-Algorithm} \addb{and \cite{condat_2013_splitting_methods,Boyd-Parikh-Chu-Peleato-Eckstein-2010}} are arbitrary, for a fair comparison, an average performance of the algorithm\addb{s} is considered. Specifically, the initialization point $\boldsymbol{\lambda}^{(0)}$ of Algorithm~\ref{Alg:Dual-Decomposition-Algorithm} and \addb{\cite{condat_2013_splitting_methods}} is chosen uniformly over the sphere centered at $\boldsymbol{\lambda}^{\star}$ with radius $\|\boldsymbol{\tilde\mu}-\boldsymbol{\lambda}^{\star}\|_2$. \addb{For \cite{Boyd-Parikh-Chu-Peleato-Eckstein-2010}, $\boldsymbol{\lambda}^{(0)}$ is chosen uniformly over the sphere centered at $\boldsymbol{\lambda}_{\texttt{admm}}^{\star}$ with radius $\|\boldsymbol{\tilde\mu}-\boldsymbol{\lambda}^{\star}\|_2$, where $\boldsymbol{\lambda}_{\texttt{admm}}^{\star}$ is the dual optimal point of the corresponding dual problem of \cite{Boyd-Parikh-Chu-Peleato-Eckstein-2010}.}

\addb{\subsection{Comparisons with Algorithm~\ref{Alg:Dual-Decomposition-Algorithm}}}
The comparisons \addb{with Algorithm~\ref{Alg:Dual-Decomposition-Algorithm}} are demonstrated with the following stepsize rules: 1) traditional stepsizes \addb{(TS)}: constant stepsize \addb{$\gamma_k=\addb{1/L_g}$, where $L_g$ represents the gradient Lipschitz constant of $g$}, polynomially decay stepsizes $\gamma_k=\gamma_0/k,\gamma_0/\sqrt{k}$, and geometrically decay stepsize $\gamma_k=\gamma_0 q^k$, where $\gamma_0>0$ is an appropriately chosen constant, 2) step decay stepsize \addb{(SD)}\cite{wang_2021}, 3) stepsize rule \cite{Yura_2020}. Moreover, the corresponding results are illustrated under two cases, \textsc{Case 1}: quadratic program formulation, and \textsc{Case 2}: regularized least squares regression with \addbb{two real datasets.}

\subsubsection{Case 1: Quadratic Program Formulation}\label{subsec:Case 1: Experiments with Quadratic Local Objectives}
We consider problem \eqref{eq:consensus-problem-distributed} with quadratic $f_i$s, i.e., $f_i(\mathbf{y}_i)=(1/m)(\mathbf{y}_i\tran \mathbf{A}_i \mathbf{y}_i+\mathbf{q}_i\tran \mathbf{y}_i),\,i\in\mathcal{S}$, where $\mathbf{A}_i\in\mathbb{S}^ {n}_{++}$ and $\mathbf{q}_i\in\R^n$ are arbitrarily chosen. %

\begin{figure}[!t]
\centering
\subfloat[]{\includegraphics[width=0.5\linewidth]%
{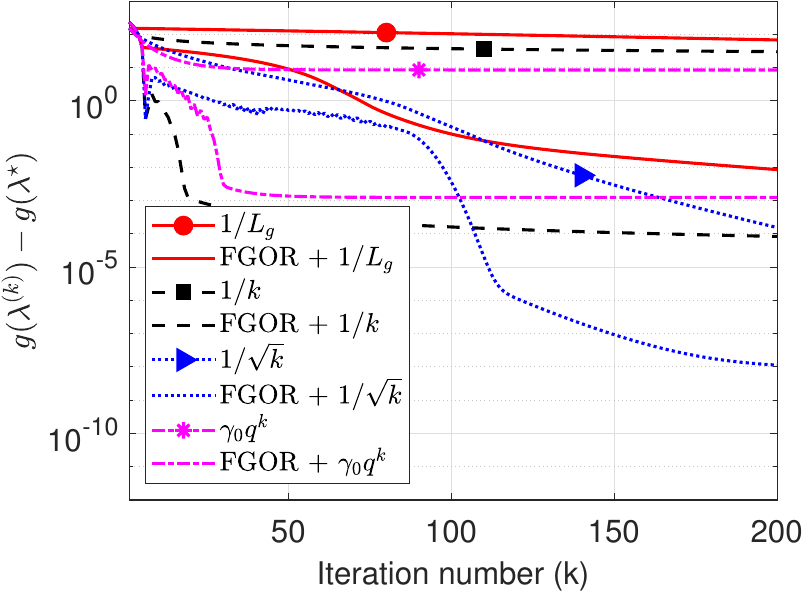}%
\label{subfig:Case1-9}}
\hfil
\subfloat[]{\includegraphics[width=0.5\linewidth]
{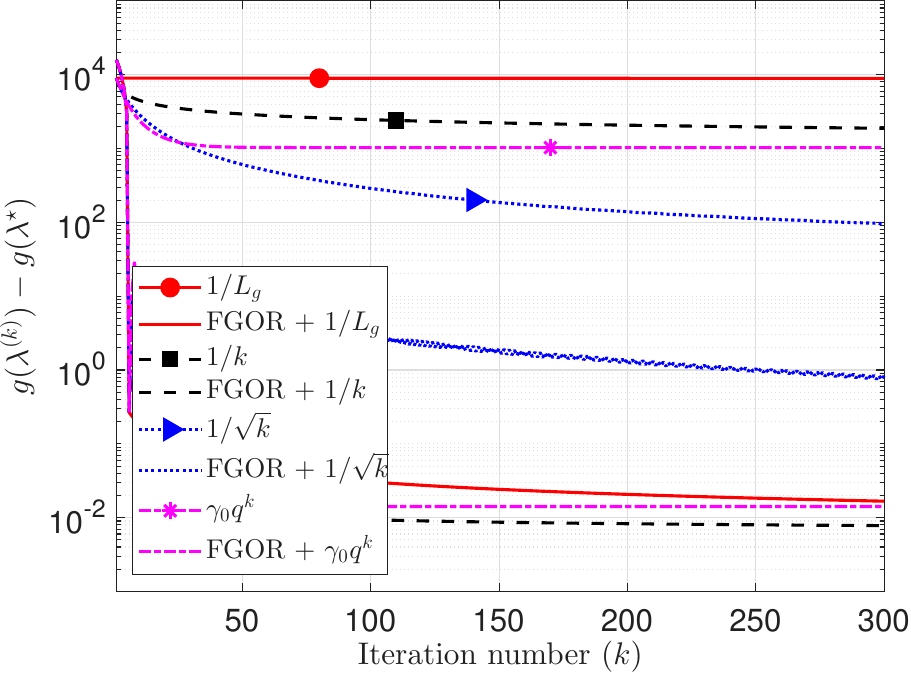}%
\label{subfig:Case1-10}}
\\
\subfloat[]{\includegraphics[width=0.5\linewidth]
{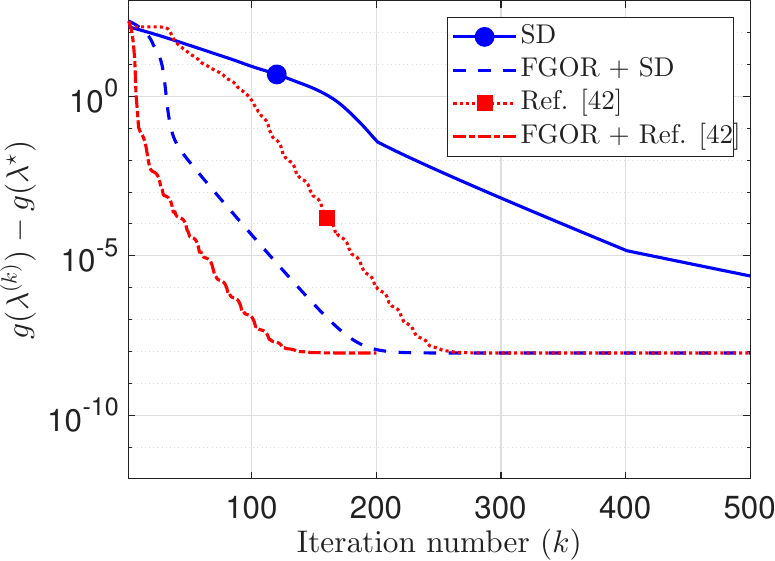}%
\label{subfig:Case1-11}}
\hfil
\subfloat[]{\includegraphics[width=0.5\linewidth]
{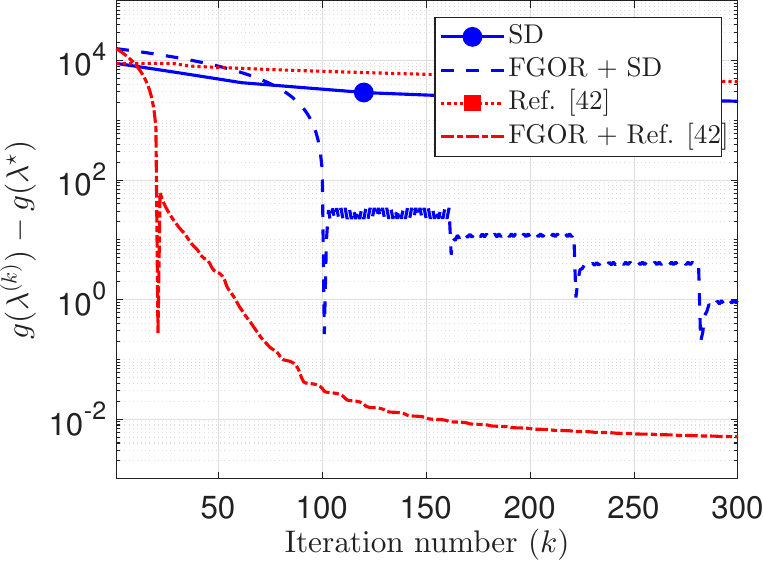}%
\label{subfig:Case1-12}}
\hfil
\caption{Comparisons with Algorithm~\ref{Alg:Dual-Decomposition-Algorithm}: \textsc{Case 1}: $\mathbf{y}^{\star}\in\texttt{int}~\mathcal{Y}$. \protect\subref{subfig:Case1-9} TS: $m=4$ and $n=10$. 
\addbb{\protect\subref{subfig:Case1-10} TS: $m=100$ and $n=10$.} \protect\subref{subfig:Case1-11} SD and Ref.~\cite{Yura_2020}: $m=4$ and $n=10$. \addbb{\protect\subref{subfig:Case1-12} SD and Ref.~\cite{Yura_2020}: $m=100$ and $n=10$.}}  
\label{Fig:Comparison-case1-2}
\end{figure}

\begin{figure}[!t]
\centering
\subfloat[]{\includegraphics[width=0.5\linewidth]{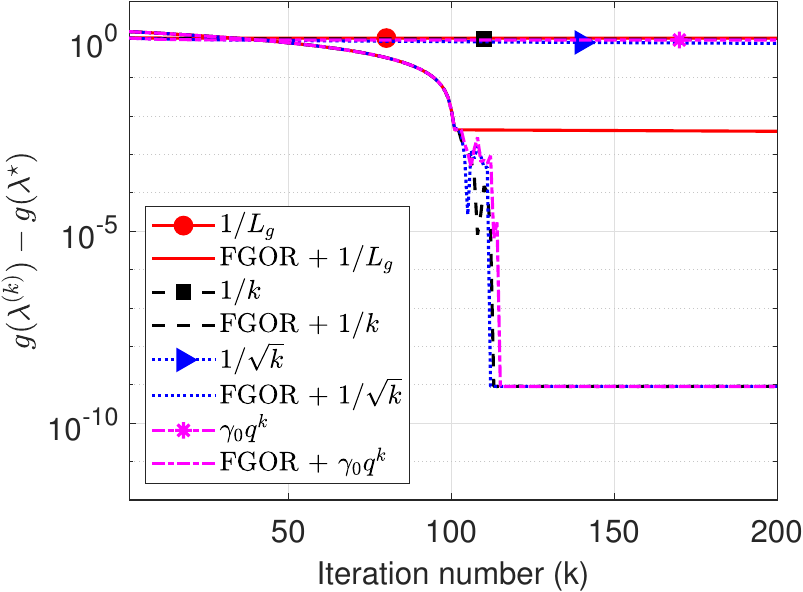}%
\label{subfig:Case1-5}}
\hfil
\subfloat[]{\includegraphics[width=0.5\linewidth]
{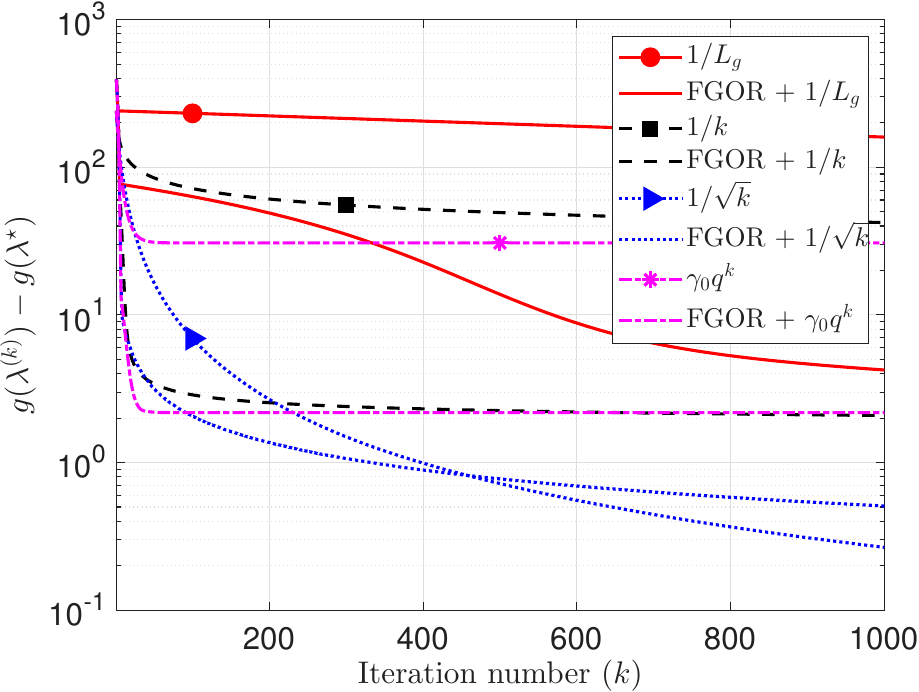}%
\label{subfig:Case1-6}}
\\
\subfloat[]{\includegraphics[width=0.49\linewidth]
{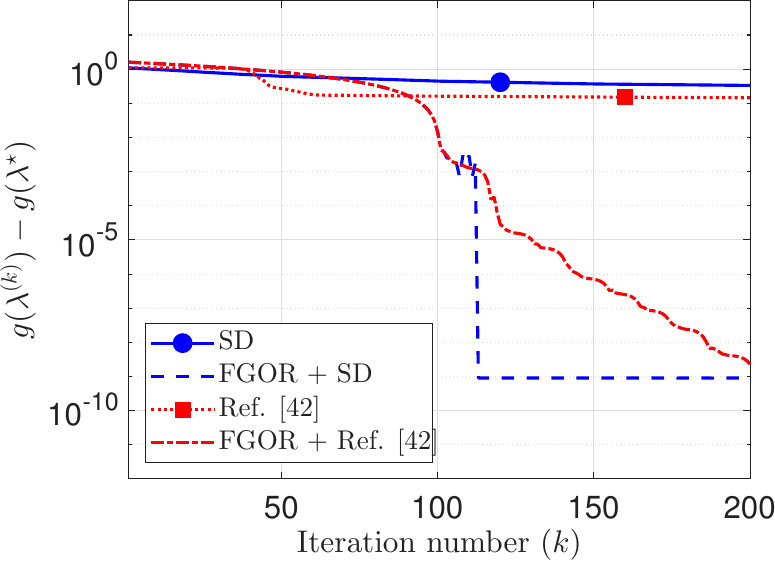}
\label{subfig:Case1-7}}
\hfil
\subfloat[]{\includegraphics[width=0.48\linewidth]
{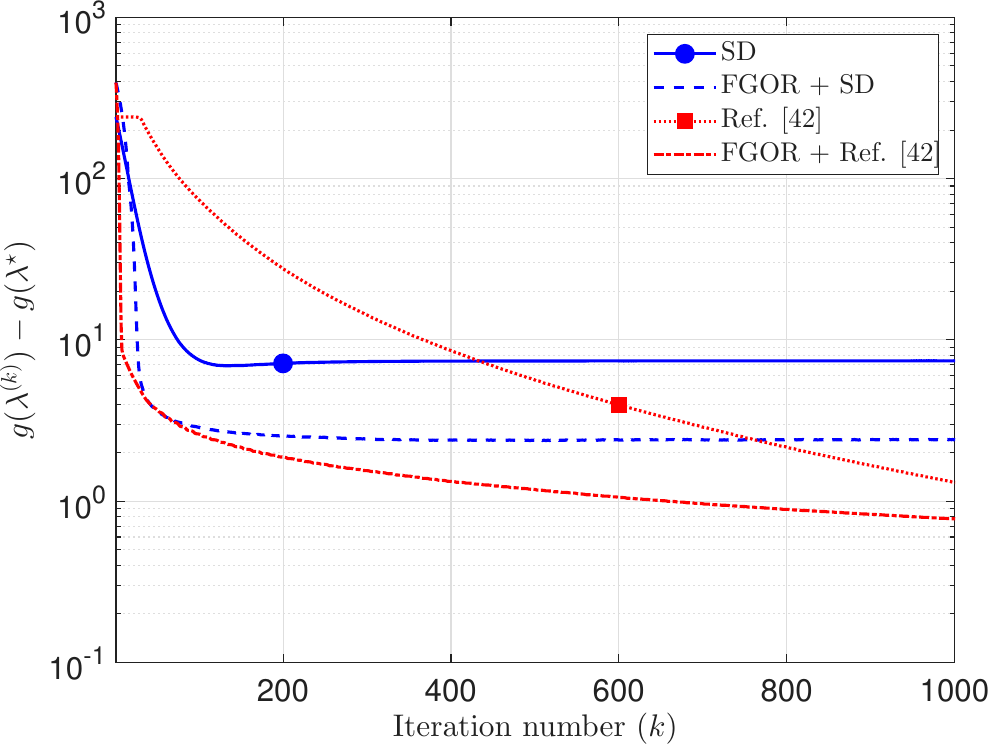}
\label{subfig:Case1-8}}
\hfil
\caption{Comparisons with Algorithm~\ref{Alg:Dual-Decomposition-Algorithm}: \textsc{Case 1}: $\mathbf{y}^{\star}\in\texttt{bnd}~\mathcal{Y}$. \protect\subref{subfig:Case1-5} TS: $m=5$ and $n=1$. 
\addbb{\protect\subref{subfig:Case1-6} TS: $m=100$ and $n=100$.} \protect\subref{subfig:Case1-7} SD and Ref.~\cite{Yura_2020}: $m=5$ and $n=1$. \addbb{\protect\subref{subfig:Case1-8} SD and Ref.~\cite{Yura_2020}: $m=100$ and $n=100$.}}  
\label{Fig:Comparison-case1-1}
\end{figure}

\figurename~\ref{Fig:Comparison-case1-2} and \figurename~\ref{Fig:Comparison-case1-1} show the comparison of results when $\mathbf{y}^{\star}\in\texttt{int}~\mathcal{Y}$ and $\mathbf{y}^{\star}\in\texttt{bnd}~\mathcal{Y}$, respectively, for different values of $m$ and $n$. Both figures show that Algorithm~\ref{Alg:Dual-decompositoin-FGOR} achieves significantly faster convergences than Algorithm~\ref{Alg:Dual-Decomposition-Algorithm}. To further clarify, the comparison of results in \figurename~\ref{Fig:Comparison-case1-2}\subref{subfig:Case1-11} is summarized in Table~\ref{Table:case1-comparison-table}. For each method, the number of iterations and the total number of bits required to achieve an accuracy of $10^{-5}$ are given by $k^{\star}$ and $b^{\star}$, respectively. Furthermore, for our method $b^{\star}=k_0(m+1)+n(k^{\star}-k_0)(2m-1) b$ and for the other methods $b^{\star}=nk^{\star}(2m-1)b$, where $k_0$ is the number of iterations that $\boldsymbol{\lambda}^{(k)}$ lies in the \addb{RFGOR} region and $b$ is the bits used to represent each component of communicated vectors $\mathbf{y}_i^{(k)}$ and $\boldsymbol{\lambda}_i^{(k)}$, \cf~\S~\ref{subsec:Efficient-Communication-using-Algorithm3}.  

\subsubsection{Case 2: Regularized Least Squares Regression}
We consider the following least squares regression problem with $\ell_{\infty}$-norm regularization: 
\begin{equation} \label{eq:least-squares-regression}
\begin{array}{ll}
\mbox{minimize} & f_0(\mathbf{y})=(1/m)\sum_{i{=}1}^{m}f_i(\mathbf{y}_{i})+d\|\mathbf{y}\|_{\infty}^2 \\
\mbox{subject to} & \mathbf{Ay}=\mathbf{0},
\end{array}
\end{equation}
where $ f_i(\mathbf{y}_i)=(1/N)\sum_{j=1}^N(\mathbf{y}_i\tran \mathbf{x}_j-t_j)^2$, $i\in\mathcal{S}$, $\mathbf{x}_j\in\R^n$ is the feature vector with $1$ in the first entry for the bias term, $t_j$ is the observation, and $N$ is the number of training examples of agent $i$. Moreover, $\mathbf{y}=[\mathbf{y}_{1}\tran \ \ldots \ \mathbf{y}_{m}\tran]\tran$, $\mathbf{A}$ is given in \eqref{eq:distributed-consenses-problem-matrix-A}, and $d>0$ is the regularization parameter. As we have already pointed out in \S~\ref{Sec:Flatness-Conjugate-Function} in connection with Assumption~\ref{Assumption:PropStCvx}, the problem of convex regularized minimization is equivalent to its corresponding regularization-constrained formulation, \addbb{\cf~Appendix~\ref{Appendix:Regularization-constrained-Formulation}.} Thus, problem~\eqref{eq:least-squares-regression} is equivalently reformulated~as:
\begin{equation} \label{eq:least-squares-regression-constrained}
\begin{array}{ll}
\mbox{minimize} & f_0(\mathbf{y})=(1/m)\sum_{i{=}1}^{m}f_i(\mathbf{y}_{i}) \\
\mbox{subject to} & \|\mathbf{y}\|_{\infty}^2\leq \bar{d} \\ 
& \mathbf{Ay}=\mathbf{0},
\end{array}
\end{equation}
where $\bar{d}>0$ is a hyperparameter, and experiments are conducted using the real estate valuation dataset (REV) \cite{Yeh-misc_real_estate_valuation} \addbb{and year prediction million song dataset (MSD)~\cite{year_prediction_msd_203}. The first dataset contains $276$ training data samples and $138$ testing data samples, where the length of the feature vector $\mathbf{x}_j$ is $n=7$. The latter dataset contains $463715$ training data samples and $51630$ testing data samples, where the length of the feature vector $\mathbf{x}_j$ is $n=91$. Both datasets were standardized before running the experiments. 

}

\addbb{\figurename~\ref{Fig:Comparison-linear-reg-2} illustrates the convergence of dual function values for two scenarios: $m=6$ with each agent holding $N=46$ samples for the REV dataset~[\cf~\figurename~\ref{Fig:Comparison-linear-reg-2}\subref{subfig:linear-regression-5} and \figurename~\ref{Fig:Comparison-linear-reg-2}\subref{subfig:linear-regression-7}], and $m=100$ with each agent holding $N=4637$ samples for the MSD dataset~[\cf~\figurename~\ref{Fig:Comparison-linear-reg-2}\subref{subfig:linear-regression-6} and \figurename~\ref{Fig:Comparison-linear-reg-2}\subref{subfig:linear-regression-8}]. As in \textsc{Case 1}, the results indicate that Algorithm~\ref{Alg:Dual-decompositoin-FGOR} converges to the optimal solution substantially faster than Algorithm~\ref{Alg:Dual-Decomposition-Algorithm}.}

\addb{\subsection{\addbb{Comparisons with Splitting Methods~\cite{condat_2013_splitting_methods,Boyd-Parikh-Chu-Peleato-Eckstein-2010,He_2021,Bai_2024}}}\label{subsec: Comparison-plitting-Methods} Let us now compare our proposed Algorithm~\ref{Alg:Dual-decompositoin-FGOR} with the state-of-the-art splitting methods~\cite[Algorithm 3.1]{condat_2013_splitting_methods}, \cite[\S~3.1]{Boyd-Parikh-Chu-Peleato-Eckstein-2010}, \addbb{\cite[\S~4.2]{He_2021}, and \cite[\S~3]{Bai_2024}}. The comparisons are illustrated when the local functions $f_i$s in \eqref{eq:consensus-problem-distributed} are quadratic, i.e., of the same form considered in \ref{subsec:Case 1: Experiments with Quadratic Local Objectives}. For the comparisons, Algorithm~\ref{Alg:Dual-decompositoin-FGOR} uses the following stepsizes as the switching stepsize rules: 1) constant stepsize $\gamma_k=1/L_g$, where $L_g$ is the gradient Lipschitz constant of $g$, 2) polynomially decay stepsize $\gamma_k=\gamma_0/k$, where $\gamma_0>0$ is an appropriately chosen constant, 3) stepsize rule~\cite{Yura_2020}. We note that the dual function formulations considered in \cite{Boyd-Parikh-Chu-Peleato-Eckstein-2010} and in our paper for the primal problem~\eqref{eq:consensus-problem-distributed} are different. Therefore, instead of comparing the convergence of dual function values, we consider the convergence of primal function values $f_0(\mathbf{y}_{\texttt{F}}^{(k)})$, where $\mathbf{y}_{\texttt{F}}^{(k)}=[{\mathbf{y}_{\texttt{f}}}^{(k)\mbox{\scriptsize T}}\,\ldots\,{\mathbf{y}_{\texttt{f}}}^{(k)\mbox{\scriptsize T}}]\tran$ and $\mathbf{y}_{\texttt{f}}$ is the average of $\mathbf{y}_i^{(k)}$s, $i=1,\ldots,m$. 

\addbb{\figurename~\ref{Fig:Comparison-case1-splitting-methods}\subref{subfig:Case1-splitting-1} and \figurename~\ref{Fig:Comparison-case1-splitting-methods}\subref{subfig:Case1-splitting-2} present a comparison of the results obtained using \cite[Algorithm 3.1]{condat_2013_splitting_methods} and \cite[\S~3.1]{Boyd-Parikh-Chu-Peleato-Eckstein-2010} for cases where $\mathbf{y}^{\star}\in\texttt{int}~\mathcal{Y}$ with $m=4$, $n=3$ and $\mathbf{y}^{\star}\in\texttt{bnd}~\mathcal{Y}$ with $m=5$, $n=1$, respectively.} Both figures show that Algorithm~\ref{Alg:Dual-decompositoin-FGOR} outperforms \cite[Algorithm 3.1]{condat_2013_splitting_methods} and \cite[\S~3.1]{Boyd-Parikh-Chu-Peleato-Eckstein-2010} on average. To further highlight the improvements in convergence and communication overhead, a summary of the results for \figurename~\ref{Fig:Comparison-case1-splitting-methods}\subref{subfig:Case1-splitting-2} is provided in Table~\ref{Table:case1-comparison-table}. \addbb{\figurename~\ref{Fig:Comparison-case1-splitting-methods}\subref{subfig:Case1-splitting-3} and \figurename~\ref{Fig:Comparison-case1-splitting-methods}\subref{subfig:Case1-splitting-4} present comparisons with benchmark methods, namely the balanced augmented Lagrangian method (ALM)~\cite[\S~4.2]{He_2021} and the double-penalty ALM~\cite[\S~3]{Bai_2024}. Both figures indicate that the proposed Algorithm~\ref{Alg:Dual-decompositoin-FGOR} outperforms the method in~\cite[\S~3]{Bai_2024}. Moreover, Algorithm~\ref{Alg:Dual-decompositoin-FGOR} demonstrates faster convergence during the initial iterations, however, as the iterations progress, the method in~\cite[\S~4.2]{He_2021} attains better overall convergence performance.


}

\begin{figure}[!t]
\centering
\subfloat[]{\includegraphics[width=0.5\linewidth]{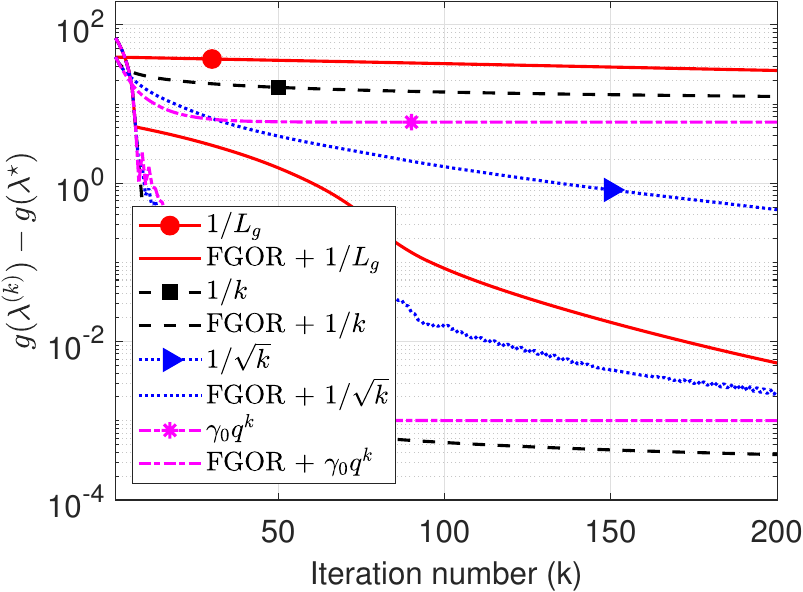}%
\label{subfig:linear-regression-5}}
 \hfil
 \subfloat[]{\includegraphics[width=0.5\linewidth]{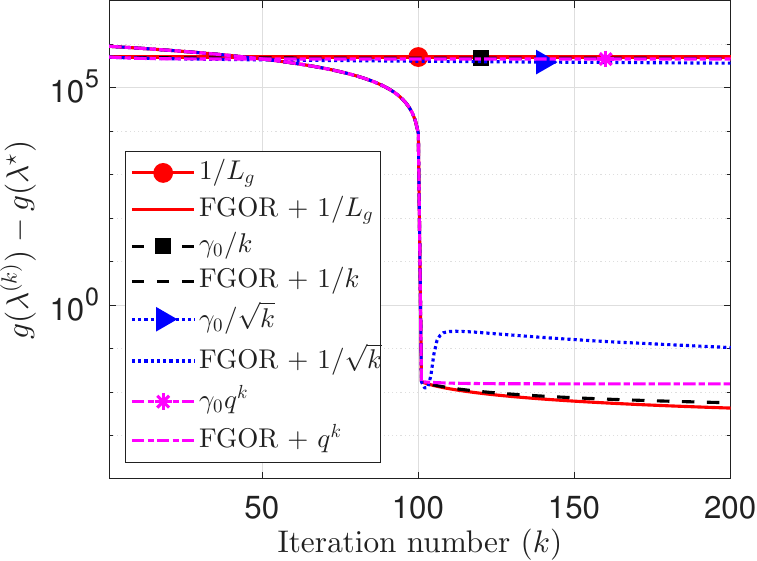}%
 \label{subfig:linear-regression-6}}
 \\
\subfloat[]{\includegraphics[width=0.5\linewidth]{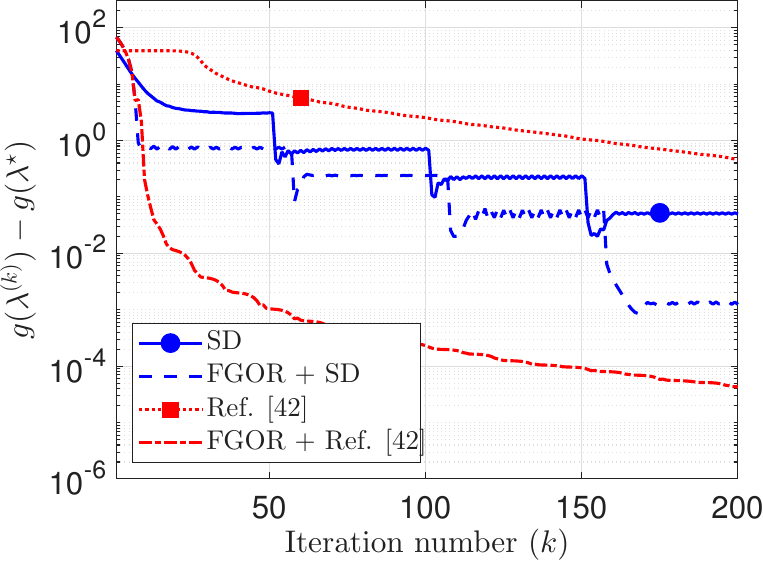}%
\label{subfig:linear-regression-7}}
 \hfil
 \subfloat[]{\includegraphics[width=0.5\linewidth]{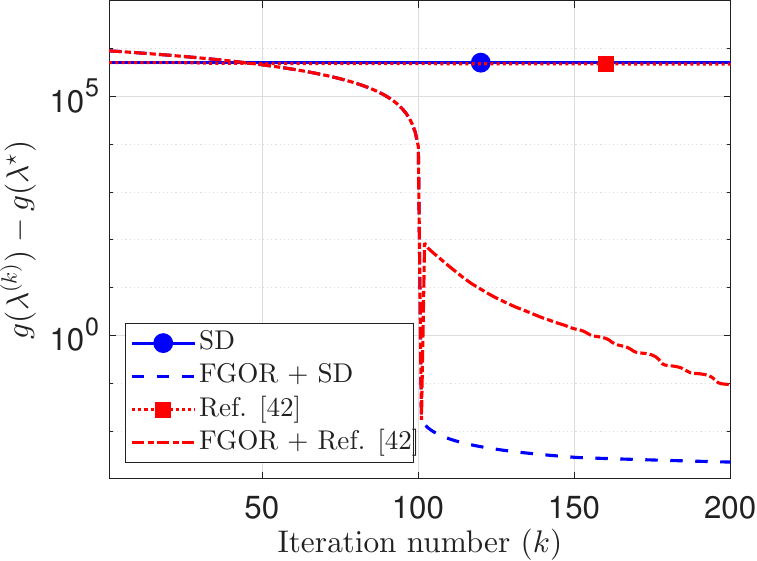}%
 \label{subfig:linear-regression-8}}
\caption{Comparisons with Algorithm~\ref{Alg:Dual-Decomposition-Algorithm}: \textsc{Case 2}. \protect\subref{subfig:linear-regression-5} TS: REV dataset with $m=6$.  \addbb{\protect\subref{subfig:linear-regression-6} TS: MSD with $m=100$. }\protect\subref{subfig:linear-regression-7} SD and Ref.~\cite{Yura_2020}: REV dataset with $m=6$. \addbb{\protect\subref{subfig:linear-regression-8} SD and Ref.~\cite{Yura_2020}: MSD with $m=100$.}
}   
\label{Fig:Comparison-linear-reg-2}
\end{figure}

%
}
%

\section{Conclusion}\label{sec:conclusion}
A novel characteristic of the conjugate function associated to a generic convex optimization problem was investigated, which is referred to as FGOR. Based on this characterization, we have also established the FGOR characteristic of the associated dual function. This characteristic, in turn, is used to devise a simple stepsize policy for dual subgradient methods that can be prepended with state-of-the-art stepsize rules, enabling faster convergences. %
Moreover, \addb{to highlight the practical implications of our theoretical results,} we have explored how the FGOR characteristics can be exploited when solving the global consensus problem using dual decomposition. More importantly, we have shown that FGOR can be leveraged to improve the performance of the dual decomposition methods not only in terms of the speed of convergence but also in terms of communication efficiency. \addbb{We have established the extension of FGOR to nonconvex formulations and demonstrated its effectiveness in stochastic optimization settings.} \addb{The numerical experiments highlight that FGOR can significantly improve the performance of existing stepsize methods and outperform state-of-the-art splitting methods on average in terms of both convergence speed and communication efficiency.} %
As such, our theoretical and empirical expositions suggested the importance of exploiting the structure of the problem, especially for large-scale problems that are common in almost all engineering application domains, such as signal processing and machine~learning.

\begin{appendices}
\section{Definitions} \label{Appendix:Definitions}
\begin{definition}[Conjugate Function]\label{Definition:Conjugate-function} 
Let $f:\R^n\to\overline{\R}$. Then the function $f^{*}:\R^n\to \overline{\R}$ defined by $f^{*}(\mathbf{y})=\sup_{\mathbf{x}\in\texttt{dom}~f}\left(\mathbf{y}\tran\mathbf{x}-f(\mathbf{x})\right)$ is called the conjugate of $f$.
\end{definition}

\begin{definition}[Indicator Function]\label{Definition-Indicator-function} 
Let $\mathcal{C}\subseteq \R^n$. Then the function $\delta_{\mathcal{C}}$ defined by
$\delta_{\mathcal{C}}(\mathbf{x})=0$ if $\mathbf{x}\in\mathcal{C}$, and $\delta_{\mathcal{C}}(\mathbf{x})=\infty$ if $\mathbf{x}\notin\mathcal{C}$ is called the indicator function of the set $\mathcal{C}$.
\end{definition}

\begin{definition}[Lipschitz Continuity] \label{definition:Lipschitz-continuity}
Let $f:\R^n\to \R$ with $\texttt{dom}~f=\mathcal{X}$. Then $f$ is Lipschitz continuous on $\mathcal{C}\subseteq \mathcal{X}$, if $\exists\, L\geq 0$, s.t $\|f(\mathbf{x})-f(\mathbf{y})\|_2\leq L\|\mathbf{x}-\mathbf{y}\|_2, \ \forall \, \mathbf{x},\mathbf{y}\in \mathcal{C}$, where $L$ is called the Lipschitz constant for $f$ on $\mathcal{C}$.
\end{definition}

\begin{figure}[!t]
\centering
\subfloat[]{\includegraphics[width=0.5\linewidth]%
{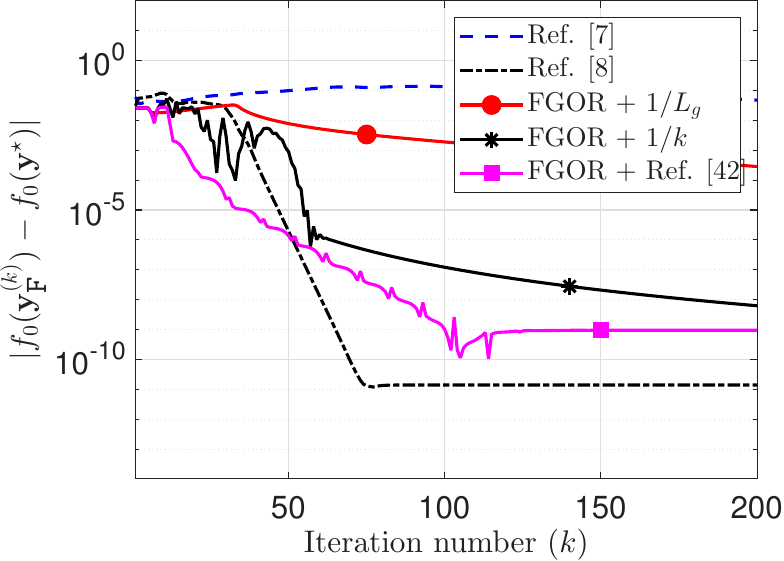}%
\label{subfig:Case1-splitting-1}}
\hfil
\subfloat[]{\includegraphics[width=0.5\linewidth]
{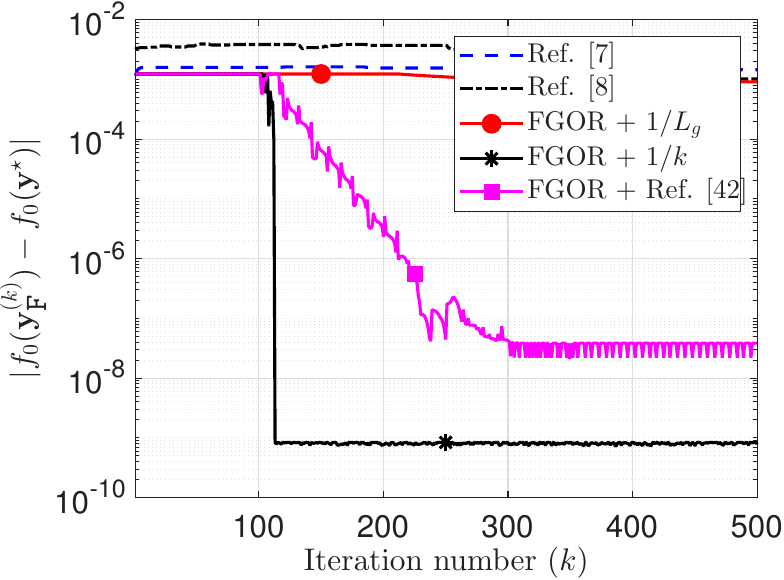}%
\label{subfig:Case1-splitting-2}}
\hfil
\subfloat[]{\includegraphics[width=0.5\linewidth]%
{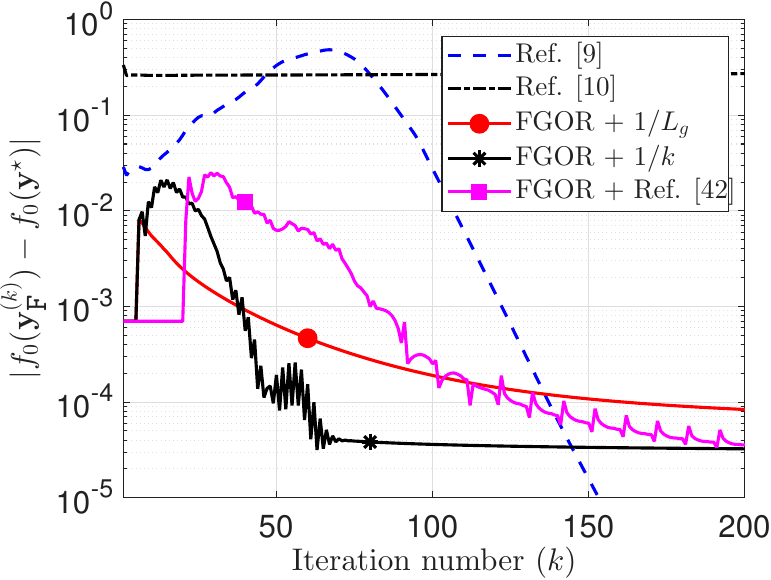}%
\label{subfig:Case1-splitting-3}}
\hfil
\subfloat[]{\includegraphics[width=0.5\linewidth]
{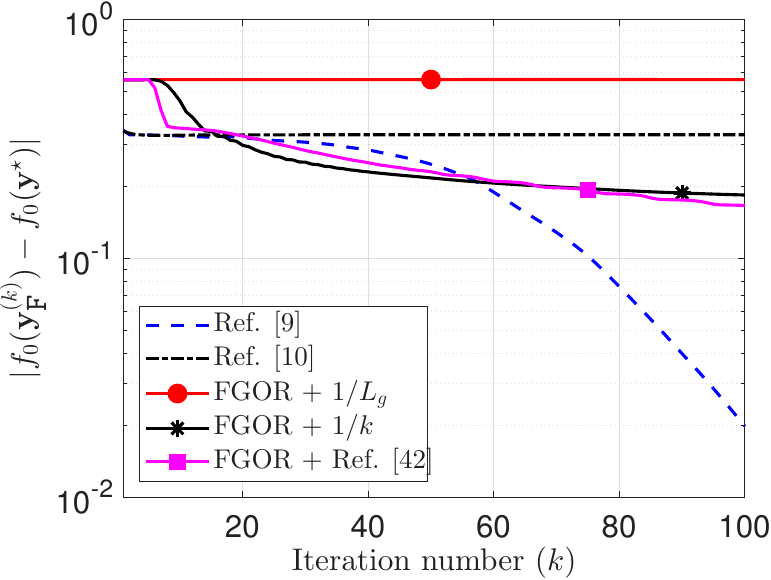}%
\label{subfig:Case1-splitting-4}}
\hfil
\caption{\addb{Comparisons with splitting methods.} \addb{\protect\subref{subfig:Case1-splitting-1} $\mathbf{y}^{\star}\in\texttt{int}~\mathcal{Y}$: $m=4$ and $n=3$.} 
\addb{\protect\subref{subfig:Case1-splitting-2} $\mathbf{y}^{\star}\in\texttt{bnd}~\mathcal{Y}$: $m=5$ and $n=1$. \addbb{\protect\subref{subfig:Case1-splitting-3} $\mathbf{y}^{\star}\in\texttt{int}~\mathcal{Y}$: $m=100$ and $n=10$.} \addbb{\protect\subref{subfig:Case1-splitting-4} $\mathbf{y}^{\star}\in\texttt{bnd}~\mathcal{Y}$: $m=100$ and $n=100$.}} 
}  
\label{Fig:Comparison-case1-splitting-methods}
\end{figure}

\begin{definition}[Strong Convexity] \label{definition:Strong-convexity}
Let $f:\R^n\to \R$ with $\texttt{dom}~f=\mathcal{X}$. Then $f$ is strongly convex on $\mathcal{C}\subseteq \mathcal{X}$, if $\exists\, l>0$, s.t $f\left(t\mathbf{x}{+}(1-t)\mathbf{y}\right)\leq tf(\mathbf{x}){+}(1-t)f(\mathbf{y}){-}\frac{1}{2}lt(1-t)\|\mathbf{x}-\mathbf{y}\|_2^2, \ \forall \, \mathbf{x},\mathbf{y}\in \mathcal{C}$, when $0< t< 1$, where $l$ is called the strong convexity constant for $f$ on $\mathcal{C}$.
\end{definition}
\begin{definition}[Relative Interior of a Convex Set] \label{Definition:relative-interior}
The relative interior of a convex set $\mathcal{X}$ relative to its affine hull is called the relative interior of $\mathcal{X}$, and is denoted by $\texttt{rint}~\mathcal{X}$. 
\end{definition}

\addbb{\section{Regularized versus Constrained Formulations} \label{Appendix:Regularization-constrained-Formulation}
The following remark substantiates the equivalence between the convex regularized minimization and its corresponding regularization-constrained formulation.
\begin{remark}
    Let $h:\R^{n}\rightarrow \R$ and $r:\R^{n}\rightarrow \R$ be proper, lower semicontinuous, and convex functions. Consider the regularized problem
    \begin{equation} \label{eq:regularized_problem}
\begin{array}{ll}
\mbox{minimize} & h(\mathbf{y})+ur(\mathbf{y}),
\end{array}
\end{equation}
and its corresponding regularization-constrained formulation
\begin{equation} \label{eq:regularization-constrained_problem}
\begin{array}{ll}
\mbox{minimize} & h(\mathbf{y})\\
\mbox{subject to} & r(\mathbf{y})\leq w 
\end{array}
\end{equation}
where $u\geq 0$ and $w>0$. Then $\mathbf{y}^\star$ is a solution of \eqref{eq:regularization-constrained_problem} for some $w>0$ if and only if it is a solution of \eqref{eq:regularized_problem} for some~$u\geq 0$.
\end{remark}
\begin{IEEEproof}
    The Lagrangian $L:\R^n\times\R\to\R$ associated with problem~\eqref{eq:regularization-constrained_problem} is given by
    \begin{equation}\label{eq:Lagrangian-regularization-constrained_problem}     L(\mathbf{y},\lambda)=h(\mathbf{y})+\lambda(r(\mathbf{y})-w),  
    \end{equation}
where $\lambda\geq 0$. Let $\mathbf{y}^\star$ be any optimal solution of \eqref{eq:regularization-constrained_problem} and $\lambda^\star$ be any optimal solution of the dual problem associated with \eqref{eq:regularization-constrained_problem} for some $w>0$. Then $\mathbf{y}^\star$ minimizes $L(\mathbf{y},\lambda^{\star})$, and thus, from \eqref{eq:Lagrangian-regularization-constrained_problem}, we have 
\begin{equation}\label{eq:minimum_regularized_problem}
\mathbf{y}^\star=\underset{\mathbf{y}\in\R^n}{\arg\min}~h(\mathbf{y})+\lambda^\star(r(\mathbf{y})-w).
\end{equation}
From \eqref{eq:minimum_regularized_problem}, we note that $\mathbf{y}^\star=\arg\min_{\mathbf{y}\in\R^n}~h(\mathbf{y})+\lambda^\star r(\mathbf{y})$. Thus, $\mathbf{y}^\star$ is an optimal solution of \eqref{eq:regularized_problem} when $u=\lambda^\star$. The reverse direction is established in a similar manner. 
\end{IEEEproof}

\section{Proof of Remark~\ref{rem:Equivalent-FGOR-for-Distribution}}\label{App:Equivalent-FGOR-for-Distribution}

Recall that $\mathbf{y}^{(k)}$ is the minimizer of the Lagrangian associated to problem~\eqref{eq:optimization_prob_main_extended} with Lagrange multiplier $\boldsymbol{\lambda}^{(k)}$. Then from \cite[Prop.~11.3]{Rockafellar-98}, we have
\begin{equation}\nonumber
 \mathbf{y}^{(k)}=\nabla f^*(\mathbf{A}\tran{\boldsymbol{\lambda}^{(k)}}){\iff}   \mathbf{A}\tran{\boldsymbol{\lambda}^{(k)}}\in \partial f_0\big({\mathbf{y}}^{(k)}\big) + \partial\delta_{\mathcal{Y}}\big({\mathbf{y}}^{(k)}\big).
\end{equation}
Replacing $\mathbf{y}^{(k)}$ by $\bar{\mathbf{y}}$ yields the result. \QE

}


\begin{table}[t]
\begin{center}
\begin{threeparttable}
\caption{Summary of Results in \figurename~\ref{Fig:Comparison-case1-2}\protect\subref{subfig:Case1-11} \addb{and \figurename~\ref{Fig:Comparison-case1-splitting-methods}\protect\subref{subfig:Case1-splitting-2}}} \label{Table:case1-comparison-table}
\vspace{-0.5cm}
          
\begin{tabular}{|m{3em}|c|c|c|c|}
 \hline
 Figure & Method & $k^{\star}$ & $k_0$ & $b^{\star}$\\ 
 \hline
\multirow{4}{*}{\figurename~\ref{Fig:Comparison-case1-2}\protect\subref{subfig:Case1-11}} & FGOR ${+}$ SD & $115$ & $21$ & $105+6580b$\\
  \cline{2-5}
& SD & $421$ & $NA$ & $29470b$\\
  \cline{2-5}
&   FGOR $+$ Ref~\cite{Yura_2020} & $61$ & $5$ & $25+3920b$\\
  \cline{2-5}
& Ref~\cite{Yura_2020} & $183$ & $NA$ & $12810b$\\
  \hline
 \multirow{3}{*}{\addb{\figurename~\ref{Fig:Comparison-case1-splitting-methods}\protect\subref{subfig:Case1-splitting-2}}} & \addb{FGOR $+$ Ref~\cite{Yura_2020}} & \addb{$184$}  & \addb{$100$} & \addb{$600+756b$} \\
 \cline{2-5}   
& \addb{Ref~\cite{condat_2013_splitting_methods}} & \addb{$>1000$}~\tnote{\addb{$a$}} & \addb{$NA$} & \addb{$>10000b$}\\
  \cline{2-5}
& \addb{Ref~\cite{Boyd-Parikh-Chu-Peleato-Eckstein-2010}} & \addb{$>1000$}~\tnote{\addb{$b$}} & \addb{$NA$} & \addb{$>5000b$} \\
  \hline
\end{tabular}
 \begin{tablenotes}
            \item[\addb{$a$}] \addb{For $k=1000$ iterations, $f_0(\mathbf{y}_{\texttt{F}}^{(k)})-f_0(\mathbf{y}^{\star})=0.0014$.}
            \item[\addb{$b$}] \addb{For $k=1000$ iterations, $f_0(\mathbf{y}_{\texttt{F}}^{(k)})-f_0(\mathbf{y}^{\star})=0.0021$.}
        \end{tablenotes}
    \end{threeparttable}
\end{center}
\end{table}
\


\addbb{\section{Extension of FGOR for Nonconvex Problems} \label{Sec:Appendix:Extensions-nonconvex-case}
In this section, we address the case in which the objective function $f_0$ in problem~\eqref{eq:optimization_prob_main} is nonconvex. Within this generalized setting, we first extend the FGOR properties established in Section~\ref{Sec:Flatness-Conjugate-Function} and Section~\ref{sec:Flatness-of-g} for the conjugate function $f^*$ and the dual function $g$, respectively. We then analyze the optimality conditions of problem~\eqref{eq:optimization_prob_main} and its corresponding dual problem~\eqref{eq:dual-problem}. Next, we consider the global consensus problem~\eqref{eq:consensus-problem-distributed} under additional structural properties. In particular, we examine the case where $f_0$ is composed of both convex and concave functions, leading to a nonconvex global objective function. Finally, empirical evaluations for this nonconvex setting are provided.



\subsection{Generalization of FGOR to Nonconvex Settings}\label{subsec:Generalization-of-FGOR-to-Nonconvex-Settings}

Recall that in Sections~\ref{Sec:Flatness-Conjugate-Function} and~\ref{sec:Flatness-of-g}, the FGOR properties of $f^*$ and $g$ are established with $f_0$ being strictly convex, which ensures the differentiability of both $f^*$ and $g$. Throughout this section, we consider $f_0$ to be nonconvex and $\mathcal{Y}$ to be convex. When $f_0$ is nonconvex, the differentiability of $f^*$ and $g$ can no longer be guaranteed. In such cases, we refer to the corresponding properties of $f^*$ and $g$ as \emph{fixed subgradient over rays} (FSGOR) instead of FGOR. Next, we present a couple of results useful for deriving the FSGOR characteristics of $f^*$ and $g$. 
\begin{remark}\label{Lemma:intV_mabs_subgrads_of_intY-nonconvex-case}
Let $ \mathcal{V}=\texttt{cl}\{\boldsymbol{{\nu}}\in\partial\texttt{con}(f_0)(\mathbf{y}) \ | \ \mathbf{y} \in \texttt{int}~\mathcal{Y}\}.$
Then 
$\exists~\mathbf{y}\in\texttt{int}~\mathcal{Y}~ \mbox{s.t.,}~\boldsymbol{\nu}_0\in\partial\texttt{con}(f)(\mathbf{y})\implies\boldsymbol{\nu}_0\in\mathcal{V}$.
\end{remark}
\begin{IEEEproof}
    The proof follows from the definition of $\mathcal{V}$.
\end{IEEEproof}
Note that Remark~\ref{Lemma:intV_mabs_subgrads_of_intY-nonconvex-case} is weaker than the corresponding convex counterpart Lemma~\ref{Lemma:intV_mabs_subgrads_of_intY}.
\begin{remark}\label{Remark:subdifferential-con(f)-at-boundary}
    Let $\mathbf{\bar y}\in\texttt{bnd}~\mathcal{Y}$. Then $\partial\texttt{con}(f)(\mathbf{\bar y})=\partial\texttt{con}(f_0)(\mathbf{\bar y})+N_{\mathcal{Y}}(\mathbf{\bar y})$.
\end{remark}
\begin{IEEEproof}
    We have that $\texttt{con}f=\texttt{con}f_0+\delta_{\mathcal{Y}}$. Then the convexity of $\delta_{\mathcal{Y}}$, together with Remark~\ref{Remark: subdifferential-indicator-function-at-boundary} yields the result.
\end{IEEEproof}

Let us next establish the FSGOR characteristic of $f^*$.
\begin{prop}
    Suppose $f_0$ is proper and lsc, the set $\mathcal{Y}$ is compact, and Assumption~\ref{Assumption:Exclusion-of-Some-Regions} holds. Then, for all $\boldsymbol{\nu}\in\R^n\setminus\mathcal{V}$, there exist $\boldsymbol{\eta}\in\R^n$ such that $\exists~\mathbf{\bar y}\in\texttt{bnd}~\mathcal{Y}$, $\forall$ $\alpha\geq 0$, $\boldsymbol{\nu}+\alpha\boldsymbol{\eta}\in\R^n\setminus\mathcal{V}$ and $\mathbf{\bar y}\in\partial f^*(\boldsymbol{\nu}+\alpha\boldsymbol{\eta})$.
\end{prop}
\begin{IEEEproof}
    Let $\boldsymbol{\nu}\in\R^n\setminus\mathcal{V}$. Then, from the contrapositive of the implication in Remark~\ref{Lemma:intV_mabs_subgrads_of_intY-nonconvex-case} and \cite[Prop.~11.3]{Rockafellar-98}, we have $\boldsymbol{\nu}\in\partial\texttt{con}(f)(\mathbf{\bar y})$ for some $\mathbf{\bar y}\in\texttt{bnd}~\mathcal{Y}$. Then $\boldsymbol{\nu}=\boldsymbol{\bar \nu}+\boldsymbol{\eta}$, where $\boldsymbol{\bar \nu}\in\partial\texttt{con}(f_0)(\mathbf{\bar y})$ and $\boldsymbol{\eta}\in N_{\mathcal{Y}}(\mathbf{\bar y})$, \cf~Remark~\ref{Remark:subdifferential-con(f)-at-boundary}. Moreover, we have $\forall$ $\alpha\geq 0$, $\boldsymbol{\nu}+\alpha\boldsymbol{\eta}=\boldsymbol{\bar \nu}+(1+\alpha)\boldsymbol{\eta}\in\partial\texttt{con}(f)(\mathbf{\bar y})$, because $\boldsymbol{\bar \nu}\in\partial\texttt{con}(f_0)(\mathbf{\bar y})$ and $(1+\alpha)\boldsymbol{\eta}\in N_{\mathcal{Y}}(\mathbf{\bar y})$, \cf~Remark~\ref{Remark:subdifferential-con(f)-at-boundary}. This, together with \cite[Prop.~11.3]{Rockafellar-98}, we have $\mathbf{\bar y}\in\partial(\texttt{con}(f))^*(\boldsymbol{\nu}+\alpha\boldsymbol{\eta})$. Then, since $(\texttt{con}(f))^*=f^*$ [\cf~\cite[Theorem 11.1]{Rockafellar-98}], $\mathbf{\bar y}\in\partial f^*(\boldsymbol{\nu}+\alpha\boldsymbol{\eta})$. 
\end{IEEEproof}

To derive FSGOR characteristic of the corresponding dual function $g$, we again rely on the fact that $g$ is based on a restriction of $-f^*$ to a linear space, \cf~\eqref{eq:dual-as-a-restriction-3}. We start by first imposing the following assumption.
\begin{assump}\label{Assumption:dual-flat-assump-noncvx}
   Assumption~\ref{Assumption:dual-flat-assump}, where the FGOR condition is replaced by the FSGOR condition. 
\end{assump}
As such, the FSGOR characteristic of $g$ is given by the following corollary.
\begin{corollary}
    Suppose $f_0$ is proper and lsc, the set $\mathcal{Y}$ is compact, Assumption~\ref{Assumption:Exclusion-of-Some-Regions} holds, and Assumption~\ref{Assumption:dual-flat-assump-noncvx} holds. Then $\exists~\boldsymbol{\lambda}\in\R^m$, $\exists~\boldsymbol{\mu}\in\R^m$ such that $\exists~\mathbf{a}\in\R^m$, $\forall$ $\alpha\geq 0$, $\mathbf{a}\in\partial g(\boldsymbol{\lambda}+\alpha\boldsymbol{\mu})$.
\end{corollary}
\begin{IEEEproof}
    The proof is similar to that of Corollary~\ref{Prop:flat-region-dual-function}, except that the subgradients of $f^*$ and $g$ are considered instead of their gradients.
\end{IEEEproof}

Incorporating the above generalized results, the subsequent analyses presented in Section~\ref{sec:Dual-Subgradient-Algorithm-The-Proposed- Approach} and Section~\ref{sec:Application:The-General-Consensus-Problem} continue to hold in the considered nonconvex setting, with the gradients of $f^*$ and $g$ replaced by their corresponding subgradients.

\subsection{Analysis of the Optimality of Problem~\eqref{eq:optimization_prob_main} and Problem~\eqref{eq:dual-problem}}

The dual function $g$ is always concave, even when the primal objective function $f_0$ [\cf~\eqref{eq:optimization_prob_main}] is nonconvex~\cite[\S~5]{Boyd-Vandenberghe-04}. Consequently, if a subgradient of $-g$ can be computed at each iteration of Algorithm~\ref{Alg:Dual-Subgradient-Method}, the algorithm converges to the dual optimal value $d^\star$ of the dual problem~\eqref{eq:dual-problem}. Nevertheless, the corresponding primal optimal value $p^\star$ is not necessarily attained unless strong duality holds~\cite[\S~5.2.3]{Boyd-Vandenberghe-04}. In general, it follows that $p^\star \neq d^\star$. Accordingly, it is straightforward to see that the proposed method guarantees convergence to the primal optimal solution, as summarized in the following remark.

\begin{remark}\label{Remark:Analysis-of-optimality}
Suppose that $f_0$ is proper and lsc, the set $\mathcal{Y}$ is compact, Assumption~\ref{Assumption:Exclusion-of-Some-Regions} holds, Assumption~\ref{Assumption:dual-flat-assump-noncvx} holds, and strong duality holds. Then, Algorithm~\ref{Alg:Dual-Subgradient-Method} yields the primal optimal solution $\mathbf{y}^\star$ of problem~\eqref{eq:optimization_prob_main}, provided that, at each iteration $k$, the Lagrangian 
$L(\mathbf{y},\boldsymbol{\lambda}^{(k)}) = f_0(\mathbf{y}) - \boldsymbol{\lambda}^{(k)\mbox{\scriptsize T}}(\mathbf{A}\mathbf{y} - \mathbf{b})$
associated with problem~\eqref{eq:optimization_prob_main_extended} is minimized optimally.
\end{remark}

It is important to emphasize that, in order to compute a subgradient of $-g$ at $\boldsymbol{\lambda}^{(k)}$, the Lagrangian must be minimized optimally at each iteration $k$ of Algorithm~\ref{Alg:Dual-Subgradient-Method}. Specifically, if $L(\mathbf{y},\boldsymbol{\lambda}^{(k)})$ is minimized and the corresponding $\mathbf{y}^{(k)}$ is obtained, then $\mathbf{A}\mathbf{y}^{(k)}-\mathbf{b}$ constitutes a subgradient of $-g$ at $\boldsymbol{\lambda}^{(k)}$. Since $f_0$ is generally nonconvex in $\mathbf{y}$, the Lagrangian $L$ is also nonconvex, which may, in general, impede exact Lagrangian minimization at each iteration. Nevertheless, in certain cases, the specific structure of the nonconvex functions can be exploited to enable optimal minimization of the Lagrangian. To illustrate this point, we consider the global consensus problem [\cf~\S~\ref{sec:Application:The-General-Consensus-Problem}, \eqref{eq:consensus-problem-distributed}], which is widely used across various application domains.

\subsection{Global Consensus Problem~\eqref{eq:consensus-problem-distributed} with Nonconvex $f_0$}

We consider a special case of \eqref{eq:consensus-problem-distributed} in which $m_1$ local functions are convex while the remaining $(m-m_1)$ local functions are concave. Consequently, the global objective function $f_0(\mathbf{y}_1,\ldots,\mathbf{y}_m) =(1/m)\sum_{i=1}^m f_i(\mathbf{y}_i)$ is generally nonconvex. In addition, we assume that the constraint set $\mathcal{Z}$ is defined by an $\ell_{\infty}$-norm bound, \cf~\S~\ref{subsec:Consensus-problem-ell-infinity-norm-constraints}.

To solve this problem, we employ Algorithm~\ref{Alg:Dual-decompositoin-FGOR}, a variant of Algorithm~\ref{Alg:Dual-Subgradient-Method} that is tailored for the global consensus problem form~\eqref{eq:consensus-problem-distributed}. At each iteration of Algorithm~\ref{Alg:Dual-decompositoin-FGOR}, every agent $i$ is required to optimally minimize its associated partial Lagrangian. For agents with convex local functions $f_i$, $i=1,\ldots,m_1$, the corresponding subproblems are convex, ensuring that this minimization can be performed optimally. For the remaining $(m-m_1)$ agents, the local functions $f_i$ are concave, and thus the associated partial Lagrangian minimization problems are nonconvex. Nevertheless, the minimization can still be performed optimally because the minimization of a concave function over a polyhedral set occurs at one of its vertices. Since the constraint set $\mathcal{Z}$ is an $\ell_{\infty}$-norm bound (i.e., a box) with explicitly known vertices, the Lagrangian minimization can be carried out exactly.

Consequently, all agents can solve their respective subproblems optimally, ensuring that Algorithm~\ref{Alg:Dual-decompositoin-FGOR} attains the dual optimal value of problem~\eqref{eq:dual-problem-consensus-problem}. Moreover, if strong duality holds, the algorithm also recovers the primal optimal solution of~\eqref{eq:consensus-problem-distributed}, \cf~Remark~\ref{Remark:Analysis-of-optimality}. To further illustrate these results, a couple of examples are provided in the following section.

%

\subsection{Empirical Evaluations for Nonconvex Settings} \label{subsec:FGOR-Empirical-Evaluations-for-Nonconvex-Setting}

In this section, we consider solving two nonconvex problems: the first yields strong duality, whereas the second exhibits weak duality. For both cases, we compare the performance of Algorithm~\ref{Alg:Dual-decompositoin-FGOR} with Algorithm~\ref{Alg:Dual-Decomposition-Algorithm}, where the latter is implemented using three stepsize rules: constant stepsize, polynomially decay stepsize $\gamma_0/k$, where $\gamma_0>0$ is an appropriately chosen constant, and stepsize rule in \cite{Yura_2020}. Specifically, the two nonconvex problems are as follows.
%
\begin{example}[A Problem with Strong Duality]
We consider the global consensus problem~\eqref{eq:consensus-problem-distributed} with $m=2$, where the local functions are given by $f_1(\mathbf{y}_1) = (1/16)\mathbf{y}_1\tran\mathbf{y}_1 + \mathbf{1}_n\tran \mathbf{y_1}+5$, $f_2(\mathbf{y}_2) = -\mathbf{y}_2\tran\mathbf{y}_2+2\mathbf{1}_n\tran\mathbf{y}_2+2$. Moreover, the constraint set $\mathcal{Z}$ is defined by an $\ell_\infty$-norm bound. \figurename~\ref{Fig:example-nonconvex-global-consensus-problem}\subref{subfig:example-nonconvex-1} shows that Algorithm~\ref{Alg:Dual-decompositoin-FGOR}, which exploits FSGOR characteristics, can achieve significantly faster convergence compared to Algorithm~\ref{Alg:Dual-Decomposition-Algorithm}.
\end{example}
%

%

\begin{figure}[!t]
\centering
\subfloat[]{\includegraphics[width=0.48\linewidth]%
{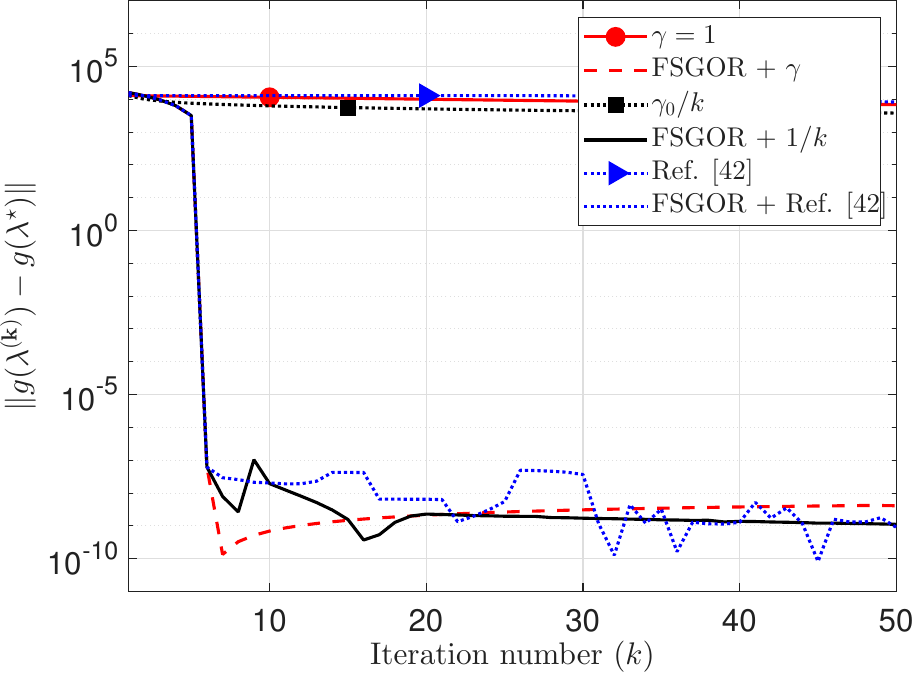}%
\label{subfig:example-nonconvex-1}}
\hfil
\subfloat[]{\includegraphics[width=0.48\linewidth]
{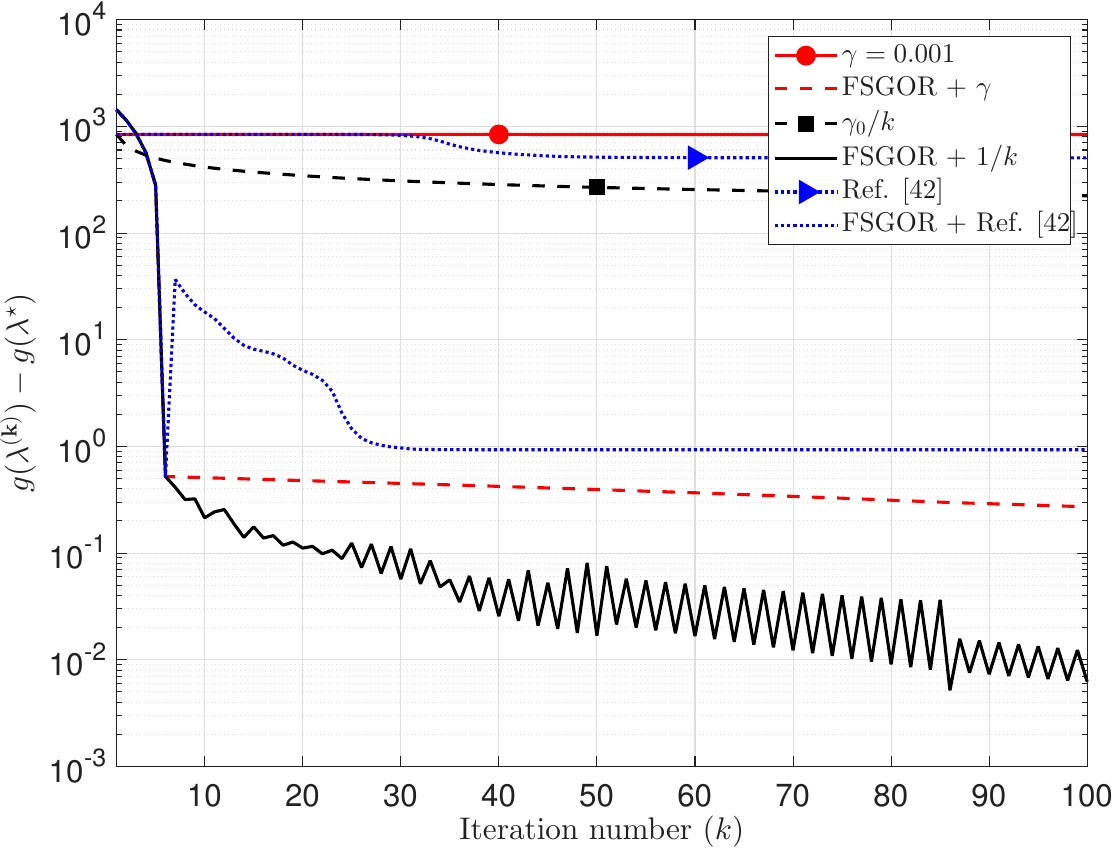}%
\label{subfig:example-nonconvex-2}}
\hfil
\addbb{\caption{Comparisons with Algorithm~\ref{Alg:Dual-Decomposition-Algorithm}: problem~\eqref{eq:consensus-problem-distributed} includes both convex and concave local objective functions. \protect\subref{subfig:example-nonconvex-1} Strong duality holds: $m=2,\, n=10$. 
\protect\subref{subfig:example-nonconvex-2} Strong duality does not hold: $m=10, \, n=10$. 
}  
\label{Fig:example-nonconvex-global-consensus-problem}}
\end{figure}
\begin{example}[A Problem with Weak Duality]
We consider solving the global consensus problem~\eqref{eq:consensus-problem-distributed} with $m=10$, where $5$ users have convex quadratic objectives and the remaining $5$ users have concave quadratic objectives, resulting in an overall nonconvex global objective. Moreover, we chose $n=10$. For each quadratic function, the Hessian matrices and the coefficients of the affine terms are generated randomly. The constraint set $\mathcal{Z}$ is again defined by an $\ell_\infty$-norm bound.

\figurename~\ref{Fig:example-nonconvex-global-consensus-problem}\subref{subfig:example-nonconvex-2} demonstrates that FSGOR enables fast convergence even in the presence of a nonzero duality gap between the primal and dual optimal values. However, due to weak duality, the primal optimal solution is not attained by any of the algorithms. In addition, although not shown in the figure, our numerical results indicate that the primal feasible points produced by the two algorithms are generally comparable, with neither algorithm consistently outperforming the other. 


\end{example}

}

\addbb{\section{FGOR on Stochastic Subgradient Methods}\label{subsec:FGOR-for-Stochastic-Subgradient-Method}

First we note that the stochastic subgradient method is used to solve the dual problem~\eqref{eq:dual-problem}. As the nondifferentiable case proceeds analogously [\cf~Appendix~\ref{subsec:Generalization-of-FGOR-to-Nonconvex-Settings}], the subsequent exposition focuses on the differentiable case. 

Let us start by considering a generic scenario in which we have a large dataset composed of sample pairs $(\mathbf{x}_j, t_j)$ ${j\in\{1,\ldots,N\}}$. 
The corresponding optimization problem is then expressed as
\begin{equation} \label{eq:optimization_prob_large_datasets}
\begin{array}{ll}
\mbox{minimize} & f(\mathbf{y}) = \frac{1}{N}\sum_{j=1}^N f(\mathbf{y}; (\mathbf{x}_j, t_j)) + \delta_\mathcal{Y}(\mathbf{y})\\
\mbox{subject to} & \mathbf{A}\mathbf{y} = \mathbf{b}.
\end{array}
\end{equation}

Following the terminology commonly adopted in the machine learning community~(\cf~\cite[\S~8.1.3]{goodfellow2016deep}), we refer to the subgradient method \eqref{eq:Lambda-Update} applied to the dual problem of~\eqref{eq:optimization_prob_large_datasets} as the \emph{batch} or \emph{deterministic} gradient method when the entire dataset $\{(\mathbf{x}_j, t_j)\}_{j\in\{1,\ldots,N\}}$ is utilized. More specifically, the gradient of the dual function $g$ at $\boldsymbol{\lambda}^{(k)}$ is given by 
\begin{equation}
   \nabla g(\boldsymbol{\lambda}^{(k)}) = -(\mathbf{A}\mathbf{y}(\boldsymbol{\lambda}^{(k)}) -\mathbf{b}), 
\end{equation}
where $\mathbf{y}(\boldsymbol{\lambda}^{(k)}) = \underset{\mathbf{y}\in\R^n}{\arg\min} \ f(\mathbf{y}) - \boldsymbol{\lambda}^{(k)\mbox{\scriptsize T}} \mathbf{A}\mathbf{y}$.

Next, consider the problem in which the entire dataset $\{(\mathbf{x}_j, t_j)\}_{j\in\{1,\ldots,N\}}$ is \emph{not} utilized. Specifically, in every iteration $k$, we select a batch of size $n_0 < N$ by uniformly sampling indices $j_1,\ldots,j_{n_0}$ from $\{1,\ldots,N\}$, either with or without replacement. We denote by $f^{(n_0,k)}$ and $g^{(n_0,k)}$ the corresponding primal objective function and the associated dual function, respectively. In this case, the subgradient method applied to the resulting dual problem is commonly referred to as a \emph{minibatch stochastic method}, and it is customary to simply refer to such approaches as \emph{stochastic methods}~\cf~\cite[\S~8.1.3]{goodfellow2016deep}. As such, we have 
\begin{equation}
   \nabla g^{(n_0,k)}(\boldsymbol{\lambda}^{(k)}) = -\left(\mathbf{A}\mathbf{y}^{(n_0,k)}(\boldsymbol{\lambda}^{(k)}) -\mathbf{b}\right), 
\end{equation}
where $\mathbf{y}^{(n_0,k)}(\boldsymbol{\lambda}^{(k)}) = \underset{\mathbf{y}\in\R^n}{\arg\min} \ f^{(n_0,k)}(\mathbf{y}) - \boldsymbol{\lambda}^{(k)\mbox{\scriptsize T}} \mathbf{A}\mathbf{y}$.

\begin{figure}[!t]
\centering
\subfloat[]{\includegraphics[width=0.5\linewidth]%
{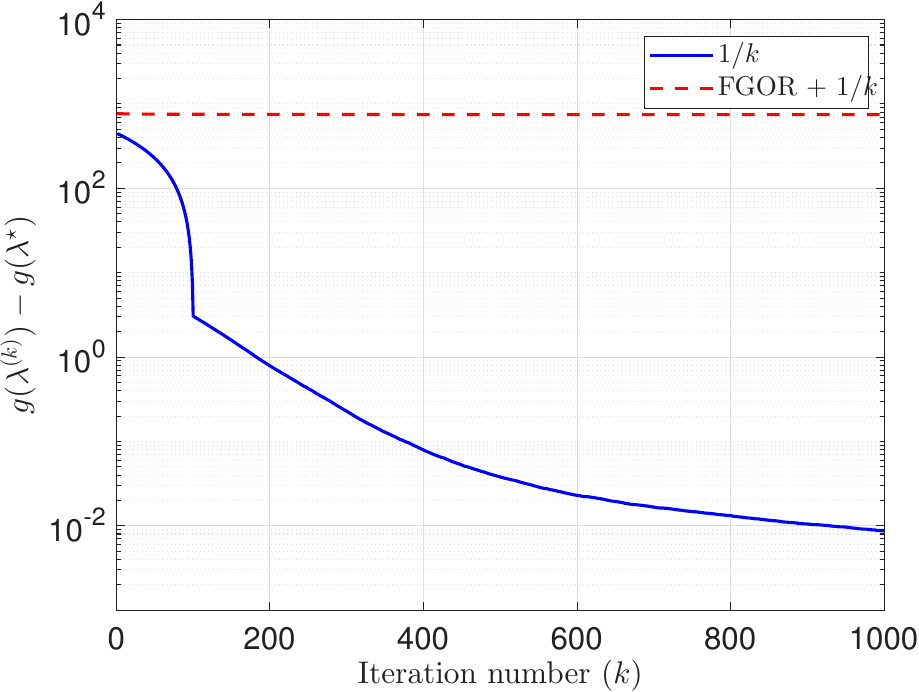}%
\label{subfig:stochastic-2}}
\hfil
\subfloat[]{\includegraphics[width=0.5\linewidth]
{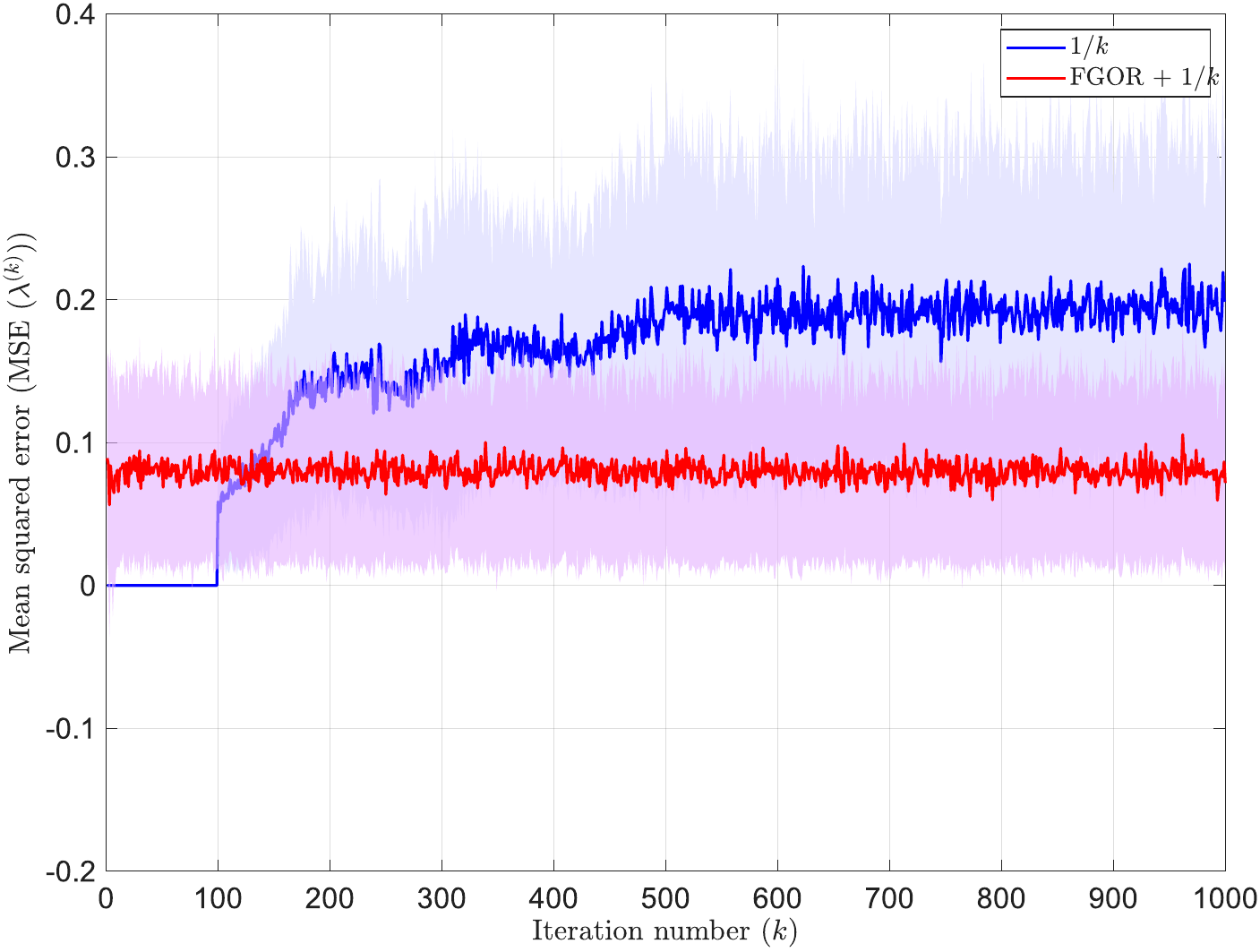}%
\label{subfig:stochastic-1}}
\hfil
\subfloat[]{\includegraphics[width=0.5\linewidth]%
{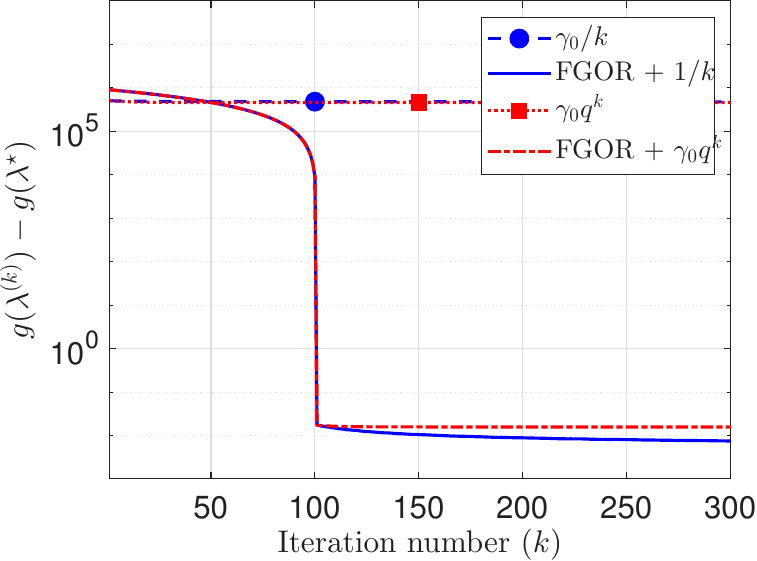}%
\label{subfig:stochastic-3}}
\hfil
\subfloat[]{\includegraphics[width=0.5\linewidth]
{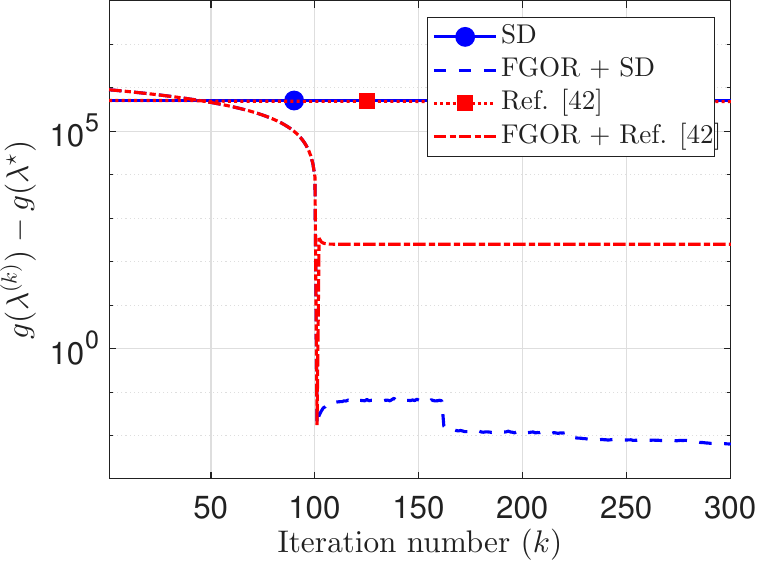}%
\label{subfig:stochastic-4}}
\hfil
\addbb{\caption{Comparisons using stochastic gradient method. \protect\subref{subfig:stochastic-2} REV dataset: $m=2$. \protect\subref{subfig:stochastic-1} REV dataset: $m=2$: evolution of squared error $||\mathbf{e}^{(n_0,k)}(\cdot)||^2_2$. \protect\subref{subfig:stochastic-3} TS: MSD: $m=100$. \protect\subref{subfig:stochastic-4} SD and Ref. \cite{Yura_2020}: MSD: $m=100$. 
}  
\label{Fig:Comparison-case1-Stochastic-methods }}
\end{figure}
It is worth pointing out that $\nabla g^{(n_0,k)}(\boldsymbol{\lambda}^{(k)})$ is a stochastic gradient for $g$ at $\boldsymbol{\lambda}^{(k)}$. As such, the error $\mathbf{e}^{(n_0,k)}(\boldsymbol{\lambda}^{(k)})$ given~by 
\begin{align}
  \mathbf{e}^{(n_0,k)}(\boldsymbol{\lambda}^{(k)}) &= \nabla g(\boldsymbol{\lambda}^{(k)})-\nabla g^{(n_0,k)}(\boldsymbol{\lambda}^{(k)}) \\
  & = \mathbf{A}\left(\mathbf{y}^{(n_0,k)}(\boldsymbol{\lambda}^{(k)}) -\mathbf{y}(\boldsymbol{\lambda}^{(k)})\right)
\end{align}
is generally a random vector~\footnote{\addbb{The general tractability of characterizing $\mathbf{e}^{(n_0)}(\boldsymbol{\lambda})$ is nontrivial and lies outside the main scope of this paper.}}.

We next gain empirical insights into the characteristics of the noise of the stochastic gradients. To this end, we first execute the stochastic gradient algorithm with a diminishing stepsize and the proposed stepsize, generating two sequences of dual variables, $\{\boldsymbol{\lambda}^{(k)}_{\textrm{dim}}\}_{k=0}^{1000}$ and $\{\boldsymbol{\lambda}^{(k)}_{\textrm{prop}}\}_{k=0}^{1000}$, respectively. Here we consider the regularized least squares regression with the REV dataset and $m=2$ with each agent holding $N=200$ samples and minibatch size $n_0=40$. The corresponding convergence plots are shown in \figurename~\ref{Fig:Comparison-case1-Stochastic-methods }\subref{subfig:stochastic-2}. Results suggests that the stochastic algorithm stays in a \emph{common}~\footnote{\addbb{Since the objective $f^{(n_0,k)}$ varies at each iteration $k$ due to the uniformly sampled indices $j_1, \ldots, j_{n_0}$, the term `common' is used to reflect this shared yet iteration-dependent structure. Technically speaking, such a common RFGOR region is to be characterized either probabilistically or statistically. However, such expositions are extraneous to the focus of this paper.}} RFGOR region for iteration indices $k \in \{0, \ldots, 100\}$. This behavior follows from the equivalence between \eqref{eq:find-ray-direction-eta---} and \eqref{eq:find-ray-direction-eta-2}, and from the fact that the solution to problem~\eqref{eq:find-ray-direction-eta-2} is independent of~$f$.

\figurename~\ref{Fig:Comparison-case1-Stochastic-methods }\subref{subfig:stochastic-1} shows the evolution of squared error $||\mathbf{e}^{(n_0,k)}(\cdot)||^2_2$ for the \emph{fixed} sequences $\{\boldsymbol{\lambda}^{(k)}_{\textrm{prop}}\}_{k=0}^{1000}$ and $\{\boldsymbol{\lambda}^{(k)}_{\textrm{dim}}\}_{k=0}^{1000}$ of dual variables generated earlier. In this respect, for all $k \in \{0, \ldots, 1000\}$ and for all $\boldsymbol{\lambda}^{(k)}_{\textrm{prop}}$, we randomly sampled minibatches of size $n_0 = 40$ over $100$ independent trials to generate the corresponding variance-shaded plot, illustrating both the mean performance and the variability across runs. The results indicate that as long as $\boldsymbol{\lambda}^{(k)}_{\textrm{prop}}$ lies within the common RFGOR region (e.g., for $k \in \{0,\ldots, 100\}$), the stochastic gradients of $g^{(n_0,k)}$ across different trials, when evaluated at these $\boldsymbol{\lambda}^{(k)}_{\textrm{prop}}$, \emph{coincide} with the true gradient of $g$. This is again a consequence of the equivalence between \eqref{eq:find-ray-direction-eta---} and \eqref{eq:find-ray-direction-eta-2}. On the other hand, when ${\boldsymbol{\lambda}}^{(k)}_{\textrm{prop}}$ lies outside the common RFGOR region, the results suggest that the stochastic gradients are influenced by different types of noise, as characterized in~\cite[\S~4.1.2]{Polyak_Intro_Opt}. For example, the observed errors resemble \emph{absolute deterministic noise}, i.e., for some $\epsilon > 0$, $\|\mathbf{e}^{(n_0,k)}_{\omega_j}(\cdot)\|_2^2 \leq \epsilon$, with $\omega_j$ denoting the index of the sample space corresponding to uniform sampling. Moreover, the results also indicate the presence of \emph{absolute random noise}, where for some $\sigma > 0$, $(1/100)\sum_{j=1}^{100}\|\mathbf{e}^{(n_0,k)}_{\omega_j}(\cdot)\|_2^2 \leq \sigma$. Similar observations are also obtained for the sequence $\{\boldsymbol{\lambda}^{(k)}_{\textrm{dim}}\}_{k=0}^{1000}$ generated using the diminishing stepsize rule. We further observe empirically that the convergence in \figurename~\ref{Fig:Comparison-case1-Stochastic-methods }\subref{subfig:stochastic-2} reaches a steady-state floor after a several thousands of iterations, which is expected since the noise component appears to be unbiased~\cite[\S~4]{Polyak_Intro_Opt}. 

Finally, \figurename~\ref{Fig:Comparison-case1-Stochastic-methods }\subref{subfig:stochastic-3} and 
\figurename~\ref{Fig:Comparison-case1-Stochastic-methods }\subref{subfig:stochastic-4} illustrate the convergence behavior of the stochastic methods with and without the proposed stepsize rule on the large dataset MSD. The results clearly demonstrate the advantage of exploiting the RFGOR region in accelerating convergence.
}
\section{Lipschitzian Properties of $f^*$ and $g$} \label{Appendix:Imp-Assump-Conj-Func}

We derive Lipschitzian properties of the conjugate function $f^*$ of $f$ [\cf~\eqref{eq:optimization_prob_main_extended}] and the dual function $g$ of problem \eqref{eq:optimization_prob_main_extended} under Assumption~\ref{Assumption:PropStCvx}. Since we consider strictly convex functions, our results are more general than the existing Lipschitzian properties for strongly convex objective functions~\cite[Prop.~12.60]{Rockafellar-98}. We first outline some results that are useful when asserting the intended results.
\begin{lemma}\label{Lemma:Lsc-Proper-Cvx-of-f-star}
Suppose Assumption~\ref{Assumption:PropStCvx} holds. Then the conjugate function $f^*$ of $f$ is lower semicontinuous (lsc), proper, and strictly convex with $\texttt{dom}~f^*=\R^n$.
\end{lemma}
\begin{IEEEproof}
The result follows directly from Theorem~11.1 of \cite{Rockafellar-98}, together with that $f^*(\boldsymbol{\nu})$ is finite for all $\boldsymbol{\nu}\in\R^n$.
\end{IEEEproof}

\begin{lemma}\label{Lemma:ConjugateDifferentiability}
Suppose Assumption~\ref{Assumption:PropStCvx} holds. Then the conjugate function $f^*$ of $f$ is differentiable on $\R^n$.
\end{lemma}
\begin{IEEEproof} Since Assumption~\ref{Assumption:PropStCvx} holds, the function $f$ is almost strictly convex, in the sense that $f$ is strictly convex in $\texttt{rint}~\mathcal{Y}$, \cf Definition~\ref{Definition:relative-interior}. Thus, Theorem~11.13 of \cite{Rockafellar-98} guarantees that $f^*$ is almost differentiable on $\texttt{int}~(\texttt{dom}~f^*)=\R^n$ [\cf Lemma~\ref{Lemma:Lsc-Proper-Cvx-of-f-star}], which in turn guarantees that $f^*$ is differentiable on $\R^n$. 
\end{IEEEproof}

\begin{lemma}\label{Lemma:Gradient-Boundness}
Suppose Assumption~\ref{Assumption:PropStCvx} holds. Then the gradient $\nabla f^*$ of the conjugate function $f^*$ is bounded. In particular,
\begin{equation} \label{eq:Upper-Bound-Grad-f-conj}
   \|\nabla f^*(\boldsymbol{\nu})\|_2\leq \max_{\mathbf{y}\in\mathcal{Y}} \|\mathbf{y}\|_2 \quad \forall~\boldsymbol{\nu} .
\end{equation}
\end{lemma}
\begin{IEEEproof}
	Proposition 11.3 of \cite{Rockafellar-98}, together with Assumption~\ref{Assumption:PropStCvx} and Lemma~\ref{Lemma:ConjugateDifferentiability} ensures that
	\begin{equation}\nonumber
	    \boldsymbol{\bar \nu}\in \partial f(\mathbf{\bar y}) \ \iff \ \mathbf{\bar y}\in\partial f^*(\boldsymbol{\bar \nu}) \ \iff \ \mathbf{\bar y}\in \{\nabla f^*(\boldsymbol{\bar \nu})\}.
	\end{equation}
Thus, $\mathbf{\bar y}=\nabla f^*(\boldsymbol{\bar \nu})$. The result follows immediately by using that $\mathbf{\bar y}\in\mathcal{Y}$.
\end{IEEEproof}

Next, the following proposition claims the Lipschitzian property of the conjugate function $f^*$ of $f$. 
\begin{prop}\label{Prop:Lipscitz-Continuity}
Suppose Assumption~\ref{Assumption:PropStCvx} holds. Then the gradient $\nabla f^*$ of $f^*$ is Lipschitz continuous with the constant $L=\sqrt{n}\max_{\mathbf{y}\in\mathcal{Y}}\|\mathbf{y}\|_2$.
\end{prop}
\begin{IEEEproof}
Define $\texttt{lip}~F( \boldsymbol{\nu})$, the \emph{Lipschitz modulus} of a single-valued mapping $F:\R^m\rightarrow \R^n$ at $\boldsymbol{\nu}$ \cite[Def. 9.1]{Rockafellar-98}: 
\begin{equation}\label{eq:Lip-defn}
 \texttt{lip}~F(\boldsymbol{\nu}) =  \underset{\substack{\boldsymbol{\nu}',\boldsymbol{\nu}''\rightarrow \boldsymbol{\nu} \\ \boldsymbol{\nu}'\neq \boldsymbol{\nu}''}}{\lim\sup} \quad \frac{\|F(\boldsymbol{\nu}')-F(\boldsymbol{\nu}'')\|_2}{\|\boldsymbol{\nu}'-\boldsymbol{\nu}''\|_2}.
\end{equation}
The central step in the proof is to show that $\forall{ \boldsymbol{\nu}\in\R^n}$, $\texttt{lip}~\nabla f^*(\boldsymbol{\nu})\leq L$, which in turn implies the desired result from \cite[Theorem 9.2]{Rockafellar-98}. To this end, for each vector
$\mathbf{u} = [u_1 \ u_2 \ \ldots \ u_n]\tran \in \R^n$, define the function $(\mathbf{u}\nabla f^*):\R^n\rightarrow \R$ by 
\begin{equation}\label{eq:u-Grad-f-conj-def}
    (\mathbf{u}\nabla f^*)(\boldsymbol{\nu})\triangleq \mathbf{u}\tran\nabla f^*(\boldsymbol{\nu}),
\end{equation}
which is the inner product between $\mathbf{u}$ and $\nabla f^*(\boldsymbol{\nu})$. Then,
\begingroup
\allowdisplaybreaks
\begin{align} \label{eq:Lip-1}
  \hspace{-3mm}  \sup_{\|\mathbf{u}\|_2=1}  \texttt{lip}~(\mathbf{u}\nabla f^*)(\boldsymbol{\nu}) &= \sup_{\|\mathbf{u}\|_2=1} \texttt{lip} \left(\sum_{i=1}^{n} u_i\frac{\partial f^*(\boldsymbol{\nu})}{\partial\nu_i}\right)\\  \label{eq:Lip-2}
   & \leq  \sup_{\|\mathbf{u}\|_2=1}  \sum_{i=1}^{n}  \texttt{lip}\left(
u_i\frac{\partial f^*(\boldsymbol{\nu})}{\partial\nu_i}\right)\\ \label{eq:Lip-3}
   & = \sup_{\|\mathbf{u}\|_2=1}  \sum_{i=1}^{n}  
|u_i|~\texttt{lip}~\frac{\partial f^*(\boldsymbol{\nu})}{\partial\nu_i}\\  \label{eq:Lip-4}
   &  \leq \sup_{\|\mathbf{u}\|_2=1}  \sum_{i=1}^{n}  
\left(|u_i|\max_{\mathbf{y}\in\mathcal{Y}}\|\mathbf{y}\|_2\right)\\ \label{eq:Lip-5}
    & {=} \left(\max_{\mathbf{y}\in\mathcal{Y}}\|\mathbf{y}\|_2\right)\left(\sup_{\|\mathbf{u}\|_2=1} \|\mathbf{u}\|_1\right)\\ \label{eq:Lip-6}
    &=\sqrt{n}\max_{\mathbf{y}\in\mathcal{Y}}\|\mathbf{y}\|_2.
\end{align}
\endgroup
The equality \eqref{eq:Lip-1} follows from \eqref{eq:u-Grad-f-conj-def}, \eqref{eq:Lip-2} follows from \cite[Exerc.~9.8(b)]{Rockafellar-98}, \eqref{eq:Lip-3} follows from \cite[Exerc.~9.8(a)]{Rockafellar-98}, and \eqref{eq:Lip-4} follows from Lemma~\ref{Lemma:Gradient-Boundness} and \eqref{eq:Lip-defn} since the scalar-valued function $\partial f^*(\boldsymbol{\nu})/\partial \nu_i$ is Lipschitz continuous with $\max_{\mathbf{y}\in\mathcal{Y}}\|\mathbf{y}\|$, \eqref{eq:Lip-5} follows from trivial rearrangements of terms, and \eqref{eq:Lip-6} follows directly from the fact that $\sup~\|\mathbf{u}\|_1$ is achieved when $\mathbf{u}=[1/\sqrt{n} \ \ldots \ 1/\sqrt{n}]\tran$. Finally, we note that $\texttt{lip}~\nabla f^*(\boldsymbol{\nu})= \sup_{\|\mathbf{u}\|_2=1}  \texttt{lip}~(\mathbf{u}\nabla f^*)(\boldsymbol{\nu})$ to conclude the result, \cf \cite[Exerc. 9.9]{Rockafellar-98}.
\end{IEEEproof}

It is worth noting that the  Lipschitzian properties claimed in Proposition~\ref{Prop:Lipscitz-Continuity} \emph{do not} rely on any \emph{strong convexity} properties of $f$. However, if such properties are imposed on $f$, Lipschitzian properties of $\nabla f^*$ directly follow from the results pertaining to the \emph{dualization of strong convexity}, \cf~\cite[Theorem~12.60]{Rockafellar-98}. The result is summarized in the following Remark.
\begin{remark}\label{Prop:Lipscitz-Continuity-more-strict}
Suppose the function $f_0$
is lsc, proper, and strongly convex with constant $\sigma$. Then the conjugate function $f^*$ of $f$ is differentiable and its gradient $\nabla f^*$ is Lipschitz continuous with constant $M=1/\sigma$.
\end{remark}
Note that the assertions of Proposition~\ref{Prop:Lipscitz-Continuity} are more general than those of Remark~\ref{Prop:Lipscitz-Continuity-more-strict}. This is illustrated by the following example.

\begin{example}[Limitations of Remark~\ref{Prop:Lipscitz-Continuity-more-strict}]\label{ex:y-power-4}
Let $h_0:\R\rightarrow \R$ and $l_0:\R\rightarrow \R$ are defined as $h_0(y) = y^4$ and $l_0(y)=1-\sqrt{1-y^2}$, respectively. Moreover, let $\mathcal{Y}=[-1,1]$. Note that both $h_0$ and $l_0$ and the set $\mathcal{Y}$ fulfill Assumption~\ref{Assumption:PropStCvx}. Let $h=h_0+\delta_{\mathcal{Y}}$ and $l=l_0+\delta_{\mathcal{Y}}$. Then from Proposition~\ref{Prop:Lipscitz-Continuity}, it follows that both $\nabla h^*$ and $\nabla l^*$ are Lipschitz continuous with the constant $1$. However, it is worth pointing out that $\partial^2 h_0(\bar y)/\partial y^2 = 0$ when $\bar y=0$, and $\partial^2 l_0(\bar y)/\partial y^2\rightarrow \infty$ as $\bar y\rightarrow \pm 1$. Thus, both $h_0$ and $l_0$ are \emph{not strongly convex}. As a result, Remark~\ref{Prop:Lipscitz-Continuity-more-strict} does not apply in this case. 
\end{example}
Following the Lipschitzian properties of $f^*$ [\cf~Proposition~\ref{Prop:Lipscitz-Continuity}], we establish the Lipschitzian properties of the dual function $g$ as given below.

\begin{corollary}\label{Lemma:UnconstrainedMin-of-f}
Suppose Assumption~\ref{Assumption:PropStCvx} holds. Then the dual function $g$ is differentiable on $\R^m$. The gradient $\nabla g$ of $g$ is Lipschitz continuous with the constant $G=||A||^2_2 L$.
\end{corollary}
\begin{IEEEproof} The identity~\eqref{eq:dual-grad-and-conjugate-grad} together with Lemma~\ref{Lemma:ConjugateDifferentiability} guarantees the differentiability of $g$. The same identity together with Proposition~\ref{Prop:Lipscitz-Continuity} guarantees the Lipschitz continuity of $\nabla g$.
\end{IEEEproof}
Note that it can be shown that $\nabla g$ is Lipschitz continuous with the constant $\|A\|^2_2 M$ when $f_0$ is strongly convex, \cf Remark~\ref{Prop:Lipscitz-Continuity-more-strict}. 
%


 \end{appendices}


\bibliographystyle{IEEEtran}
\bibliography{IEEEabrv, References}

\end{document}